\documentclass[12pt]{extarticle}
\usepackage{amsmath, amsthm, amssymb, mathtools, hyperref,color}
\usepackage{graphicx}
\usepackage{caption}
\usepackage{subcaption}
\usepackage{tikz}
\usepackage{graphicx}
\usepackage{caption}
\usepackage{xparse}
\usepackage{tikz}
\usepackage{tikz-cd}
\usepackage{tikz-cd}
\usepackage{tkz-graph}
\usetikzlibrary{shapes,arrows,positioning}
\usepackage{quiver}
\usepackage{algorithm}
\usepackage{algpseu docode}
\usepackage{lipsum}
\DeclareRobustCommand{\rchi}{{\mathpalette\irchi\relax}}
\newcommand{\irchi}[2]{\raisebox{\depth}{$#1\chi$}}



\newtheorem{theorem}{Theorem}
\numberwithin{theorem}{section}
\newtheorem{proposition}[theorem]{Proposition}
\newtheorem{lemma}[theorem]{Lemma}
\newtheorem{corollary}[theorem]{Corollary}
\newtheorem{definition}[theorem]{Definition}
\newtheorem{remark}[theorem]{Remark}

\newtheorem{example}[theorem]{Example}

\algnewcommand\algorithmicforeach{\textbf{for each}}
\algdef{S}[FOR]{ForEach}[1]{\algorithmicforeach\ #1\ \algorithmicdo}
\title{Rational Normal Curves, Chip Firing and Free Resolutions}
\author{Rahul Karki and Madhusudan Manjunath}

\begin{document}
\maketitle
\begin{abstract}
We study rational normal curves via a connection to the chip firing game. A key technique, introduced in this article, is  to interpret the defining ideal of the rational normal curve as an ideal associated to a generalisation of a cycle graph called a parcycle.  This association allows us to study rational normal curves by combinatorial methods. Given any Cohen-Macaulay initial monomial ideal of the rational normal curve, we  explicitly construct (via this association)  a corresponding Gr\"obner degeneration and an explicit combinatorial minimal free resolution of this Gr\"obner degeneration. 
We also construct an explicit minimal free resolution of  a certain determinantal ideal called the  toppling ideal associated to a parcycle,  a minimal cellular resolution for its initial ideal called the $G$-parking function ideal and an explicit minimal free resolution of the corresponding Gr\"obner degeneration.  Applications include minimal cellular resolutions for each Cohen-Macaulay initial monomial ideal of the rational normal curve,  explicit combinatorial formulas  for Hilbert series of certain lex-segment ideals and a combinatorial perspective on the Eagon-Northcott complex associated to the rational normal curve.

\end{abstract}
\section{Introduction}\label{intro}
Rational normal curves are among the simplest non-trivial examples of algebraic varieties, usually introduced in a first course in algebraic geometry \cite{harris2013algebraic}. The twisted cubic, an algebraic curve of degree three in $\mathbb{P}^3$, is a basic example of a non-(scheme-theoretic) complete intersection, i.e. the number of minimal generators of its defining ideal is strictly greater than its codimension. More generally, given an integer $n \geq 2$, the rational normal curve $\Gamma_n$ of degree $n$ is the image of  the map $\nu_n: \mathbb{P}^1 \rightarrow \mathbb{P}^{n}$ defined as $\nu_n([s:t])=[s^n:s^{n-1}t:s^{n-2}t^2:\ldots:t^n]$. Geometrically,  $\Gamma_n$ is the embedding of $\mathbb{P}^1$  into $\mathbb{P}^n$  given by the space of global sections of $\mathcal{O}_{\mathbb{P}^1}(n)$.  Furthermore, rational normal curves are images of the map induced by the complete canonical linear system on hyperelliptic curves. The defining ideal $P_n$ of the rational normal curve $\Gamma_n$ is minimally generated by the set of all $2 \times 2$ minors of the Hankel matrix 
\begin{center}
$
\begin{bmatrix}
x_1 & x_2 & x_3 &\dots &x_{n}\\
x_2 &x_3 & x_4 &\dots &x_{n+1}
\end{bmatrix}. 
$
\end{center}

The coordinate ring of $\Gamma_n$ is Cohen-Macaulay of depth and Krull dimension equal to two.

Minimal free resolutions offer a powerful tool for studying projective algebraic varieties, providing a bridge between the algebra and the geometry of the variety \cite{eisenbud2005geometry}. The rational normal curve being a determinantal variety (associated to the $2 \times 2$ minors of  a ``one-generic'' matrix of linear forms)  is minimally resolved by the Eagon-Northcott complex \cite{EAGONNORTHCOTT}.  In fact, the Eagon-Northcott complex resolves more general determinantal varieties \cite{EAGONNORTHCOTT}, \cite[Chapter 6]{eisenbud2005geometry}. One drawback of specialising such a general construction to rational normal curves is that their structure  beyond being determinantal varieties of one-generic matrices is not utilised.  Elucidating this aspect further, initial ideals of the rational normal curve and their minimal free resolutions are not directly evident from the Eagon-Northcott complex. 
 Initial ideals of the rational normal curve have found important applications in \cite{conca2007contracted}, \cite{deopurkar2014grobner}. These drawbacks are one source of motivation for the current work.  

Another source of motivation is the study of \emph{Gr\"obner degenerations}. Let $I$ be an ideal of $R_{n}=\mathbb{K}[x_1,\dots,x_n]$ for a field $\mathbb{K}$. Given an initial ideal $\mathcal{M}$ of $I$  with respect to a weight vector $\lambda$ on the variables, the Gr\"obner degeneration of $I$ to $\mathcal{M}$  with respect to $\lambda$ is a flat family over the affine line  obtained by homogenising $I$ with $\lambda$ as the weight (Definition \ref{GROBDEGDEF}, \cite[Section 15.8]{eisenbud2013commutative}). Gr\"obner degenerations enjoy important applications
to the study of Hilbert schemes and related objects called Betti strata \cite[Section 18.2]{miller2005combinatorial}. 

 In this article, we take a combinatorial approach to rational normal curves. More precisely, we consider combinatorial objects called \emph{parcycles} that are certain generalisations of cycles. 
 We study rational normal curves and their initial ideals in terms of certain ideals called \emph{toppling ideals and $G$-parking function ideals} associated to parcycles. Given any Cohen-Macaulay initial monomial ideal $\mathcal{M}$ of the rational normal curve $\Gamma_n$, we associate a parcycle $C_{\Pi,\mathcal{M}}$ to $\mathcal{M}$ and explicitly construct a minimal free resolution $\tilde{\mathcal{F}}_t(C_{\Pi,\mathcal{M}})$ for a
  Gr\"obner degeneration of $P_n$
 to $\mathcal{M}$,
 purely in terms of $C_{\Pi,\mathcal{M}}$. One of the main contributions of this article is the explicit construction of the minimal free resolution 
$\tilde{\mathcal{F}}_t(C_{\Pi,\mathcal{M}})$. In the following, we briefly outline its construction and refer to Section \ref{mainthmpf} for more details.

A parcycle is essentially a cycle along with a partition of the vertex set of this cycle such that the subgraph induced by each subset in the partition is connected. 
The building blocks of $\tilde{\mathcal{F}}_t(C_{\Pi,\mathcal{M}})$ are certain combinatorial objects called ``acyclic partitions'' that are derived from the parcycle  $C_{\Pi,\mathcal{M}}$. 
We first construct an intermediate complex $\mathcal{F}_t(C_{\Pi,\mathcal{M}})$ of free modules whose summands are labelled by (certain equivalence classes) of acyclic partitions  of $C_{\Pi,\mathcal{M}}$, we refer to Figure \ref{acyclic partitions} for an example and Subsection \ref{resolution of toppling ideal} for more details. The differentials of 
$\mathcal{F}_t(C_{\Pi,\mathcal{M}})$ arise from certain contractions of edges of acyclic partitions of $C_{\Pi,\mathcal{M}}$. The complex $\mathcal{F}_t(C_{\Pi,\mathcal{M}})$ also turns out to be  exact and minimal (Theorem \ref{minrehomtop}). The minimal free resolution $\tilde{\mathcal{F}}_t(C_{\Pi,\mathcal{M}})$ is obtained from $\mathcal{F}_t(C_{\Pi,\mathcal{M}})$ by applying a flat base change followed by quotienting with  a regular sequence. 
More precisely, we show the following:


\begin{theorem}\label{Grodegminres_theo}
Let $n \geq 2$ be any positive integer,  let $R_{n+1}=\mathbb{K}[x_1,\dots,x_{n+1}]$ where $\mathbb{K}$ is an algebraically closed field. Given any Cohen-Macaulay initial monomial ideal $\mathcal{M}$ of $P_n$, there exists a weight vector $\lambda$ in $\mathbb{N}^{n+1}$ such that $R_{n+1}/\mathcal{M}$ is the special fibre (fibre over $t=0$) of the Gr\"obner degeneration $R_{n+1}[t]/P_{n,t}$ of $R_{n+1}/P_n$ with respect to $\lambda$ and $\tilde{\mathcal{F}}_{t}(C_{\Pi,\mathcal{M}})$
is a minimal free resolution of  $R_{n+1}[t]/P_{n,t}$.
\end{theorem}
From our perspective, the chief utility of 
Theorem \ref{Grodegminres_theo} is that it offers a combinatorial perspective on
Gr\"obner degenerations of the rational normal curve (to Cohen-Macaulay initial monomial ideals). 
More precisely, the minimal free resolution $\tilde{\mathcal{F}}_{t}(C_{\Pi,\mathcal{M}})$ is constructed from a parcycle and allows the study of the corresponding Gr\"obner degeneration in terms of this underlying parcycle. For instance, we obtain combinatorial formulas for Hilbert series of Gr\"obner degenerations (Corollary \ref{weiHilGrob}). 

A few remarks on Theorem \ref{Grodegminres_theo} are in order.
\begin{enumerate}
\item\label{remthm1} By general principles \cite[Lemma 8.27]{miller2005combinatorial}, both $t-1$ and $t$
are  non-zero divisors on the Gr\"obner degeneration $R_{n+1}[t]/P_{n,t}$. This fact along  with specific  information about the  differentials shows that the specialisations  of  $\tilde{\mathcal{F}}_t(C_{\Pi,\mathcal{M}})$  at $t=1$ and $t=0$ are minimal free resolutions of $R_{n+1}/P_n$ and $R_{n+1}/\mathcal{M}$, respectively. Furthermore, the complex 
$\tilde{\mathcal{F}}_t(C_{\Pi,\mathcal{M}})$ can be viewed as an interpolation between these specialisations.

\item\label{rem2thm1} From Remark \ref{remthm1}, the quotient rings of $P_n$ and all its Cohen-Macaulay initial monomial ideals have the same Betti table. In fact, $R_{n+1}/P_n$ and $R_{n+1}/P_{n,t_0}$ have the same Betti table for all $t_0\in \mathbb{K}$. This can also be seen via a simple Artinian reduction argument.

\item Geometrically,  the Gr\"obner degeneration $R_{n+1}[t]/P_{n,t}$ (with the grading induced by $\lambda$) is the homogeneous  coordinate ring of an algebraic surface $\Gamma_{n,\mathcal{M}}$  in the  $\lambda$-weighted projective space over $\mathbb{K}$. Since the rational normal curve is a monomial (also known as ``toric'') curve, its Gr\"obner degenerations are all toric surfaces \cite{de2023semigroup}. Hence, $\Gamma_{n,\mathcal{M}}$ is a rational toric surface and is in fact a $\mathbb{P}^1$-fibration over $\mathbb{A}^1$. The complex $\tilde{\mathcal{F}}_t(C_{\Pi,\mathcal{M}})$ is a minimal free resolution of this embedded surface $\Gamma_{n,\mathcal{M}}$.  A complete classification of surfaces that arise as such Gr\"obner degenerations is an interesting topic for future work.  The combinatorial formulas for the Hilbert series of $\Gamma_{n,\mathcal{M}}$  may be useful in this context (Corollary \ref{weiHilGrob}).  
\end{enumerate}
The main fresh perspective in this article  is to view rational normal curves through the lens of a certain combinatorial game called the {\emph{chip firing game}}.

{\bf Chip Firing Game:} The chip firing game is a certain dynamical system associated to a graph $G$ (in our case, undirected, connected with possibly multiedges and with no loops).  The chip firing game consists of an assignment of an integer (to be thought of as the ``number of chips'') to each vertex of the graph called the \emph{initial configuration}. At each move of the chip firing game, some vertex $v$ can lend one chip  along every edge incident on it. This leads to a new configuration in which the vertex $v$ loses ${\rm val}(v)$ chips, where ${\rm val}(v)$ is its valence and every other vertex $u$ gains $a_{u,v}$ chips, where $a_{u,v}$ is the number of edges between $u$ and $v$.  A typical question is whether given an initial configuration whose total number of chips is non-negative,  there exists a finite sequence of moves that transforms it to a configuration where every vertex has a non-negative number of chips. The answer to this question is ``no'' in general  and the set of configurations which admit such a sequence depends subtly on the underlying graph, see \cite{BAKER2007766}, \cite{amini2010riemann} for more details. The chip firing game with its variants has been studied from several perspectives, for example statistical physics where it is usually referred to as the \emph{Abelian sandpile model}, combinatorics, algebraic and arithmetic geometry, and commutative algebra. We refer to the textbooks \cite{klivans2018mathematics}, \cite{corry2018divisors} for a more detailed discussion. 
 
This article mainly involves the commutative algebraic aspects of the chip firing game that we briefly describe in the following. This approach proceeds by associating certain ideals to $G$. The toppling ideal and the $G$-parking function ideal are two such ideals \cite{MR1896344,perkinson2013primer,manjunath2013monomials}. We briefly introduce them in the following.  We label the vertices of the graph $[1,\dots,n]$. We define the Laplacian matrix of $G$ to be $\Lambda_G:=\Delta-A_G$ where $\Delta$ is the diagonal matrix whose $i$-th diagonal entry is the valence of the vertex $i$ and $A_G$ is the vertex-vertex adjacency matrix of $G$. Consider the Laplacian lattice $L_G$ of $G$ \cite[Section 1]{amini2010riemann}, i.e. the sublattice of $\mathbb{Z}^{n}$ generated by the rows of the Laplacian matrix of $G$. By Kirchhoff's matrix-tree theorem, $L_G$ is a finite index sublattice of the root lattice $A_{n-1}=(1,\dots,1)^{\perp}\cap \mathbb{Z}^n$.  

Let $\mathbb{K}$ be an arbitrary field. Let $R_n=\mathbb{K}[x_1,\dots,x_n]$ be the polynomial ring with coefficients in $\mathbb{K}$, where the variable $x_i$ corresponds to the vertex $i$ of $G$.   The \emph{toppling ideal} $I_G \lhd R_n$, associated with $G$, is defined as the lattice ideal (\cite[Definition 7.2]{miller2005combinatorial}) of the Laplacian lattice $L_G$. In other words, $I_G= \langle {\bf x^u}-{\bf x^v}|~{\bf u}, ~{\bf v} \in \mathbb{N}^n, ~{\bf u}-{\bf v} \in L_G \rangle$ \footnote{Throughout the paper, we assume that $\mathbb{N}$ is the set of non-negative integers and hence contains $0$.}. Note that since $L_G$ is contained in $A_{n-1}$, the degrees of ${\bf x^u}$ and ${\bf x^v}$ are equal.  Hence, the ideal $I_G$ is $\mathbb{Z}$-graded  in the standard sense. The \emph{$G$-parking function ideal} $M_G \lhd  R_n$ is an initial monomial ideal of $I_G$ that closely reflects the properties of $I_G$.  For instance, $I_G$ and $M_G$ share their $\mathbb{Z}$-graded Betti table. Both $R_n/I_G$ and $R_n/M_G$ are Cohen-Macaulay of Krull dimension and depth one \cite[Proposition 10.1]{mohammadi2016divisors}.  However, $I_G$ and $M_G$ are not prime ideals except when $G$ is a tree, in which case both are prime ideals.  One notable feature of these ideals is that their minimal free resolutions can be described purely in terms of the underlying graph $G$, we refer to \cite[Theorem 1.1]{manjunathschwil} for the description.

The toppling ideal $I_{C_n}$ of the $n$-cycle $C_n$ is the determinantal ideal given by the maximal minors of the matrix \begin{center} $
\begin{bmatrix}
x_1 & x_2  & x_3 & \dots &x_{n-1} & x_{n}\\
x_2 &x_3 & x_4 &\dots & x_{n} & x_{1}
\end{bmatrix}. 
$
  \end{center}

The quotient rings of $I_{C_n}$ and $P_n$ (the defining ideal of the rational normal curve of degree $n$) have the same Betti table, see Table \ref{bTCn}.  

\begin{table}[h!]
\begin{center}
    \[\begin{array}{l|ccccccc}
    
       & 0& 1& 2&\dots&i&\dots & n-1\\ \hline
      \\0 &
      1 &0&0&\dots&0&\dots &0 \\ 
      \\
      1&0&\binom{n}{2}& 2\binom{n}{3}&\dots &i\binom{n}{i+1}&\dots&(n-1) \binom{n}{n}\\
     \end{array}.
    \]
\end{center}
\caption{Betti table of $I_{C_n}$}
\label{bTCn}
\end{table}

Furthermore, $I_{C_n}$ is an ideal in $R_{n}$ and $P_n$ is an ideal in $R_{n+1}$. The quotient ring $R_{n}/I_{C_n}$ is isomorphic to $R_{n+1}/(P_n + \langle x_{n+1}-x_1\rangle)$ and $x_1-x_{n+1}$ is a non-zero divisor of the quotient ring of $P_n$ (Proposition \ref{thm7}). Hence, their Betti tables are the same \cite[Proposition 1.1.5]{bruns1998cohen}. A natural question that arises in this context is whether there is a generalisation of the notion of toppling ideal such that (the defining ideal of) the rational normal curve is itself a toppling ideal in this generalised sense. 

In the current article, we take a first step towards such a generalisation. We study a generalisation of a cycle called a \emph{parcycle} and associate a toppling ideal to such a parcycle. 
A partition $\Pi=\{V_1,\dots,V_k\}$ of the vertex set of $C_n$ is called a \emph{connected $k$-partition} (or simply, a connected partition) if the subgraph induced by $V_i$ is connected for each $i$ from one to $k$.

\begin{figure}[h]
    \centering
    \[\begin{tikzcd}[sep=scriptsize]\label{figure 1}
	& {(1} & {6)} \\
	2 &&& 5 \\
	& {(3} & {4)}
	\arrow[no head, from=1-2, to=1-3]
	\arrow[no head, from=1-2, to=2-1]
	\arrow[no head, from=2-1, to=3-2]
	\arrow[no head, from=3-2, to=3-3]
	\arrow[no head, from=3-3, to=2-4]
	\arrow[no head, from=1-3, to=2-4]
\end{tikzcd}\]
    \caption{The parcycle $(C_6,\Pi)$ with $\Pi =\{\{1,6\},\{2\},\{3,4\},\{5\}\}$}
    \label{aparc6}
\end{figure}
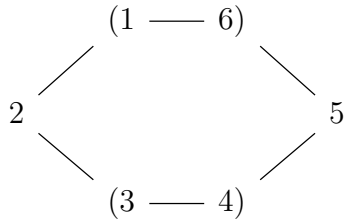

A parcycle is an ordered pair $(C_n,\Pi)$ where $C_n$  is a cycle of length $n$ and the set $\Pi=\{V_1,\dots,V_k\}$ is a connected partition.  Informally, the parcycle $(C_n,\Pi)$  is the cycle $C_n$ with its vertex set partitioned according to $\Pi$, we refer to Figure \ref{aparc6} 
for an example.
The set of \emph{relevant edges} of the parcycle is defined as the set of all edges $\{u,w\}$ of $C_n$, such that $u \in V_i$ and $w \in V_j$ for $V_i,~V_j \in \Pi$ and $i \neq j$. For example, $\{1,2\},\{2,3\},\{4,5\},\{5,6\}$ (in Figure \ref{aparc6}) are the relevant edges of $(C_6,\Pi)$.

In the following, we define the notion of toppling ideal of a parcycle. Given a connected $k$-partition $\Pi=\{V_1,\dots,V_k\}$ of the vertex set of $C_n$, let $C_{\Pi,n}$ be the parcycle $(C_n,\Pi)$.  Let $b_j$ be the row of the Laplacian matrix of $C_n$ corresponding to the vertex $j$. We define the \emph{Laplacian lattice} $L_{\Pi}$ of the parcycle $C_{\Pi,n}$ as the sublattice of the Laplacian lattice of $C_n$ generated by the vectors $b_{V_1},\dots,b_{V_k}$, where $b_{V_i}=\sum_{j \in V_i}b_j$. We define the \emph{toppling ideal} $I_\Pi$ of $C_{\Pi,n}$ as the lattice ideal of $L_{\Pi}$.   In contrast to the case of graphs, the Laplacian lattice of the parcycle $C_{\Pi,n}$ is not necessarily a finite index sublattice of $A_{n-1}$. More precisely, it has rank $k-1$, where $k$ is the size of $\Pi$.
Hence, the quotient ring $R_n/{I_{\Pi}}$ has Krull dimension $n-k+1$ (\cite[Proposition 7.5]{miller2005combinatorial}).
Furthermore, it is also Cohen-Macaulay (Corollary \ref{prime and cohen macaulayness of toppling ideal}, \cite[Proposition 10.1]{mohammadi2016divisors}). 
Unlike the toppling ideal of a cycle, the toppling ideal of a parcycle can be a prime ideal (Corollary \ref{prime and cohen macaulayness of toppling ideal}). For instance, the defining ideal $P_n$ of the rational normal curve of degree $n$ is the toppling ideal of the parcycle $(C_{n+1},\Pi)$ where $\Pi=\{ \{1,n+1\}, \{2\},\{3\},\dots,\{n\}\}$.  Furthermore,  the notion of $G$-parking function ideal (with respect to a fixed vertex  called the \emph{sink}) also naturally  generalises to the setting of parcycles, we refer to Subsection \ref{Gb and construction of G-parking} for more details. Extending \cite{manjunathschwil}, we explicitly construct minimal free resolutions of the $G$-parking function ideal and the toppling ideal of a parcycle, we refer to Subsection \ref{resolution of G-parking} and Subsection \ref{resolution of toppling ideal}, respectively for more details.  In fact, we construct a complex $\mathcal{F}_{t,\Pi}$ of $R_n[t]$-modules from acyclic orientations on connected partitions of the parcycle
and show the following.

\begin{theorem}\label{minrehomtop}Given any parcycle $(C_n,\Pi)$ and any $G$-parking function ideal $M_{\Pi}$ of $(C_n,\Pi)$, there exists a Gr\"obner degeneration $R_n[t]/I_{\Pi,t}$ of $R_n/I_{\Pi}$ to $R_n/M_{\Pi}$ such that $\mathcal{F}_{t,\Pi}$ is a minimal free resolution of $R_n[t]/I_{\Pi,t}$. 
\end{theorem}

Note that the intermediate complex $\mathcal{F}_{t}(C_{\Pi,\mathcal{M}})$ that appeared in the construction of the complex $\tilde{\mathcal{F}}_{t}(C_{\Pi,\mathcal{M}})$, the protagonist of Theorem \ref{Grodegminres_theo}, is an instance of $\mathcal{F}_{t,\Pi}$ with $C_{\Pi,\mathcal{M}}$ as the underlying parcycle. 

Furthermore, the minimal free resolution of the $G$-parking function ideal of a parcycle (Definition \ref{defGparid}) is a minimal cellular resolution. This cellular resolution is supported on a certain section of the ``graphical hyperplane arrangement'' associated to the parcycle. This is inspired by the work of Dochtermann and Sanyal~\cite{ds} who constructed a minimal cellular resolution for the  $G$-parking function ideal of a graph via the graphical hyperplane arrangement of the underlying graph. We also construct a minimal cellular resolution for any Cohen-Macaulay initial monomial ideal of the rational normal curve (Remark \ref{celres}).  In the following, we discuss some key ideas underlying the cellular resolution. 

Let $(C_n,\Pi)$ be a parcycle where $\Pi=\{V_1,\dots,V_k\}$ is a connected $k$-partition of the vertex set of $C_n$ and $n\in V_k$. For every relevant edge $\{i,j\}$ of $(C_n,\Pi)$, we define a corresponding hyperplane $h_{ij}:=\{ {\mathbf p}=(p_1,\dots,p_n) \in \mathbb{R}^n ~|~ p_{i} =p_{j}\} \subseteq \mathbb{R}^n$. The hyperplane arrangement $\mathcal{A}_{\Pi} = \{h_{ij}~|~ \{i,j\}\, \text{is a relevant}\newline \text{edge of } (C_n,\Pi)\}$ is called the \emph{graphical hyperplane arrangement} of $(C_n,\Pi)$. Let  
$U_{\Pi}=\{ {\bf p} \in \mathbb{R}^n ~|~ p_n=0,~p_{1}+\dots+p_{n-1} = 1, ~p_{u}=p_{w}~\text{whenever}~  \{u,w\}\subseteq V_i
~\text{for some} ~i~ \,\text{from } 1 \text{~to~} k\}
$ be the affine subspace associated with the parcycle $(C_n,\Pi)$. Consider $\tilde{\mathcal{A}}_{\Pi} = \{h_{ij} \cap U_{\Pi} ~|~ \{i,j\} \,\text{is a relevant edge of } (C_n,\Pi)\}$, i.e. the restriction of the hyperplane arrangement $\mathcal{A}_{\Pi}$ to the affine subspace $U_{\Pi}$.
The arrangement $\tilde{\mathcal{A}}_{\Pi}$ divides $U_{\Pi}$ into polyhedra of various dimensions called \emph{cells}. 
Let $B_{\Pi}$ be the collection of all bounded cells of $\tilde{\mathcal{A}}_{\Pi}$ on $U_{\Pi}$. This collection $B_{\Pi}$ forms a polyhedral complex and is contained in $U_{\Pi}$. We  assign the minimal generators of $M_{\Pi}$ to the $0$-cells (vertices) of $B_{\Pi}$, we refer to Subsection \ref{resolution of G-parking} for more details. This assignment gives us a labelled polyhedral complex and we show the following.

\begin{theorem}
    \label{star convexity theorem}
Let $(C_n,\Pi)$ be a parcyle and $B_{\Pi}$ be its labelled polyhedral complex obtained from the hyperplane arrangement $\tilde{\mathcal{A}}_{\Pi}$. The cellular free complex $\mathcal{H}_{\Pi}$ obtained from the labelled polyhedral complex $B_{\Pi}$ is a minimal free resolution of $R/M_{\Pi}$.
\end{theorem}

In fact, the complex $\mathcal{H}_{\Pi}$ is precisely the specialisation $\mathcal{F}_{0,\Pi}$ of $\mathcal{F}_{t,\Pi}$ at $t=0$ (Theorem \ref{complex of G-parking ideal}).

In the following, we put our results into the context of current research.
 \subsection{The Context}
 \begin{itemize}
 \item  {\bf Explicit Minimal Free Resolutions:} The explicit construction of minimal free resolutions of modules over the polynomial ring is a central topic in commutative algebra \cite{miller2005combinatorial}, \cite{eisenbud2013commutative}. Apart from the Koszul complex and complexes such as the Eagon-Northcott complex that are derived from it, there are hardly such explicit minimal free resolutions.  We  refer to \cite{manjunathschwil}, \cite{o2018minimal}, \cite{eagon2019minimal}, \cite{li2022minimal} for recent progress on this topic.  A related theme in algebraic geometry is to extract intrinsic information about a variety from minimal free resolutions of its embeddings in projective space, we refer to \cite{eisenbud2005geometry}, \cite{farkas2022minimal} for more information on this topic. We believe that Theorems \ref{Grodegminres_theo}, \ref{minrehomtop} and \ref{star convexity theorem} are contributions to this theme.

\item {{\bf Multigraded Betti Numbers and Hilbert Series:}}  Multigraded Betti Numbers and Hilbert series of monomial ideals often encode subtle combinatorial information about them. Theorem \ref{Grodegminres_theo} allows us to explicitly compute these invariants for any Cohen-Macaulay initial monomial ideal of a rational normal curve in terms of the underlying parcycle (Subsection \ref{mulgrabnum}). As an application, we explicitly compute a variant of the  Hilbert series of certain lex-segments ideals in two variables continuing a theme initiated by Conca, Negri and Rossi \cite{conca2007contracted}. In the same spirit, via Theorems \ref{minrehomtop} and \ref{star convexity theorem}, we can explicitly determine these invariants for the $G$-parking function ideal of a parcycle.
\end{itemize}

\subsection{Other Applications}

\begin{itemize}

\item {\bf Combinatorics of the Eagon-Northcott Complex Associated to the Rational Normal Curve:}   The defining ideal $P_n$ of the rational normal curve $\Gamma_n$ can be realised as the toppling ideal of the  parcycle $(C_{n+1},\Pi)$ where $\Pi=\{ \{1,n+1\}, \{2\},\{3\},\dots,\{n\}\}$. Hence, the 
 minimal free resolution $\mathcal{F}_{1,\Pi}$ of the  toppling ideal of this parcycle is isomorphic to the Eagon-Northcott complex that minimally resolves the quotient ring of $P_n$.  In Appendix \ref{Eagon northcott complex and toppling complex}, we explicitly construct this isomorphism. Recently, there has been considerable interest in subcomplexes of free resolutions, particularly of fundamental complexes such as the Eagon-Northcott complex \cite{banks2021subcomplexes}. The minimal free resolution $\mathcal{F}_{1,\Pi}$ contains  minimal free resolutions of the toppling ideals of each of its (sub-)parcycles 
as subcomplexes 
(Remark \ref{remsubcom}).  This description is an apt representative for advantages of the combinatorial approach to minimal free resolutions (as adopted in this article).

\item {\bf{Minimal Free Resolution of the Binomial Edge Ideal of $K_n$:}} 
The binomial edge ideal $J_{K_n}$ of $K_n$ \cite{herzog2010binomial}, the complete graph on $n$-vertices, is generated by the set of all $2\times 2$ minors of the matrix 
\begin{equation*}
\begin{bmatrix}
x_1 & x_2 &\dots & x_{n}\\
z_1 & z_2 &\dots &z_{n} 
\end{bmatrix}.
\end{equation*} 

This ideal can  be realised as a toppling ideal of a parcyle, we refer to Example \ref{binedid} for more details.  
This interpretation yields an explicit minimal free resolution of $J_{K_n}$ and the fact that it shares its Betti table with an initial ideal cf. (\cite[Proposition 3.2]{ene2011cohen}, Corollary \ref{equality of betti numbers}).  
\end{itemize}

\subsection{A Future Direction}
The algebraic surface $\Gamma_{n,\mathcal{M}}$ is contained in the universal family over the Hilbert scheme of $\Gamma_n$. Since $P_n$ and $\mathcal{M}$ share their Betti table (Remark \ref{rem2thm1}), $\Gamma_{n,\mathcal{M}}$ is also contained in the universal family over the Betti stratum of $\Gamma_n$. Using this information to study the Hilbert scheme and the Betti stratum of $\Gamma_n$ is a direction for future work.

{\bf Acknowledgements:} We thank Daniele Agostini, R.V. Gurjar, Manoj K. Keshari, Maria A. Mathew, Satoshi Murai, Dipendra Prasad and Kohji Yanagawa for interesting discussions on this topic.
Many thanks to Hidefumi Ohsugi and Kenta Mori for their interest in this work. A part of this work was done while the authors were visiting the International Centre for Theoretical Sciences (ICTS), Bangalore as a part of the program ``Combinatorial Algebraic Geometry: Tropical and Real" in June-July, 2022. We thank ICTS for their kind hospitality. 
MM thanks Waseda University for their kind hospitality in a visit during which this work was presented. The examples were computed using the computer algebra package Macaulay 2. 

\section{Preliminaries}\label{Preliminaries}
This section consists of three subsections. In Subsection \ref{cminit_subsect}, we state a theorem of Conca, Negri and Rossi that classifies Cohen-Macaulay initial monomial ideals of the rational normal curve. We introduce the notion of Gr\"obner degeneration in Subsection \ref{grobdegidealsection}, and
Subsection \ref{preparcyc_subsect} treats basic notions related to parcycles.  

\subsection{Cohen-Macaulay Initial Monomial Ideals of the Rational Normal Curve}\label{cminit_subsect}
Let $\Gamma_n$ be the rational normal curve of degree $n$. Consider the polynomial ring $R_{n+1}=\mathbb{K}[x_1,\dots,x_{n+1}]$ over the field $\mathbb{K}$ where $\mathbb{K}$ is an algebraically closed field. Recall that the defining ideal $P_n$ of the rational normal curve $\Gamma_n$ is generated by the set of all $2\times 2$ minors of the matrix 
\begin{equation}\label{matTn}
T_n = \begin{bmatrix}
x_1 &x_2 &\dots &x_{n}\\
x_2 & x_3 &\dots &x_{n+1} 
\end{bmatrix}.
\end{equation}
The ideal $P_n$ is a Cohen-Macaulay ideal of $R_{n+1}$ (\cite[Theorem 6.4]{eisenbud2005geometry}). The following classification of Cohen-Macaulay initial monomial ideals of $P_n$, due to Conca, Negri and Rossi \cite[Theorem 4.9]{conca2007contracted}, plays an important role in our proof of Theorem \ref{Grodegminres_theo}.

\begin{theorem}\label{Concanegrirossi main theorem}
Let $\mathcal{M}$ be a monomial ideal of $R_{n+1}$. The following conditions are equivalent:
\begin{enumerate}
    \item $\mathcal{M}$ is a Cohen-Macaulay initial ideal of $P_n$.
    \item There exists a sequence $i_0=1<i_1<\dots<i_k=n+1$ such that $\mathcal{M}$ is generated by 
    \begin{enumerate}
        \item the main diagonals $x_{v}x_{r+1}$ of the $2\times 2$ minors of the matrix $T_n$ with column indices $v,r$ such that $v+1\leq i_j<r+1$ for some $j$,
        \item the anti-diagonals $x_{v+1}x_{r}$ of the $2\times 2$ minors of the matrix $T_{n}$ with column indices $v,r$ such that $i_j < v+1<r+1\leq i_{j+1}$ for some $j$.
    \end{enumerate}

\end{enumerate}
\end{theorem}
There is one more equivalent condition in \cite[Theorem 4.9]{conca2007contracted} that we have not mentioned in Theorem \ref{Concanegrirossi main theorem}. In our work, we only need these two equivalent conditions.

\subsection{Gr\"obner Degeneration of an Ideal}\label{grobdegidealsection} 
Let $I$ be an ideal of the polynomial ring $R_n=\mathbb{K}[x_1,\dots,x_n]$ and let $\lambda \in \mathbb{N}^{n}$ be a weight vector on $R_n$. 
Let $\mathcal{M}$ be the initial ideal of $I$ with respect to $\lambda$ (\cite[Subsection 2.1.3]{herzog2011monomial}, \cite[Section 15.2]{eisenbud2013commutative}). For an element $g=\sum a_{\bf u} {\bf x}^{\bf u} \in I $ where $a_{\bf u} \in \mathbb{K}$ and ${\bf u}\in \mathbb{N}^{n}$, we define $g^{h}:=t^{m}g(t^{-\lambda(x_1)}x_1,\dots,t^{-\lambda(x_n)}x_n)$ where $m={\rm{max}}\{{\bf u}\cdot \lambda \mid a_{\bf u}\neq 0\}$.
\begin{definition}\label{GROBDEGDEF}
  The Gr\"obner degeneration of $I$ to $\mathcal{M}$ with respect to $\lambda$ is defined as the ideal $I_{t}=\langle g^{h}\mid g\in I\rangle \subseteq R[t]$ and the Gr\"obner degeneration of
  $R/I$ to $R/\mathcal{M}$ with respect to $\lambda$ is defined as the quotient ring $R[t]/I_t$. 
\end{definition}
For more details on Gr\"obner degenerations of ideals and free $R_n$-modules, 
we refer to \cite[Section 15.8]{eisenbud2013commutative} and \cite[Section 8.3]{miller2005combinatorial}.

\subsection{Parcycles and Related Combinatorics}\label{preparcyc_subsect}
Consider the $n$-cycle graph $C_n=(V(C_{n}),E(C_{n}))$, where $V(C_n)=\{1,\dots,n\}$ is the vertex set and $E(C_n)=\{\{i,i+1\},\{1,n\}\, \mid \,\text{for}\,\, i \,\,\text{from $1$ to }n-1\}$ is the edge set of $C_n$. A \emph{graph isomorphism} between two graphs $G_1$ and $G_2$ is a bijective map $\phi$ from $V(G_1)$ onto $V(G_2)$ with the property that $\{i,j\}\in E(G_1)$ if and only if $\{\phi(i),\phi(j)\}\in E(G_2)$.
We call $C$ a  \emph{cycle graph} on $n$-vertices if it is isomorphic to the $n$-cycle $C_n$. In the following, we generalise cycles to 
objects called \emph{parcycles}.

Let $C$ be a cycle graph on $n$-vertices. A partition $\Pi=\{V_1,\dots,V_k\}$ of the vertex set of $C$ is called a \emph{connected $k$-partition} (or simply, a connected partition) if the subgraph induced by $V_i$ is connected for each $i$ from one to $k$. A \emph{parcycle} is an ordered pair $(C,\Pi)$ where $C$  is a cycle graph on $n$-vertices and $\Pi=\{V_1,\dots,V_k\}$ is a connected partition of the vertex set of $C$. The vertices of $C$ are called \emph{basic vertices} of $(C,\Pi)$ and the elements of $\Pi$ are called \emph{main vertices} of $(C,\Pi)$. An edge of $C$ that connects two basic vertices of $(C,\Pi)$ that lie in the same main vertex of $(C,\Pi)$ is called a \emph{basic edge} of $(C,\Pi)$.
An edge $\{i,j\}$ of $C$ is called a \emph{relevant edge} of $(C,\Pi)$ if $i\in  V_i, j\in V_j$, and $V_i \neq V_j$. 
Relevant edges of $(C,\Pi)$ give  information about the connectivity  of basic vertices lying in two
different main vertices of $(C,\Pi)$. 
We refer to Figure \ref{paroeg} for an example.
\begin{example}\label{parcycmonM}\rm
    Consider the parcycle $C_{\Pi,\mathcal{M}}=(C,\Pi)$ depicted in Figure \ref{paroeg}, where $C$ is a cycle graph on the vertices $1,v_0,3,2,v_1,s_1,s_2$, and $\Pi=\{\{1,v_0\},\{3\},\{2,v_1\},\{s_1,s_2\}
    \}$. The edges $\{s_2,1\},\{v_0,3\},\{3,2\},\{v_1,s_1\}$ are the relevant edges of $C_{\Pi,\mathcal{M}}$ and $\{1,v_0\},\{2,v_1\},\{s_1,s_2\}$ are the basic edges of $C_{\Pi,\mathcal{M}}$.
     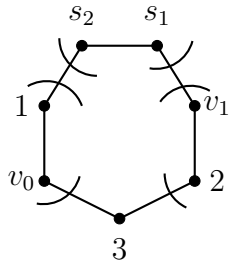
\begin{figure}[ht]
        \centering
       \begin{tikzpicture}
\draw[fill=black](0,0) circle (2pt);
\node at (-0.3,0){$1$};
\draw[fill=black](0,-1) circle (2pt);
\node at (-0.3,-1){$v_0$};
\draw[thick] (0,0)--(0,-1);
\draw[fill=black](1,-1.5) circle (2pt);
\node at (1,-1.9){$3$};
\draw[thick] (0,-1)--(1,-1.5);
\draw[fill=black](2,-1) circle (2pt);
\node at (2.3,-1){$2$};
\draw[thick] (1,-1.5)--(2,-1);
\draw[fill=black](2,0) circle (2pt);
\node at (2.3,0){$v_1$};
\draw[thick] (2,-1)--(2,0);
\draw[fill=black](1.5,0.8) circle (2pt);
\node at (1.5,1.2){$s_1$};
\draw[fill=black](0.5,0.8) circle (2pt);
\node at (0.5,1.2){$s_2$};
\draw[thick] (2,0)--(1.5,0.8);
\draw[thick] (1.5,0.8)--(0.5,0.8);
\draw[thick] (0.5,0.8)--(0,0);
\draw[thick, -](0.5,0) arc (20:120:0.5);
\draw[thick, -](-0.1,-1.3) arc (260:345:0.5);
\draw[thick, -](1.6,-0.9) arc (180:245:0.5);
\draw[thick, -](2.2,0.2) arc (60:145:0.5);
\draw[thick, -](1.4,0.5) arc (260:345:0.5);
\draw[thick, -](0.2,0.9) arc (180:270:0.5);
\end{tikzpicture}
        \caption{The parcycle $C_{\Pi,\mathcal{M}}$.}
        \label{paroeg}
    \end{figure}
    \end{example}

We call $S\subseteq V(C_n)$ a \emph{parset} of $(C_n,\Pi)$ if it is of the form $\cup_{i \in U} V_i$, where $U \subseteq [1,\dots,k]$. We call $S\subseteq V(C_n)$ a \emph{connected-parset} of $(C_n,\Pi)$ if $S$ is a parset of $(C_n,\Pi)$ and the subgraph induced by $S$ on $C_n$ is connected \footnote{The definitions and results on $(C_n,\Pi)$ readily generalise to $(C,\Pi)$ where $C$ is a cycle graph on $n$-vertices. For simplicity of exposition, we deal with the former case throughout.}. The sets $\{1,v_0,2,v_1\}$ and $\{1,v_0,3\}$ in Figure \ref{paroeg} are examples of a parset that is not connected and a connected-parset, respectively.
Given a parset $S\subset V(C_n)$,  we define  ${\bf x}^{S \rightarrow \overline{S}}\in \mathbb{K}[x_1,\dots,x_n]$  to be the monomial $\prod_{i \in S}x_i^{d_i}$, where $d_i=\sum_{j \in \overline{S}}a_{i,j}$ and $\overline{S}$ is the complement of $S$ in $V(C_n)$. Here, $(a_{i,j})$ is the adjacency matrix of $C_{n}$. 
The differentials of minimal free resolutions of the $G$-parking function ideal and the toppling ideal of $(C_n,\Pi)$ can be naturally described in terms of objects called ``acyclic partitions'' of a parcycle $(C_n,\Pi)$ defined in the following. 


{\bf Acyclic Partitions:} Let $\tilde{\Pi}=\{W_1,\dots,W_j\}$ be a connected $j$-partition of the vertex set of $C_n$. The graph $(C_{n})_{\tilde{\Pi}}$ is defined as a graph whose vertex set $V((C_{n})_{\tilde{\Pi}})=\{W_1,\dots,W_j\}$ and the $(W_s,W_l)$-entry of its adjacency matrix is given by $\sum_{u\in W_s,v\in W_l}a_{u,v}$, where $a_{u,v}$ is the $(u,v)$-th entry of the adjacency matrix of $C_n$. The graph $(C_{n})_{\tilde{\Pi}}$ is called the \emph{associated partition graph} of $C_{n}$ with respect to the partition $\tilde{\Pi}$.
An \emph{acyclic $j$-partition} $\mathcal{C}$ of $C_n$ is an ordered pair $(\tilde{\Pi},\mathcal{O})$ where $\tilde{\Pi}$ is a connected $j$-partition of the vertex set of $C_n$ and $\mathcal{O}$ is an acyclic orientation on the associated partition graph $(C_{n})_{\tilde{\Pi}}$.

We define a \emph{connected  $j$-partition} of the parcycle $(C_{n},\Pi)$ as a connected $j$-partition of the vertex set of $C_{n}$ that coarsens $\Pi$, i.e. the elements of this partition are  connected-parsets of $(C_n,\Pi)$. In this context, we refer to the associated connected partition graph as the \emph{connected partition graph} associated with the given connected partition of $(C_n,\Pi)$. An ordered pair $(\tilde{\Pi},\mathcal{O})$ is called an \emph{acyclic $j$-partition} of the parcycle $(C_n,\Pi)$ if $\tilde{\Pi}$ is a connected $j$-partition of $(C_n,\Pi)$ and $\mathcal{O}$ is an acyclic orientation on $(C_{n})_{\tilde{\Pi}}$. An acyclic $j$-partition $\mathcal{C}=(\tilde{\Pi},\mathcal{O})$ of $(C_n,\Pi)$ is called a \emph{$V_k$-acyclic $j$-partition} if $\mathcal{O}$ is an acyclic orientation on $(C_n)_{\tilde{\Pi}}$ with a unique sink at the vertex containing $V_k$.
We refer the reader to Figure \ref{acyclic partitions} for examples of acyclic partitions.

{\bf Chip Firing on Parcycles:}  The chip firing game on a parcycle $(C_n,\Pi)$ is identical to playing chip firing game on $C_n$, except that only specific chip firing moves are permitted.
In this restricted chip firing game on $C_n$, we are only allowed to fire from the parsets of $(C_n,\Pi)$. A \emph{configuration} on the parcycle $(C_n,\Pi)$ is an assignment of an integer to the vertices of $C_n$, i.e. the basic vertices of $(C_n,\Pi)$. We also refer to it as a \emph{divisor} on $(C_n,\Pi)$.
Two divisors $D_1$ and $D_2$ on $(C_n,\Pi)$ are called \emph{equivalent} if we can obtain $D_2$ from $D_1$ by a sequence of restricted chip firing moves on $C_n$. As elements of $\mathbb{Z}^n$, two elements $D_1,D_2\in \mathbb{Z}^n$ are equivalent if and only if $D_1-D_2\in L_{\Pi}$, where $L_{\Pi}$ is the Laplacian lattice of the parcycle $(C_n,\Pi)$. Recall that $L_{\Pi}$
is the lattice generated by vectors of the form $\sum_{j \in V_i}b_j$
where $b_j$ a typical row of the Laplacian matrix of $C_n$.
For more details on the chip firing game and acyclic orientations on $C_n$, the reader is referred to \cite[Section 2]{manjunathschwil}.

Let $\mathcal{C}=(\tilde{\Pi},\mathcal{O})$ be an acyclic $j$-partition of $(C_n,\Pi)$, where
$\tilde{\Pi} = \{W_1,\dots,W_j\}$ is the connected $j$-partition of $(C_{n},\Pi)$ associated to $\mathcal{C}$. Let $i\in W_s$, then ${\rm{out}}_{\mathcal{C}}(i)=\sum_{(W_s,W_l)\in \mathcal{C}}(\sum_{l \in W_l}a_{i,l})$ is called the out-degree of $i$ in $\mathcal{C}$, where $(a_{i,l})$ is the adjacency matrix of $C_n$ and $(W_s,W_l)$ is a directed edge in $\mathcal{C}$. We associate a divisor $D(\mathcal{C})$ on $(C_n,\Pi)$ as follows:
$D(\mathcal{C})=\sum_{i\in C_n}{\rm{out}}_{\mathcal{C}}(i)[i]$, where ${\rm out}_{\mathcal{C}}(i)$ is the out-degree of $i$ in $\mathcal{C}$.
 
Note that acyclic $j$-partitions of $(C_n,\Pi)$ are also acyclic $j$-partitions of $C_n$. Two acyclic $j$-partitions $\mathcal{C}_1,\mathcal{C}_2$ of $(C_n,\Pi)$ are \emph{chip firing equivalent} if they are chip firing equivalent as acyclic $j$-partitions of $C_n$ \cite[Definition 2.2]{manjunathschwil}. As per the definition of the chip firing equivalence of two acyclic partitions of $(C_n,\Pi)$, the chip firing equivalence class $[\mathcal{C}]$ of an acyclic $j$-partition $\mathcal{C}$ of $(C_n,\Pi)$ is the same as the chip firing equivalence class of $\mathcal{C}$ in $C_n$. From \cite[Lemma 2.3]{manjunathschwil}, for a fixed vertex $n\in V_k$, every chip firing equivalence class of acyclic $j$-partitions of $(C_n,\Pi)$ contains a unique $n$-acyclic $j$-partition. From \cite[Lemma 2.4]{manjunathschwil}, if $\mathcal{C}_1,\mathcal{C}_2$ are two acyclic $j$-partitions of $(C_n,\Pi)$ such that $\tilde{\Pi}(\mathcal{C}_1)=\tilde{\Pi}(\mathcal{C}_2)$ and $\mathcal{C}_1 \sim \mathcal{C}_2$, then $D(\mathcal{C}_1)$ and  $D(\mathcal{C}_2)$ are equivalent as divisors of $(C_n,\Pi)$. 
\begin{example}\rm
   The parcycle $C_{\Pi,\mathcal{M}}$ defined in Example \ref{parcycmonM} has a unique $1$-partition, six chip firing equivalence classes of acyclic $2$-partitions, eight chip firing equivalence classes of acyclic $3$-partitions, and three chip firing equivalence classes of acyclic $4$-partitions. The acyclic representative of these equivalence classes with the sink vertex that contains $\{s_1,s_2\}$ is illustrated in Figure \ref{acyclic partitions}.
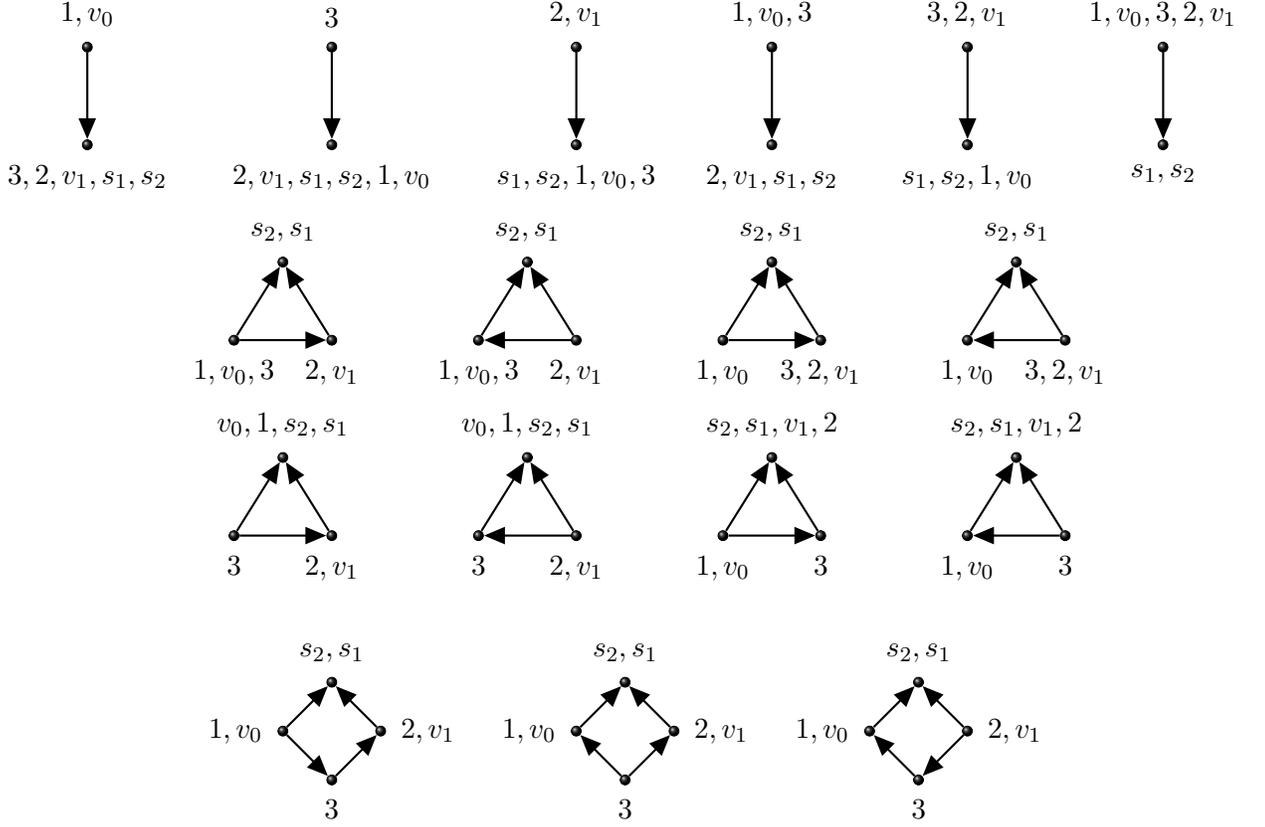
\begin{figure}[ht] 
\begin{center}
\begin{tikzpicture}[scale=1.3]
\SetVertexMath
\GraphInit[vstyle=Art]
\SetUpVertex[MinSize=3pt]
\SetVertexLabel
\tikzset{VertexStyle/.style = {%
shape = circle,
shading = ball,
ball color = black,
inner sep = 1.5pt
}}
\SetUpEdge[color=black]
\tikzstyle{every node}=[font=\small]

\Vertex[LabelOut,Lpos=90, Ldist=.05cm,x=-1,y=9,L={1,v_0}]{a}
\Vertex[LabelOut,Lpos=270, Ldist=.05cm,x=-1,y=8,L={3,2,v_1,s_1,s_2}]{b}
\Edge[style={-triangle 45}](a)(b)

\Vertex[LabelOut,Lpos=90, Ldist=.05cm,x=1.5,y=9,L={3}]{a}
\Vertex[LabelOut,Lpos=270, Ldist=.05cm,x=1.5,y=8,L={2,v_1,s_1,s_2,1,v_0}]{b}
\Edge[style={-triangle 45}](a)(b)

\Vertex[LabelOut,Lpos=90, Ldist=.05cm,x=4,y=9,L={2,v_1}]{a}
\Vertex[LabelOut,Lpos=270, Ldist=.05cm,x=4,y=8,L={s_1,s_2,1,v_0,3}]{b}
\Edge[style={-triangle 45}](a)(b)

\Vertex[LabelOut,Lpos=90, Ldist=.05cm,x=6,y=9,L={1,v_0,3}]{a}
\Vertex[LabelOut,Lpos=270, Ldist=.05cm,x=6,y=8,L={2,v_1,s_1,s_2}]{b}
\Edge[style={-triangle 45}](a)(b)

\Vertex[LabelOut,Lpos=90, Ldist=.05cm,x=8,y=9,L={3,2,v_1}]{a}
\Vertex[LabelOut,Lpos=270, Ldist=.05cm,x=8,y=8,L={s_1,s_2,1,v_0}]{b}
\Edge[style={-triangle 45}](a)(b)

\Vertex[LabelOut,Lpos=90, Ldist=.05cm,x=10,y=9,L={1,v_0,3,2,v_1}]{a}
\Vertex[LabelOut,Lpos=270, Ldist=.05cm,x=10,y=8,L={s_1,s_2}]{b}
\Edge[style={-triangle 45}](a)(b)


\Vertex[LabelOut,Lpos=270, Ldist=.05cm,x=0.5,y=6,L={1,v_0,3}]{2}
\Vertex[LabelOut,Lpos=270, Ldist=.05cm,x=1.5,y=6,L={2,v_1}]{1}
\Vertex[LabelOut,Lpos=90, Ldist=.05cm,x=1,y=6.8,L={s_2,s_1}]{34}
\Edge[style=-{triangle 45}](2)(1)
\Edge[style={-triangle 45}](1)(34)
\Edge[style={-triangle 45}](2)(34)

\Vertex[LabelOut,Lpos=270, Ldist=.05cm,x=3,y=6,L={1,v_0,3}]{2}
\Vertex[LabelOut,Lpos=270, Ldist=.05cm,x=4,y=6,L={2,v_1}]{1}
\Vertex[LabelOut,Lpos=90, Ldist=.05cm,x=3.5,y=6.8,L={s_2,s_1}]{34}
\Edge[style=-{triangle 45}](1)(2)
\Edge[style={-triangle 45}](1)(34)
\Edge[style={-triangle 45}](2)(34)

\Vertex[LabelOut,Lpos=270, Ldist=.05cm,x=5.5,y=6,L={1,v_0}]{2}
\Vertex[LabelOut,Lpos=270, Ldist=.05cm,x=6.5,y=6,L={3,2,v_1}]{1}
\Vertex[LabelOut,Lpos=90, Ldist=.05cm,x=6,y=6.8,L={s_2,s_1}]{34}
\Edge[style=-{triangle 45}](2)(1)
\Edge[style={-triangle 45}](1)(34)
\Edge[style={-triangle 45}](2)(34)

\Vertex[LabelOut,Lpos=270, Ldist=.05cm,x=8,y=6,L={1,v_0}]{2}
\Vertex[LabelOut,Lpos=270, Ldist=.05cm,x=9,y=6,L={3,2,v_1}]{1}
\Vertex[LabelOut,Lpos=90, Ldist=.05cm,x=8.5,y=6.8,L={s_2,s_1}]{34}
\Edge[style=-{triangle 45}](1)(2)
\Edge[style={-triangle 45}](1)(34)
\Edge[style={-triangle 45}](2)(34)


\Vertex[LabelOut,Lpos=270, Ldist=.05cm,x=0.5,y=4,L={3}]{2}
\Vertex[LabelOut,Lpos=270, Ldist=.05cm,x=1.5,y=4,L={2,v_1}]{1}
\Vertex[LabelOut,Lpos=90, Ldist=.05cm,x=1,y=4.8,L={v_0,1,s_2,s_1}]{34}
\Edge[style=-{triangle 45}](2)(1)
\Edge[style={-triangle 45}](1)(34)
\Edge[style={-triangle 45}](2)(34)

\Vertex[LabelOut,Lpos=270, Ldist=.05cm,x=3,y=4,L={3}]{2}
\Vertex[LabelOut,Lpos=270, Ldist=.05cm,x=4,y=4,L={2,v_1}]{1}
\Vertex[LabelOut,Lpos=90, Ldist=.05cm,x=3.5,y=4.8,L={v_0,1,s_2,s_1}]{34}
\Edge[style=-{triangle 45}](1)(2)
\Edge[style={-triangle 45}](1)(34)
\Edge[style={-triangle 45}](2)(34)

\Vertex[LabelOut,Lpos=270, Ldist=.05cm,x=5.5,y=4,L={1,v_0}]{2}
\Vertex[LabelOut,Lpos=270, Ldist=.05cm,x=6.5,y=4,L={3}]{1}
\Vertex[LabelOut,Lpos=90, Ldist=.05cm,x=6,y=4.8,L={s_2,s_1,v_1,2}]{34}
\Edge[style=-{triangle 45}](2)(1)
\Edge[style={-triangle 45}](1)(34)
\Edge[style={-triangle 45}](2)(34)

\Vertex[LabelOut,Lpos=270, Ldist=.05cm,x=8,y=4,L={1,v_0}]{2}
\Vertex[LabelOut,Lpos=270, Ldist=.05cm,x=9,y=4,L={3}]{1}
\Vertex[LabelOut,Lpos=90, Ldist=.05cm,x=8.5,y=4.8,L={s_2,s_1,v_1,2}]{34}
\Edge[style=-{triangle 45}](1)(2)
\Edge[style={-triangle 45}](1)(34)
\Edge[style={-triangle 45}](2)(34)

\Vertex[LabelOut,Lpos=90, Ldist=.05cm,x=1.5,y=2.5,L={s_2,s_1}]{1}
\Vertex[LabelOut,Lpos=180, Ldist=.05cm,x=1,y=2,L={1,v_0}]{2}
\Vertex[LabelOut,Lpos=0, Ldist=.05cm,x=2,y=2,L={2,v_1}]{3}
\Vertex[LabelOut,Lpos=270, Ldist=.05cm,x=1.5,y=1.5,L={3}]{4}
\Edge[style={-triangle 45}](2)(1)
\Edge[style={-triangle 45}](3)(1)
\Edge[style={-triangle 45}](2)(4)
\Edge[style={-triangle 45}](4)(3)

\Vertex[LabelOut,Lpos=90, Ldist=.05cm,x=4.5,y=2.5,L={s_2,s_1}]{1}
\Vertex[LabelOut,Lpos=180, Ldist=.05cm,x=4,y=2,L={1,v_0}]{2}
\Vertex[LabelOut,Lpos=0, Ldist=.05cm,x=5,y=2,L={2,v_1}]{3}
\Vertex[LabelOut,Lpos=270, Ldist=.05cm,x=4.5,y=1.5,L={3}]{4}
\Edge[style={-triangle 45}](2)(1)
\Edge[style={-triangle 45}](3)(1)
\Edge[style={-triangle 45}](4)(2)
\Edge[style={-triangle 45}](4)(3)

\Vertex[LabelOut,Lpos=90, Ldist=.05cm,x=7.5,y=2.5,L={s_2,s_1}]{1}
\Vertex[LabelOut,Lpos=180, Ldist=.05cm,x=7,y=2,L={1,v_0}]{2}
\Vertex[LabelOut,Lpos=0, Ldist=.05cm,x=8,y=2,L={2,v_1}]{3}
\Vertex[LabelOut,Lpos=270, Ldist=.05cm,x=7.5,y=1.5,L={3}]{4}
\Edge[style={-triangle 45}](2)(1)
\Edge[style={-triangle 45}](3)(1)
\Edge[style={-triangle 45}](4)(2)
\Edge[style={-triangle 45}](3)(4)
\end{tikzpicture}
\end{center}
\caption{\label{acyclic partitions} The acyclic $j$-partitions of the parcycle $C_{\Pi,\mathcal{M}}$ for $j=2,3,4$, where $\{s_1,s_2\}$ is contained in the sink vertex.}
\end{figure}\qed
\end{example}

{\bf Edge contractions:} A directed edge $f$ of an acyclic $j$-partition $\mathcal{C}$ of $(C_n,\Pi)$ is called a \emph{contractible edge} of $\mathcal{C}$ if the directed graph $\mathcal{C}/f$ (obtained from contracting the edge $f$ of $\mathcal{C}$) is acyclic. An edge $f$ is called a \emph{contractible edge} in the chip firing equivalence class $\textbf{c}$ of acyclic $j$-partitions of $(C_n,\Pi)$ if there exists an acyclic $j$-partition $\mathcal{C}\in \textbf{c}$ such that $f\in \mathcal{C}$ is a contractible edge. For more details on contractible edges of a graph, we refer the reader to  \cite[Section 2]{manjunathschwil}.

Let $\textbf{c}$ be a chip firing equivalence class of acyclic $j$-partitions  of $(C_n,\Pi)$ and let $f$ be a contractible edge in $\textbf{c}$. Let $\mathcal{C}\in \textbf{c}$ be an acyclic $j$-partition such that $f \in \mathcal{C}$ is a contractible edge. The monomial associated with the edge contraction of $f$ in $\textbf{c}$ is defined as $m_{\textbf{c}}(f)= \textbf{x}^{D(\mathcal{C})-D(\mathcal{C}/f)}$. The well-definedness of the monomial $m_{\textbf{c}}(f)$ holds from \cite[Section 2]{manjunathschwil}. From \cite[Proposition 2.9]{manjunathschwil}, there exist functions ${\rm{sign}}_{\textbf{c}}$ for every chip firing equivalence class $\textbf{c}$ of acyclic $j$-partitions of $(C_n,\Pi)$, taking contractible edges of $\textbf{c}$ to $\{1,-1\}$. These ${\rm{sign}}_{\textbf{c}}$ functions satisfy the following conditions:
\begin{itemize}
    \item ${\rm{sign}}_{\textbf{c}}(f_1){\rm{sign}}_{\textbf{c}/f_1}(f_2)=-   {\rm{sign}}_{\textbf{c}}(f_2){\rm{sign}}_{\textbf{c}/f_2}(f_1)$ where $f_1,f_2$ are distinct contractible edges in $\textbf{c}$.
    \item ${\rm{sign}}_{\textbf{c}}(f)=-{\rm{sign}}_{\textbf{c}}(\tilde{f})$ where $f=(U,W), \tilde{f}=(W,U)$ are contractible edges in $\textbf{c}$.
\end{itemize}

\section{Commutative Algebra of Parcycles}\label{ComAlgPar}

In the opening subsection (Subsection \ref{relation between one-generic matrix and toppling ideal}), we associate a certain matrix of linear forms to a parcycle. This matrix will come in handy in various arguments in subsequent sections. We prove results related to Gr\"obner bases of toppling ideals in Section \ref{Gb and construction of G-parking}. In Sections \ref{resolution of G-parking} and \ref{resolution of toppling ideal}, we construct combinatorial minimal free resolutions for the $G$-parking function ideals and the toppling ideals of parcycles, respectively.

\subsection{One-Generic Matrices and Parcycles}\label{relation between one-generic matrix and toppling ideal}
We define a relation between the toppling ideal of a parcycle and the ideal generated by the set of all $2\times 2$ minors of a certain matrix associated with the parcycle. Before this, we recall the definition of a one-generic matrix.

Let $x_1,x_2,\dots,x_s$ be indeterminates over a field $\mathbb{K}$. Let $T$ be a matrix of order $m\times n$ such that $m\leq n$, and its entries are linear forms in $R_s(=\mathbb{K}[x_1,\dots,x_s])$ where $s\geq m+n-1$. A \emph{generalised row} of $T$ is a nonzero $\mathbb{K}$-linear combination of the rows of $T$. A \emph{generalised entry} of $T$ is a nonzero $\mathbb{K}$-linear combination of the entries of a generalised row of $T$.

\begin{definition}
A matrix $T$ comprising entries from the polynomial ring $R_s$ of linear forms is defined as $1$-generic if all of its generalised entries are nonzero.
\end{definition}

One-generic matrices have many nice properties. For example, if $T$ is a one-generic matrix of linear forms in $R_{s}$ of order $m\times n$ with $m\leq n$ over an algebraically closed field $\mathbb{K}$, then the ideal generated by the set of all maximal minors of $T$ is a Cohen-Macaulay prime ideal in $R_{s}$ (\cite[Theorem 6.4]{eisenbud2005geometry}).

In the following, we relate matrices of linear forms to parcycles.
Let $(C_n,\Pi)$ be a parcycle where $\Pi=\{V_1,\dots,V_k\}$ is a connected $k$-partition of $V(C_n)$. Without loss of generality, we assume that $n\in V_k$. 
Let $\tilde{E}(C_n)=\{(i,i+1),(n,1)\,\mid 1\leq i \leq n-1\}$ be the set of all directed edges of $C_n$. We associate a matrix $A_{(C_n,\Pi)}$ of order $2 \times k$ to the parcycle $(C_n,\Pi)$ as follows: for a fixed main vertex of $(C_n,\Pi)$, say $V_j=\{j_1,\dots,j_l\}$, we get a directed closed path $j_1\dots j_l \dots j_1$ in $C_n$ obtained from following the directed edges of $C_n$ with  $j_1$ as the starting vertex. This directed closed path gives an enumeration of the relevant edges of $(C_n,\Pi)$. The $i$-th column of $A_{(C_n,\Pi)}$ is $(x_{j_p}\,\,x_{j_q})^{t}$ if $(j_p,j_q)$ is the $i$-th directed relevant edge in the directed closed path $j_1\dots j_l \dots j_1$.

For example, consider the parcycle $C_{\Pi,\mathcal{M}}$ defined in Example \ref{parcycmonM}. Let $w_1=1,w_2=v_0,w_3=3,w_4=2,w_5=v_1,w_6=s_1,w_7=s_2$ be the labelling of the basic vertices of $C_{\Pi,\mathcal{M}}$ and $\{(w_i,w_{i+1}),(w_7,w_1)\mid 1\leq i\leq 6\}$ is the set of all directed edges of $C$. The matrix $A_{(C,\Pi)}$ associated with $(C,\Pi)$ (considering $\{s_1,s_2\}$ as the fixed main vertex) is given by

We now provide a necessary and sufficient condition for the matrix $A_{(C_n,\Pi)}$ to be a $1$-generic matrix. 

\begin{lemma}\label{one genereic of A}
Let $\Pi=\{V_1,\dots,V_k\}$ be a connected $k$-partition  of the vertex set of $C_n$ such that there exists an $i$ from $1$ to $k$ with $|V_i|\geq 2$. The matrix $A_{(C_n,\Pi)}$ is a $1$-generic matrix.
\end{lemma}
\begin{proof} Without loss of generality, we fix a main vertex, say $V_k$, of $(C_n,\Pi)$. Let 
\begin{equation}\label{matT}
T=
\begin{bmatrix}
x_{1} &x_{2} &\dots &x_{n}\\
x_{2} & x_{3} &\dots &x_{1} 
\end{bmatrix}.
\end{equation} From the construction of $A_{(C_n,\Pi)}$, if $|V_i|\geq 2$ for some $i$, then $A_{(C_n,\Pi)}$ is obtained from the matrix $T$ by removing and permuting certain columns. We now prove that if we remove a column from $T$, then it becomes a $1$-generic matrix. The $1$-genericity of a matrix obtained from removing a column of $T$ is equivalent to the $1$-genericity of the matrix $$
\tilde{T}=\begin{bmatrix}
x_{2} &x_{3} &\dots &x_{n}\\
x_{3} & x_{4} &\dots &x_{1} 
\end{bmatrix}.
$$
A generalised entry of $\tilde{T}$ is of the form \begin{equation}\label{eq1}
    \sum_{i=2}^{n-1}a_i(k_1x_{i}+k_2x_{{i+1}})+a_n(k_1x_{{n}}+k_2x_{{1}}) =0,
\end{equation} where $k_1,k_2\in \mathbb{K}$ and at least one of them is nonzero. Equation (\ref{eq1}) gives us a family of equations $a_2k_1=0, a_nk_2=0, a_ik_1+a_{i-1}k_2=0,$ for all $3\leq i \leq n$. If both $k_1,k_2$ are nonzero, then this new family of equations gives us $a_i=0$, for all $2\leq i\leq n.$ In the case of $k_1=0$ or $k_2=0$, $a_i=0$ for all $i$ since there is no repetition of variables in any of the rows of $\tilde{T}$ and $x_{i}$ is an indeterminate over $\mathbb{K}$ for all $i$ from $1$ to $n$. Hence, $\tilde{T}$ is a $1$-generic matrix. This implies that removing any column of $T$ converts it into a $1$-generic matrix. As the matrix $A_{(C_n,\Pi)}$ is a sub-matrix of a $1$-generic matrix obtained by removing certain columns, it is also $1$-generic.
\end{proof}
\begin{remark}{\rm
    Note that the matrix $T$ defined in Equation (\ref{matT}) is not $1$-generic. The matrix $T$ is the matrix $A_{(C_n,\tilde{\Pi})}$ associated with the parcycle $(C_n,\tilde{\Pi})$ where $\tilde{\Pi}$ is the connected $n$-partition of $V(C_n)$ and $\{1\}$ is the fixed main vertex of $(C_n,\tilde{\Pi})$}.\qed
\end{remark}

Next, we define a bijection between the set of all $2\times 2$ minors of the matrix $A_{(C_n,\Pi)}$, and a certain type of connected-parsets of $(C_n,\Pi)$. Let $S=\{j_{s},\dots,j_{t}\}$ be a connected-parset of $(C_n,\Pi)$ such that $S$ does not contain $V_{k}$. Let $(j_{u},j_{s}), (j_{t},j_{w})$ be the directed relevant edges in $(C_n,\Pi)$ connecting $S$ and $\overline{S}$. Let $A_{(C_n,\Pi)}$ be the matrix associated with the parcycle $(C_n,\Pi)$ with the fixed main vertex $V_k$. Let $A_{ij}$ be the minor corresponding to the $i$-th and $j$-th columns of $A_{(C_n,\Pi)}$. 
Consider the map $\phi$ between the set  $\{\mathbf{x}^{S\rightarrow \overline{S}}-\mathbf{x}^{\overline{S}\rightarrow S}\mid~S~ \text{is a connected-parset of }\, (C_n,\Pi)\, \text{not containing}\, V_{k} \}$ and the set of all $2\times 2$ minors $A_{ij}$ (with $i<j$) of $A_{(C_n,\Pi)}$ that is defined as \begin{equation*}
 \phi(\mathbf{x}^{S\rightarrow \overline{S}}-\mathbf{x}^{\overline{S}\rightarrow S})={\rm{det}}(\begin{bmatrix}
x_{j_{u}}&x_{j_{t}}\\
x_{j_{s}}&x_{j_{w}}
\end{bmatrix}).
\end{equation*}
Note that the map $\phi$ is a well-defined map.

\begin{lemma}\label{bijection lemma} Let $(C_n,\Pi)$ be a parcycle and $V_k$ is the fixed main vertex in $(C_n,\Pi)$.
The map $\phi$ is a bijection.
\end{lemma}
\begin{proof} 
Note that for a connected-parset $S=\{j_{s},\dots,j_{t}\}$ of $(C_n,\Pi)$ that does not contain $V_{k}$
\begin{equation}\label{eq 3}
 (\mathbf{x}^{S\rightarrow \overline{S}}-\mathbf{x}^{\overline{S}\rightarrow S})
 =-{\rm{det}}(\begin{bmatrix}
x_{j_{u}}&x_{j_{t}}\\
x_{j_{s}}&x_{j_{w}}
\end{bmatrix}),
\end{equation} where $(j_{u},j_{s}), (j_{t},j_{w})$ are the directed relevant edges in $(C_n,\Pi)$ connecting $S$ and $\overline{S}$.
For two distinct connected-parsets $S_1,S_2$ of $(C_n,\Pi)$ not containing $V_k$, the binomials ${\bf{x}}^{S_1\rightarrow \overline{S}_1}-{\bf{x}}^{ \overline{S}_1\rightarrow S_1}$, ${\bf{x}}^{S_2\rightarrow \overline{S}_2}-{\bf{x}}^{ \overline{S}_2\rightarrow S_2}$ are different. This implies $\phi$ is an injective map.

Let $A_{ij}$ be the minor corresponding to the $i$-th and $j$-th columns of $A_{(C_n,\Pi)}$ (note that $i<j$). Let $(x_{j_{a}}\,\,x_{j_{b}})^{t},(x_{j_u}\,\,x_{j_{w}})^{t}$ be the $i$-th and $j$-th columns of $A_{(C_n,\Pi)}$, respectively. This implies $(j_{a},j_{b})$ and $(j_{u},j_{w})$ are directed relevant edges in $(C_n,\Pi)$.
From the construction of $A_{(C_n,\Pi)}$, there exists a directed path $j_{b}\dots j_{u}$ in $(C_n,\Pi)$
such that the set
$S=\{j_{b},\dots,j_{u}\}$ is a connected-parset of $(C_n,\Pi)$ and $V_{k}$ is contained in the connected-parset $\overline{S}$. The minor corresponding to $i$-th, $j$-th columns is $A_{ij} =x_{j_{a}}x_{j_{w}}-x_{j_{b}}x_{j_{u}} = - (\mathbf{x}^{S\rightarrow \overline{S}}-\mathbf{x}^{\overline{S}\rightarrow S})$ and this holds for all $2\times 2$ minors of $A_{(C_n,\Pi)}$. Hence, $\phi$ is a surjective map.
\end{proof}
Next, we identify the set of monomials $\{\mathbf{x}^{S\rightarrow \overline{S}} \mid~S~ \text{ is a connected-parset of }\newline (C_n,\Pi)\, \text{not containing}\, V_{k} \}$  with the set of anti-diagonals of $2\times 2$ minors of $A_{(C_n,\Pi)}$.

\begin{proposition}\label{one generic matrix to parking ideal}
Consider the parcycle $(C_n,\Pi)$ and let $V_{k}$ be the fixed main vertex of $(C_n,\Pi)$.
The set $\tilde{M}=\{\mathbf{x}^{S\rightarrow \overline{S}} \mid~S~ \text{is a connected-parset of}~ (C_n,\Pi)~ \text{not}\newline \text{containing}\,\, V_{k} \}$ is equal to the set of anti-diagonals of $2\times 2$ minors $A_{ij}$ of $A_{(C_n,\Pi)}$.
\end{proposition}
\begin{proof} The proof of Proposition \ref{one generic matrix to parking ideal} is immediate from Equation (\ref{eq 3}).
 \end{proof}
 \subsection{Gr\"obner Basis of the Toppling Ideal of a Parcycle}\label{Gb and construction of G-parking}
  In this subsection, we construct a Gr\"obner basis of the toppling ideal of a parcycle in terms of its connected-parsets. Henceforth, throughout this subsection, we work with the assumption that $(C_n,\Pi)$ is a parcycle where $\Pi=\{V_1,\dots,V_k\}$ is a connected $k$-partition of $V(C_n)$ and $n \in V_k$.

We begin by constructing a generating set for the toppling ideal $I_{\Pi}$, utilising the connected-parsets of $(C_n,\Pi)$.
We know that the toppling ideal of a graph $G$ is generated by the set $\{\textbf{x}^{S \rightarrow \overline S} - \textbf{x}^{ \overline S  \rightarrow S} \mid~S~ \text{ is a connected-parset of}\, G \newline \text{not containing}~ n \}$ (\cite[Theorem 25]{manjunath2013monomials}, \cite{MR1896344}). Hence, in case of the connected $n$-partition $\Pi$ of the vertex set of $C_n$, the toppling ideal $I_{\Pi}=
 \langle \textbf{x}^{S \rightarrow \overline S} - \textbf{x}^{ \overline S  \rightarrow S} \mid~S~ \text{ is a connected-parset of}\, (C_n,\Pi)~ \text{not containing}~ n \rangle$.

Let $|V_i|\geq 2$, for some $i$ and let $b_j$ be the row corresponding to the vertex $j$ in $\Lambda_{ C_n}$.
Recall that the toppling ideal $I_{\Pi}$ of $(C_n,\Pi)$ is defined as the lattice ideal of $L_{\Pi}$, where $L_{\Pi}$ is the Laplacian lattice of the parcycle $(C_{n},\Pi)$. The set $L=\{\sum_{i\in \tilde{S}}b_{V_i}\mid \tilde{S}\subset \{1,\dots,k\},\, S=\cup_{i\in \tilde{S}}V_i\,\, \text{is a connected-parset of }\,  (C_n,\Pi)\}$ forms a generating set of the lattice $L_{\Pi}$. Let $U\subset V(C_n)$. We define a vector $(U \rightarrow \overline{U})=(u_1,\dots,u_n)\in \mathbb{Z}^n$ where $u_i$ is the total number of edges between the vertex $i$ and the vertices in $\overline{U}$, i.e. $u_i=\sum_{j\in \overline{U}}a_{i,j}\rchi_{U}(i)$. 
Let $ \tilde{S}\subset \{1,\dots,k\}$ be such that $S=\cup_{i\in \tilde{S}}V_i$ is a connected-parset of $(C_n,\Pi)$. Consequently, $(S\rightarrow \overline{S})-(\overline{S}\rightarrow S) = \sum_{i\in \tilde{S}}b_{V_i}$. Let $I_{L}$ be the ideal generated by the set $\{\mathbf{x}^{S\rightarrow \overline{S}}-\mathbf{x}^{\overline{S}\rightarrow S}\mid S~\text{is a connected-parset of}\, (C_n,\Pi)\}$.
From (\cite[Lemma 7.6]{miller2005combinatorial}), $I_{\Pi}=(I_{L}\colon \langle x_1\cdots x_n\rangle ^{\infty})$. Let $k\in \tilde{S}(\subset \{1,\dots,k\})$ and $S=\cup_{i\in \tilde{S}}V_{i}$ is a connected-parset of $(C_n,\Pi)$, then $\mathbf{x}^{S\rightarrow \overline{S}}-\mathbf{x}^{\overline{S}\rightarrow S}=-(\mathbf{x}^{\overline{S}\rightarrow S}-\mathbf{x}^{S\rightarrow \overline{S}}),$ where $\overline{S}$ is a connected-parset of $(C_n,\Pi)$ that does not contain $V_k$. 
This implies that $I_{L}=\langle\mathbf{x}^{S\rightarrow \overline{S}}-\mathbf{x}^{\overline{S}\rightarrow S}\mid~S~ \text{is a connected-parset of}\, (C_n,\Pi) \,\text{not containing}\, V_k \rangle $. From Lemma \ref{bijection lemma}, the ideal $I_{L}$ is generated by the set of all $2\times 2$ minors of the $1$-generic matrix $A_{(C_n,\Pi)}$. Hence, $I_L$ is a prime ideal (Lemma \ref{one genereic of A}, \cite[Theorem 6.4]{eisenbud2005geometry})).

\begin{lemma}\label{generating lemma}
Let $(C_n,\Pi)$ be a parcycle where $\Pi=\{V_1,\dots,V_k\}$ is a connected $k$-partition of $C_n$ and $n\in V_k$. The toppling ideal $I_{\Pi}=\langle \mathbf{x}^{S\rightarrow \overline{S}}-\mathbf{x}^{\overline{S}\rightarrow S}\mid~S~ \text{ is a connected-parset of }\, (C_n,\Pi) \,\text{not containing}\, V_k \rangle$.
\end{lemma}
\begin{proof}
In the case of the connected $n$-partition $\Pi$, the result holds from \cite[Theorem 25]{manjunath2013monomials}. Let $|V_i|\geq 2$, for some $i$. In this case, $I_{\Pi}=(I_{L}\colon \langle x_1\cdots x_n\rangle ^{\infty})$ where $I_{L}=\langle\mathbf{x}^{S\rightarrow \overline{S}}-\mathbf{x}^{\overline{S}\rightarrow S}\mid~S~ \text{ is a connected-parset of} \,\,(C_n,\Pi)\,\text{not containing}\newline V_k \rangle$. Furthermore, the ideal $I_{L}$ is a prime ideal. To prove that $I_{\Pi}\subseteq I_L$, it enough to show that $(x_1\cdots x_n)^{m}\notin I_L$ for any $m\geq 1$. Suppose that $(x_1\cdots x_n)^{m}\in I_L$, this implies that $\mathbb{V}(I_{L})\subseteq \mathbb{V}(x_i)$ for some $i$. As $I_L$ is generated by the set of all $2\times 2$ minors of the matrix $A_{(C_n,\Pi)}$, this implies
$(1,1,\dots,1)\in \mathbb{V}(I_{L})$, but $(1,1,\dots,1) \notin \mathbb{V}(x_i)$. This yields a contradiction to $\mathbb{V}(I_{L})\subseteq \mathbb{V}(x_i)$.
Hence, the toppling ideal is generated by the set $\{\mathbf{x}^{S\rightarrow \overline{S}}-\mathbf{x}^{\overline{S}\rightarrow S}\mid~S~ \text{ is a connected-parset of }\, (C_n,\Pi) \,\text{not containing}\,\, V_k\}$. \end{proof}

Next, we prove that the set $\{\mathbf{x}^{S\rightarrow \overline{S}}\mid~S \text{ is a connected-parset of} \,(C_{n},\Pi)\,\, \text{not}\newline \text{containing } V_{k}\}$ is an initial ideal of $I_{\Pi}$ with respect to an integral weight vector $\lambda \in \mathbb{N}^{n}$ that satisfies $b_{V_i}\cdot\lambda>0$, for all $i$ from $1$ to $k-1$.
One example of such a weight function is an integral solution $\tilde{\lambda}$ of the equation $\Lambda_{C_n}\tilde{\lambda} = y$, where $y$ is a positive integral multiple of $(1,1,\dots,1,-(n-1))\in \mathbb{Z}^{n}$ and $\Lambda_{C_n}$ is the Laplacian matrix of $C_n$. This weight function $\tilde{\lambda}$, up to scaling, is the potential vector $b_q$ in \cite{baker2013chip}, \cite[Section 6]{manjunathschwil}. 
Let $T^{s}$ be a spanning tree of $C_n$ that is rooted at the vertex $n$. The \emph{spanning tree order} $\rm{\tilde{rev}}$ on $R_n$ induced by $T^{s}$ is the reverse lexicographic order on $R_n$ induced by the ordering $<$ on the variables where $x_i< x_j$  if the vertex $j$ is a descendent of the vertex $i$ in $T^{s}$.
We consider a weighted monomial order $m_{\rm{rev}}$ obtained from the weight vector $\lambda$ followed by $\rm{\tilde{rev}}$ on $R_n$. In the following, we compute the initial terms of the binomials ${\mathbf x}^{S \rightarrow \overline{S}}-{\mathbf x}^{\overline{S} \rightarrow S}$ with respect to the monomial order $m_{\rm{rev}}$, where $S$ varies over all the parsets of $(C_n,\Pi)$.

Let $S$ be a parset of $(C_n,\Pi)$ that does not contain $V_k$. This implies that $S=\cup_{i\in \tilde{S}}V_i$ where $\tilde{S}\subseteq \{1,\dots,k-1\}$. Since $b_{V_i}\cdot\lambda>0$ for all $i$ from $1$ to $k-1$, this implies that $  ((S\rightarrow \overline{S})-(\overline{S}\rightarrow S))\cdot\lambda = (\sum_{i\in \tilde{S}}b_{V_i})\cdot\lambda>0$. Hence, ${\rm{in}}_{\lambda}({\bf x}^{S \rightarrow \overline{S}}-{\bf x}^{\overline{S} \rightarrow S}) ={\bf x}^{S \rightarrow \overline{S}} ={\rm{in}}_{m_{\rm{rev}}}({\bf x}^{S \rightarrow \overline{S}}-{\bf x}^{\overline{S} \rightarrow S})$, where $S$ is a parset of $(C_n,\Pi)$ not containing $V_k$.
The monomial order $m_{\rm{rev}}$ on $R_{n}=\mathbb{K}[x_1,\dots,x_n]$ is graded with respect to the weight vector $\lambda$ (\cite[Subsection 3.2.1]{herzog2011monomial}). From the proof of \cite[Theorem 3.1.2]{herzog2011monomial}, to prove that the set $\{\mathbf{x}^{S\rightarrow \overline{S}}\mid~S~ \text{is a connected-parset of}\, (C_{n},\Pi) \,\text{not containing}\,\, V_k\}$ is an initial ideal of $I_{\Pi}$ with respect to the weight vector $\lambda$, it is enough to show that it is an initial ideal of $I_{\Pi}$ with respect to the monomial order $m_{\rm{rev}}$.

\begin{theorem}\label{construction of parking ideal}
The set $\{\mathbf{x}^{S\rightarrow \overline{S}}-\mathbf{x}^{\overline{S}\rightarrow S}\mid ~S~\text{ is a connected-parset of }(C_{n},\Pi) \text{ not}\newline \text {containing } V_k \}$ forms a Gr\"obner basis of $I_{\Pi}$ with respect to the monomial order $m_{\rm{rev}}$ on $R_n$.
\end{theorem}
\begin{proof}
The set $L'=\{\sum_{i\in \tilde{S}}b_{V_i}\mid \tilde{S}\subset \{1,\dots,k\} \}$ is a generating set of the Laplacian lattice $L_{\Pi}$ of the parcycle $(C_n,\Pi)$, and it contains $L=\{\sum_{i\in \tilde{S}}b_{V_i}\mid \tilde{S}\subset \{1,\dots,k\},\, \text{the subgraph induced by } S=\cup_{i\in \tilde{S}}V_i~ \text{is connected}\}$.
This implies $I_L \subseteq I_{L'}\subseteq I_{\Pi}$ where $I_{L'}=\langle \mathbf{x}^{S\rightarrow \overline{S}}-\mathbf{x}^{\overline{S}\rightarrow S} \mid S=\cup_{i\in \tilde{S}}V_i\,\, \text{and}\,\,  \tilde{S}\subseteq \{1,\dots,k-1\}\rangle$. Recall that $I_{L}=\langle\mathbf{x}^{S\rightarrow \overline{S}}-\mathbf{x}^{\overline{S}\rightarrow S}\mid~S~ \text{ is a connected-parset of }\newline (C_n,\Pi) \,\text{not}~ \text{containing}\, V_k \rangle$. 

From Lemma  \ref{generating lemma}, $I_{\Pi}=I_L$ and hence, $I_{\Pi}=I_{L'}$. First, we prove that the generating set $ \{\mathbf{x}^{S\rightarrow \overline{S}}-\mathbf{x}^{\overline{S}\rightarrow S}\mid~S~ \text{ is a parset of }(C_{n},\Pi) \text{ not containing } \,V_k\}$ forms a Gr\"obner basis of $I_{\Pi}$ with respect to $m_{\rm{rev}}$. We prove this by using Buchberger's criterion (\cite[Theorem 2.3.2]{herzog2011monomial}). Let $S_1$ and $S_2$ be two parsets of $(C_n,\Pi)$ that do not contain $V_k$. Then $S_1\setminus S_2$ and $ S_2\setminus S_1$ (when non-empty) are also parsets of $(C_n,\Pi)$ that do not contain $V_k$. Let $B(S) ={\mathbf x}^{S \rightarrow \overline{S}}-{\mathbf x}^{\overline{S} \rightarrow S}$ and $\Lambda_{C_n}\rchi_{S}= ({S \rightarrow \overline{S}})-({\overline{S} \rightarrow S})$, where $S$ is a parset of $(C_n,\Pi)$ and $\Lambda_{C_n}$ is the Laplacian matrix of $C_n$. 
Now, we will show that the $S$-polynomial of $(B(S_1),B(S_2))$ reduces to $0$ (as per the division algorithm \cite[Theorem 2.2.1]{herzog2011monomial}) with respect to $B(S_1\setminus S_2)$ and $B(S_2\setminus S_1)$. 
The $S$-polynomial of $B(S_1)$ and $B(S_2)$ is 
$$S(B(S_1),B(S_2)) = m_2({\bf x}) {\bf{x}}^{\overline{S_{2}} \rightarrow S_{2}}-m_1({\bf{x}}) {\bf{x}}^{\overline{S_{1}} \rightarrow S_{1}},
$$
where $
m_1({\bf{x}}){\bf{x}}^{{S_{1}} \rightarrow {\overline S_{1}}}= {\rm{lcm}}({\bf{x}}^{{S_{1}} \rightarrow {\overline S_{1}}},{\bf{x}}^{{S_{2}} \rightarrow {\overline S_{2}}})
=m_2({\bf{x}}) {\bf{x}}^{{S_{2}} \rightarrow {\overline S_{2}}}$. Let $u\in \mathbb{N}^n $ be the exponent of the monomial ${\rm{lcm}}({\bf{x}}^{{S_{1}} \rightarrow {\overline S_{1}}},{\bf{x}}^{{S_{2}} \rightarrow {\overline S_{2}}})$. The exponents of the monomials $m_1({\bf{x}}) {\bf{x}}^{\overline{S_{1}} \rightarrow S_{1}}, m_2({\bf x}) {\bf{x}}^{\overline{S_{2}} \rightarrow S_{2}}$ are $u- \Lambda_{C_n}\rchi_{S_1}$ and $u- \Lambda_{C_n}\rchi_{S_2}$, respectively.
Without loss of generality, let us suppose that the initial term of $S(B(S_1),B(S_2))$ with respect to $m_{\rm{rev}}$ is $m_1({\bf{x}}) {\bf{x}}^{\overline{S_{1}} \rightarrow S_{1}}$. We divide the monomial $m_1({\bf{x}}) {\bf{x}}^{\overline{S_{1}} \rightarrow S_{1}}$ by $B(S_2\setminus S_1)$. The resulting remainder $r_{1}({\bf{x}})$ obtained after dividing the monomial $m_1({\bf{x}}) {\bf{x}}^{\overline{S_{1}} \rightarrow S_{1}}$ by $B(S_2\setminus S_1)$ has the exponent 
$u-\Lambda_{C_n}\rchi_{S_{1}}-\Lambda_{C_n}\rchi_{(S_{2}\setminus S_{1})}$. If $S_1\subset S_2$, then $u-\Lambda_{C_n}\rchi_{S_{1}}-\Lambda_{C_n}\rchi_{(S_{2}\setminus S_{1})}=u-\Lambda_{C_n}\rchi_{S_{2}}$. Hence, the remainder of $S(B(S_1),B(S_2))$ after dividing by $B(S_2\setminus S_1)$ becomes zero. Suppose that $S_1$ is not a subset $S_2$.
As $$(\Lambda_{C_n}\rchi_{S_{1}}+\Lambda_{C_n}\rchi_{(S_{2}\setminus S_{1})})\cdot\lambda > (\Lambda_{C_n}\rchi_{S_{2}})\cdot\lambda,$$ this implies $(u-\Lambda_{C_n}\rchi_{S_{1}}-\Lambda_{C_n}\rchi_{(S_{2}\setminus S_{1})})\cdot \lambda < (u-\Lambda_{C_n}\rchi_{S_{2}})\cdot \lambda$. Hence, the weight of the monomial $m_2({\bf{x}}) {\bf{x}}^{\overline{S_{2}} \rightarrow S_{2}}$ is greater than the weight of the monomial $r_1({\bf{x}}).$ Now, after dividing the monomial $m_2({\bf{x}}) {\bf{x}}^{\overline{S_{2}} \rightarrow S_{2}}$ by $B(S_1\setminus S_2)$, we get a remainder $r_2({\bf{x}})$ with the exponent $u-\Lambda_{C_n}\rchi_{S_{2}}-\Lambda_{C_n}\rchi_{(S_{1}\setminus S_{2})}$. From the chip-firing game on a graph, we know that $$\Lambda_{C_n}\rchi_{S_1}+\Lambda_{C_n}\rchi_{(S_2\setminus S_1)} = \Lambda_{C_n}\rchi_{(S_1 \cup S_2)} = \Lambda_{C_n}\rchi_{S_2}+\Lambda_{C_n}\rchi_{(S_1\setminus S_2)}.$$ Hence, the remainder monomials $r_1({\bf{x}}),r_2({\bf{x}}) $ are the same. Thus, the remainder of $S(B(S_1), B(S_2))$ after dividing by $B(S_1\setminus S_2),B(S_2\setminus S_1)$ is $r_{2}(\mathbf{x})-r_1(\mathbf{x})=0 $. By Buchberger's criterion, the set $\{{\bf x}^{S \rightarrow \overline{S}}-{\bf x}^{\overline{S} \rightarrow S}\mid~S~ \text{ is a parset of }\newline (C_n,\Pi) \,\text{not containing}\, V_k\}$ forms a Gr\"obner basis of $I_{\Pi}$.
Now, we reduce the size of the above Gr\"obner basis. We carry this out by showing that the initial term of a $B(S)$ corresponding to a parset of $(C_n,\Pi)$ is divisible by some $B(U)$, where $U$ is a connected-parset of $(C_n,\Pi)$ that does not contain $V_k$. If $\overline{S}$ is a connected-parset of $(C_n,\Pi)$, then we are done. Suppose that $\overline{S}$ is not a connected-parset of $(C_n,\Pi)$, then the subgraph induced on $C_n$ by $S$ is not connected. In this case, there exists a connected component $U$ of $S$ which is a connected-parset of $(C_n,\Pi)$. For every $i\in U$, ${\rm{out}}_{U}(i)=\sum_{j\in \overline{U}}a_{i,j}=\sum_{j\in \overline{S}}a_{i,j} ={\rm{out}}_{S}(i)$ where $(a_{i,j})$ is the adjacency matrix of $C_n$. This shows that ${\bf x}^{U \rightarrow \overline{U}}$ divides ${\bf x}^{S \rightarrow \overline{S}}$ and hence, the set $\{{\bf x}^{S \rightarrow \overline{S}}\mid~S~ \text{ is a connected-parset of }\,(C_n,\Pi) \,\text{not containing }\,V_{k}\}$ forms a generating set of ${\rm{in}}_{m_{\rm{rev}}}(I_{\Pi})$. As a consequence of it, the set $\{{\bf x}^{S \rightarrow \overline{S}}-{\bf x}^{\overline{S} \rightarrow S}\mid~S~ \text{ is a connected-parset of}~ (C_n,\Pi) \text{ not containing}\,V_k\}$ forms a Gr\"obner basis of $I_{\Pi}$ with respect to the monomial order $m_{\rm{rev}}$.\end{proof}
\begin{remark}\rm The reader may wonder if results on Gr\"obner bases of determinantal ideals available in literature can be used to prove Theorem \ref{construction of parking ideal}. These results usually assume that the underlying matrix is either a Hankel matrix or is an ordinary matrix of indeterminates. However, the matrix $A_{(C_n,\Pi)}$ associated with a parcycle $(C_n,\Pi)$ is, in general, neither of these. \qed
\end{remark}

Let $\tilde{E}(C_n)=\{(j,j+1),(n,1)\mid 1\leq j\leq n-1\}$ be the set of the direct edges of $C_n$ and $A_{(C_n,\Pi)}$ be the matrix associated with $(C_n,\Pi)$ obtained by fixing the main vertex $V_k$.  
\begin{corollary}\label{prime and cohen macaulayness of toppling ideal}
Let $\Pi=\{V_1,\dots,V_k\}$ be a connected $k$-partition of the vertex set of $C_n$. The toppling ideal $I_{\Pi}$ is generated by the set of all $2\times 2$ minors of $A_{(C_n,\Pi)}$. Moreover, $R_n/I_{\Pi}$ is a Cohen-Macaulay domain if $|V_i|\geq 2$ for some $i$.
\end{corollary}
\begin{proof}
From Lemma \ref{bijection lemma} and Theorem \ref{construction of parking ideal}, the toppling ideal $I_{\Pi}$ is generated by the set of all $2\times 2$ minors of the matrix $A_{(C_n,\Pi)}$. If $|V_i|\geq 2$ for some $i$, then $A_{(C_n,\Pi)}$ is one-generic (Lemma \ref{one genereic of A}).
Hence, $I_{\Pi}$ is a prime ideal and $R_n/I_{\Pi}$ is a Cohen-Macaulay domain (\cite[Theorem 6.4]{eisenbud2005geometry}).
\end{proof}
\begin{example}\label{binedid}\rm{
    Consider the parcycle $(C,\Pi)$ where $V(C)=\{u_1,\dots,u_{n},w_1,\newline \dots, w_n\}, E(C)=\{\{u_i,w_i\},\{w_{j},u_{j+1}\},\{w_n,u_1\}\mid 1\leq i\leq n, 1\leq j\leq n-1\}
    $ and $\Pi=\{\{w_n,u_1\},\{w_j,u_{j+1}\}\mid 1\leq j\leq n-1\}
    $. Let $\{w_n,u_1\}$ be the fixed main vertex (sink) of $(C,\Pi)$ and $u_1 w_1 u_2\dots w_n u_1$ be the directed closed path. From Corollary \ref{prime and cohen macaulayness of toppling ideal}, the toppling ideal $I_{\Pi}$ of $(C,\Pi)$ is the binomial edge ideal $J_{K_n}$ of $K_{n}$.}\qed
\end{example}

\begin{definition}\label{defGparid}
    Let $(C_n,\Pi)$ be a parcycle where $\Pi=\{V_1,\dots,V_k\}$ is a connected $k$-partition of the vertex set of $C_n$ and $n\in V_k$. The $G$-parking function ideal of $(C_n,\Pi)$ (with sink $V_k$) is the ideal generated by the set $\{\mathbf{x}^{S\rightarrow \overline{S}}\mid S\text{ is a connected-parset of } (C_{n},\Pi)\text{ not containing }\,V_k\}$. 
\end{definition}
From Theorem \ref{construction of parking ideal}, the $G$-parking function ideal of $(C_{n},\Pi)$ is an initial ideal of $I_{\Pi}$. In the case of the  connected $n$-partition $\Pi$ of $C_n$, the notion of the $G$-parking function ideal of $(C_n,\Pi)$ and $C_n$ are the same.
\begin{corollary}\label{GparantiD}
 The $G$-parking function ideal $M_{\Pi}$ of the parcycle $(C_n,\Pi)$ is generated by the set of all anti-diagonals of the $2\times 2$ minors of the matrix $A_{(C_n,\Pi)}$.
\end{corollary}
\begin{proof}
    The proof of the corollary is immediate from the definition of $M_{\Pi}$ and Proposition \ref{one generic matrix to parking ideal}.
\end{proof}

\subsection{Minimal Free Resolution of the $G$-Parking Function Ideal of a Parcycle}\label{resolution of G-parking}
In Subsection \ref{Gb and construction of G-parking}, we introduced an initial monomial ideal $M_{\Pi}=\langle {\mathbf x}^{S \rightarrow \overline{S}} \mid ~S ~\text{ is a}\newline \text{connected-parset of}\, (C_n,\Pi)\, \text{not containing}\, V_k\rangle$ of the toppling ideal $I_{\Pi}$.

In this subsection, we construct a complex $\mathcal{F}_{0,\Pi}$ of $R_{n}$-modules and prove that the complex $\mathcal{F}_{0,\Pi}$ minimally resolves $R_n/M_{\Pi}$. 
 The complex $\mathcal{F}_{0,\Pi}$ is defined as follows:
\begin{equation}\label{minresGpareq}
\mathcal{F}_{0,\Pi}: F_{0,\Pi, k-1} \xrightarrow{\delta_{0,k-1}} \dots \xrightarrow{\delta_{0,2}} F_ {0, \Pi,1} \xrightarrow{\delta_{0,1}}  F_{0,\Pi,0},
\end{equation}
where $F_{0,\Pi,j}=\oplus_{\mathcal{C}} R_n(-D(\mathcal{C}))$ and $\mathcal{C}$ varies over all $V_k$-acyclic $(j+1)$-partitions of $(C_{n},\Pi)$ which are identified with the standard basis elements of $F_{0,\Pi,j}$. Recall that the divisor $D(\mathcal{C})=\sum_{i\in C_n}{\rm{out}}_{\mathcal{C}}(i)[i]$. The differentials are defined as $$\delta_{0,j}(e_\mathcal{C}) = \sum_{f\in \mathcal{C}} {\rm sign}_{\mathcal{C}}(f) {\bf x}^{D(\mathcal{C})-D(\mathcal{C}/f)}\cdot e_{(\mathcal{C}/f)},$$
where the sum is taken over all the contractible edges $f$ of $\mathcal{C}$ and  $e_{\mathcal{C}}$ denotes the basis element of $R_n(-D(\mathcal{C}))$. 
\begin{example}\rm
    For the parcycle $C_{\Pi,\mathcal{M}}$ illustrated in Figure \ref{paroeg}, the complex $\mathcal{F}_{0,\Pi}$ is as follows:

$$\mathcal{F}_{0,\Pi}:\tilde{R}^3\overset{\delta_{0,3}}\longrightarrow \tilde{R}^8 \overset{\delta_{0,2}}\longrightarrow \tilde{R}^6 \overset{\delta_{0,1}} \longrightarrow \tilde{R}^1$$
    
where $\tilde{R}=\mathbb{K}[x_1,x_2,x_3,x_{v_0},x_{v_1},x_{s_1},x_{s_2}]$.
The matrices of the differentials are
$$\delta_{0,1}=
\begin{bmatrix}
    x_1 x_{v_0}& x_3^{2} &x_2 x_{v_1}& x_1 x_3&x_3 x_{v_1}&x_1 x_{v_1}
\end{bmatrix}
$$

    $\delta_{0,2}=
\begin{bmatrix}
    0&0&-x_{v_1}&0&0&0&-x_3&0\\
    0&0&0&0&-x_{v_1}&0&0&-x_1\\
    0&-x_1&0&0&0&-x_3&0&0\\
    -x_{v_1}&0&0&0&0&0&x_{v_0}&x_3\\
    0&0&0&-x_1&x_3&x_2&0&0\\
    x_3&x_2&x_{v_0}&x_3&0&0&0&0\\
\end{bmatrix}$
\\
\\

$
\delta_{0,3}=
\begin{bmatrix}
x_{v_0}&x_3&0\\
0&0&-x_3\\
-x_3&0&0\\
0&-x_3&x_2\\
0&-x_1&0\\
0&0&x_1\\
x_{v_1}&0&0\\
0&x_{v_1}&0\\
\end{bmatrix}
$

The basis elements of the free modules in $\mathcal{F}_{0,\Pi}$ correspond to six $\{s_1,s_2\}$-acyclic $2$-partitions, eight $\{s_1,s_2\}$-acyclic $3$-partitions, and three $\{s_1,s_2\}$-acyclic $4$-partitions, in the order from left to right are illustrated in Figure \ref{acyclic partitions}.\qed
\end{example}


We prove the exactness of the complex  $\mathcal{F}_{0,\Pi}$ by defining an isomorphism of complexes between $\mathcal{F}_{0,\Pi}$ and an exact cellular complex $\mathcal{H}_{\Pi}$ which is obtained from a certain hyperplane arrangement associated to $(C_{n},\Pi)$. Before moving to the construction of $\mathcal{H}_{\Pi}$, we study the graphical hyperplane arrangement of a parcycle.

\subsubsection{Graphical Hyperplane Arrangements of Parcycles}\label{grphyperofparsubsec}
Let $\Pi=\{V_1,\dots,V_k\}$ be a connected $k$-partition of $V(C_{n})$ and $n\in V_k$. Let $(C_{n},\Pi)$ be the associated parcycle. For every relevant edge $\{i,j\}$ of $(C_n,\Pi)$, we define a corresponding hyperplane $h_{ij}= \{ {\bf{p}}=(p_1,\dots,p_n) \in \mathbb{R}^n ~|~ p_{i} =p_{j}\}$ in $\mathbb{R}^n.$ We call the arrangement $\mathcal{A}_{\Pi} = \{h_{ij}~|~ \{i,j\}\, \text{is a relevant edge of } (C_n,\Pi)\}$ of hyperplanes in $\mathbb{R}^n$ as the \emph{graphical arrangement} of $(C_{n},\Pi)$. Consider the affine subspace $U_{\Pi}=\{ {\bf p} \in \mathbb{R}^n ~|~ p_n=0,p_{1}+\dots+p_{n-1} = 1, p_{u}=p_{w}~ ,\text{whenever}~  \{u,w\}\subseteq V_i~ \text{for some} ~i~ \,\text{from } $1$ \text{ to}\, k\}
$
of $\mathbb{R}^n.$ 
Let $\tilde{\mathcal{A}}_{\Pi} = \{h_{ij} \cap U_{\Pi} ~|~ \{i,j\} \,\text{is a relevant edge of } (C_n,\Pi)\}$ be the restriction of the arrangement $\mathcal{A}_{\Pi}$ to the affine subspace $U_{\Pi}$. Note that the arrangement $\tilde{\mathcal{A}}_{\Pi}$ is an \emph{essential hyperplane arrangement} in the sense of \cite[Subsection 1.1]{stanley2004introduction}.
In the case of the connected $n$-partition $\Pi$ of $V(C_n)$, the arrangements $\tilde{\mathcal{A}}_{\Pi}$ and $\tilde{\mathcal{A}}_{C_{n}}$ are the same where $\tilde{\mathcal{A}}_{C_{n}} = \{h_{ij} \cap U ~|~ \{i,j\} \in E(C_{n})\}$ and $U=\{ {\bf p} \in \mathbb{R}^n ~|~ p_n=0,p_{1}+\dots+p_{n-1} = 1\}$ (introduced in \cite[Section 6]{ds}). The arrangement $\tilde{\mathcal{A}}_{\Pi}$ divides the space $U_{\Pi}$ into polyhedra of various dimensions, called \emph{cells}. Let $B_{\Pi}$ be the collection of all bounded cells of $\tilde{\mathcal{A}}_{\Pi}$ on $U_{\Pi}$.  
Note that $ h_{ij} \cap U_{\Pi} = h_{ij} \cap U_{\Pi} \cap U = U_{\Pi} \cap (U \cap h_{ij}).$ This shows that the hyperplane arrangement induced on $U_{\Pi}$ by $\tilde{\mathcal{A}}_{C_{n}}$ is the same as $\tilde{\mathcal{A}}_{\Pi}$. Hence, each bounded cell of $\tilde{\mathcal{A}}_{\Pi}$ is the intersection of some bounded cell of $\tilde{\mathcal{A}}_{C_{n}}$ with $U_{\Pi}$. The collection $B_{\Pi}$ of all bounded cells in $\tilde{A}_{\Pi}$ forms a polyhedral complex in $U_{\Pi}$.
This implies that each $k$-dimensional bounded cell of $B_{\Pi}$ is also a $k$-dimensional bounded cell of $B_{C_{n}}$ where $B_{C_n}$ is the polyhedral complex of all the bounded cells in $\tilde{\mathcal{A}}_{C_{n}}$.


     We now assign monomials to the cells of $B_{\Pi}$. We start this by assigning monomials to the 0-cells of $B_{\Pi}$. From \cite[Proposition 3]{ds}, \cite[Corollary 2]{ds}, we know that $0$-cells of $B_{C_{n}}$ are of the form $u={\rchi}_{J}/|J|$ where ${\rchi}_{J} \in \{0,1\}^n$ is the characteristic vector of the non-empty subset $J \subseteq [1,\dots,n-1]$. Moreover, $u={\rchi}_{J}/|J|$ is a $0$-cell of $B_{C_{n}}$ if and only if the subgraphs induced on $C_{n}$ by $J$ and $\overline{J}$ are both connected. Since $0$-cells of $B_{\Pi}$ are $0$-cells of $B_{C_{n}}$ that lie on $U_{\Pi}$, every $0$-cell of $B_{\Pi}$ has the same value for all basic vertices lying in a main vertex of $(C_n,\Pi)$.
This implies that $u={\rchi}_{J}/|J|$ is a 0-cell of $B_{\Pi}$ if and only if $J\subseteq[1,\dots,n-1]$ is a connected-parset of $(C_{n},\Pi)$.
Corresponding to every $0$-cell $u$ of $B_{C_n}$, we assign a vector $a_{u} \in \mathbb{N}^{n-1}$ as defined in \cite[Section 6]{ds}. The vector $a_{u} \in \mathbb{N}^{n-1}$ is defined as $a_u(i) = d_J(i)$, where $d_{J}(i)$ is the number of edges in $C_{n}$ between $i$ and vertices in $\overline{J}$ if $i\in J$, and $d_J(i)=0$ if $i\in \overline{J}$. With the help of vectors assigned to $0$-cells of $B_{C_n}$, we assign a vector $a_{\tilde{B}}$  to every bounded cell $\tilde{B} \in B_{C_n}$ and it is defined as $a_{\tilde{B}}(i) = {\rm{max}}\{a_{u}(i)~|~ u \in \tilde{B}~ \text{is a 0-cell}\} $. As every vector $w \in
\mathbb{N}^{n-1}$ corresponds to a monomial $m_{w} = {\bf x}^{w} \in R_{n-1}$, this correspondence turns the polyhedral complex $B_{C_{n}}$ into a labelled polyhedral complex where each bounded cell $\tilde{B} \in B_{C_{n}}$ is labelled by the monomial $m_{\tilde{B}}= {\bf x}^{a_{\tilde{B}}} \in R_{n-1}$. Now, we assign monomials to 0-cells, in fact to all the bounded cells, of $B_{\Pi}$ in the same fashion as they were assigned in the case of $B_{C_n}$. This assignment of monomials to $0$-cells of $B_{\Pi}$ makes $B_{\Pi}$ a labelled polyhedral complex. The monomials assigned to $0$-cells of $B_{\Pi}$ are of the form ${\bf x}^{J \rightarrow \overline{J}}$, where $J$ is a connected-parset of $(C_n,\Pi)$ not containing $V_k$.
Let $\mathcal{H}_{\Pi}$ be the complex of $R$-modules obtained from the labelled polyhedral complex $B_{\Pi}$. We prove the exactness of $\mathcal{H}_{\Pi}$ by proving the the star-convexity of $B_{\Pi}$.

\begin{lemma}

    \label{star convexity theorem}
Let $(C_n,\Pi)$ be a parcyle and $B_{\Pi}$ be its labelled polyhedral complex obtained from the hyperplane arrangement $\tilde{\mathcal{A}}_{\Pi}$. The cellular free complex $\mathcal{H}_{\Pi}$ obtained from the labelled polyhedral complex $B_{\Pi}$ is a minimal free resolution of $R/M_{\Pi}$.
\end{lemma}
\begin{proof}
    We start with proving that the set $|B_{\Pi}|_{\leq \sigma} = \cup \{\tilde{B} \in B_{\Pi} ~|~ a_{\tilde{B}} \leq \sigma \}$ is star-convex for every $\sigma \in \mathbb{N}^{n-1}$. As a consequence of this, $|B_{\Pi}|_{\leq \sigma}$ is contractible and, in particular, $\mathbb{K}$-acyclic for any field $\mathbb{K}$. The proof consists of the following key steps:
\begin{enumerate}
\item \label{id_step} Identify a certain $0$-dimensional cell $q$ of $B_{C_n}$ as a candidate star-center for $|B_{\Pi}|_{\leq \sigma}$.

\item\label{qinbpi_step} Show that $q \in B_{\Pi}$.

\item \label{qinbpisig_step} Show that the label $a_{q}$ of $q$ 
satisfies $a_q \leq \sigma$ and hence, conclude that $q \in |B_{\Pi}|_{\leq \sigma}$.

\item\label{qstar_step} Prove that $q$ is indeed the star-center of $|B_{\Pi}|_{\leq \sigma}$.
\end{enumerate}

Steps \ref{id_step}, \ref{qinbpisig_step}
and \ref{qstar_step} closely follow \cite[Lemma 1]{ds}, the spirit of these arguments trace back to \cite{greene1983interpretation}. Our main contribution here is in Step \ref{qinbpi_step}. In the following, we provide the details of the proof. 

Let $\sigma \in \mathbb{N}^{n-1}$ and let $u_1, \dots,u_m $ be the $0$-dimensional cells in $|B_{\Pi}|_{\leq \sigma}$. Let $a_{i} = a_{u_i} \in \mathbb{N}^{n-1}$ be the vector corresponding to the $0$-dimensional cell $u_{i} \in B_{\Pi}$, for all $i$ from $1$ to $m$. Let $\tilde{J} = \cup_{i=1}^{m}{\rm{Supp}}(a_i)$ and $K$ be the connected component that contains $n$ in the subgraph (of $C_n$) induced by $V(C_n)\setminus \tilde{J}$. Consider $q := {\rchi}_{W}/|W|$ where $W = V(C_n)\setminus K$. As the subgraphs of $C_{n}$ induced by $W$ and $\overline{W} = K$ are connected, we deduce that $q = {\rchi}_{W}/|W|$ is a $0$-dimensional cell of $B_{C_n}$ (\cite[Corollary 2]{ds}).

Next, we will show that $q \in B_{\Pi}$. Before proving that $q \in B_{\Pi}$, we observe some relations between $W_{j}={\rm{Supp}} (u_{j}), \tilde{J}$ and $W$. Note that ${\rm{Supp}}(a_{j}) \subseteq {\rm{Supp}}(u_{j})$ for all $j$ from $1$ to $m$  and hence, $\tilde{J} \subseteq \cup_{j=1}^{m}W_{j}$. Now, we prove that $W_{j} \subseteq W$ for all $j$. Let $i\in W_{j}$ and $d_{W_{j}}(i) > 0$, then $i \in {\rm{Supp}}(a_{j}) \subseteq \tilde{J} \subseteq W$. Suppose $d_{W_{j}}(i) = 0$, this implies that $i \notin {\rm{Supp}}(a_{j})$. Hence, every path joining $i$ to $n$ contains an edge $\{s,t\}$ such that $s \in {\rm{Supp}}(a_{j})$ and $t \notin {\rm{Supp}}(a_j)$. This implies that $i$ does not lie in the connected component of the subgraph (of $C_{n}$) induced by $V(C_n) \setminus {\rm{Supp}}(a_{j})$ that contains $n$. Since ${\rm{Supp}}(a_{j}) \subseteq \tilde{J}$, this shows that $i \notin K$. Hence, $i \in W$ and $W_{j} \subseteq W$ for all $j$. As $q \in U=\{ {\bf p} \in \mathbb{R}^n ~|~ p_n=0,p_{1}+\dots+p_{n-1} = 1\}$, to prove that $q \in U_{\Pi}$ it is enough to show that for each $i$ either $V_{i} \subseteq W$ or $V_{i} \subseteq K$. Suppose that there  exists a main vertex $V_{i}\in \Pi$ such that $V_{i} \nsubseteq W$ and $V_{i} \nsubseteq K$. In this case, there exists vertices $v_{i_1},v_{i_2} \in V_{i}$ such that $\{v_{i_1},v_{i_2}\} \in E(C_{n})$, and $v_{i_1} \in W$ and $v_{i_2} \in K$. This implies that $d_{W}(v_{i_1})>0$ and hence, $v_{i_1} \in \tilde{J}$. This holds because if $v_{i_1} \notin \tilde{J}$, then $\{v_{i_1},v_{i_2}\} \in E(C_{n})$ and  $v_{i_2} \in K$ implies that $v_{i_1} \in K$, which is a contradiction. Since $v_{i_1} \in \tilde{J} \subseteq \cup_{j=1}^{m} W_{j}$ and $v_{i_2}\in K$, there exists some $j$ such that $v_{i_1} \in W_{j}$ and $v_{i_2} \notin W_{j}$. This gives us a contradiction since $v_{i_1},v_{i_2} \in V_{i}$ and $W_j$ is a connected-parset of $(C_n,\Pi)$. Hence, for every $j$, either $V_j \subseteq W$ or $V_j \subseteq K$. This shows that $q \in U_{\Pi}$ and is a $0$-dimensional cell of $B_{\Pi}$.

Let $a_{q}\in \mathbb{N}^{n-1}$ be the vector corresponding to $q$. Now, we prove that $a_{q}\leq \sigma$. Suppose $a_{q}(i) = d_{W}(i)>0$, then $i \in \tilde{J} (\subseteq\cup_{l=1}^{m} W_{l})$. Hence, there exists an $W_{l}$ such that $i \in W_{l}(\subseteq W)$ and $d_{W}(i) \leq d_{W_{l}}(i)$.
This shows that, for each $i$ there exits an $a_{l}$ such that $a_{q}(i) \leq a_{l}(i) \leq \sigma(i)$.

We now prove that $q$ is a star-center of $|B_{\Pi}|_{\leq \sigma}$. Let $t \in |B_{\Pi}|_{\leq \sigma}$ be an arbitrary point. 
We prove that the line segment $[t,q]$ lies in a bounded cell of $\tilde{\mathcal{A}}_{\Pi}$. If $t({i})>t({j})$ holds for some relevant edge $\{i,j\}$ of $(C_{n},\Pi)$, then  $(a_{B_t})(i)>0$ where $B_t \in |B_{\Pi}|_{\leq \sigma}$ is the inclusion minimal cell containing $t$ {(\cite[Proposition 4]{ds})}. As ${\rm{Supp}}(a_{B_{t}}) \subseteq \tilde{J} \subseteq W$, this implies $q(i) = 1/|W|$ and hence, $q(i) \geq q(s)$ for all $s \in V(C_n)$. Therefore, there is no hyperplane $h_{ij}$ that separates $t$ and $q$. Hence, the open line segment $(t,q)$ is contained in a cell of $\tilde{A}_{\Pi}$. 
Let $B'$ be the inclusion minimal cell in $\tilde{A}_{\Pi}$ that contains $(t,q)$. We now show that $B^{'} \in B_{\Pi}$.  To prove that $B^{'} \in B_{\Pi}$, it is enough to show that there exists a bounded cell in $\tilde{\mathcal{A}}_{C_n}$ 
that contains the line segment $[t,q]$. Let $\tilde{B}$ be the inclusion minimal cell in $\tilde{\mathcal{A}}_{C_n}$ that contains $[t,q]$. Let $g$ be a point in the open line segment $(t,q) \subseteq {\rm{relint}} (\tilde{B})$. From \cite[Proposition 3]{ds}, to prove that $\tilde{B} \in B_{C_n}$, it is enough to show that the acyclic graph $C_{n}/g$ has a unique sink at the vertex class containing $n$. Suppose that $i$ lies in different vertex class that is another sink in $C_n/g$ different from the vertex class containing $n$. As $t \in |B_{\Pi}|_{\leq \sigma}$ lies in a bounded cell,  there exists a path $i=i_0i_1\cdots i_{m}=n$ such that $t({i_{r-1}}) \geq t({i_{r}})$, for all $1\leq r \leq m$ (\cite[Proposition 3]{ds}). 
 By assumption, this path is not weakly decreasing for $g$. Hence, there exists an index $j$ with $g({i_{j-1}})< g({i_{j}})$. This implies that $g({i_j})>0$. As ${\rm{Supp}}(g) = {\rm{Supp}}(t) \cup {\rm{Supp}}(q)$ and by construction ${\rm{Supp}}(t) \subseteq {\rm{Supp}}(q)$, this implies $q({i_j}) = 1/|W|$. Note that the induced acyclic orientation on $C_n/g$ has its sink at the vertex class containing $i$ and ${\rm{Supp}}(g)= {\rm{Supp}}(q)$. Thus, the path $i=i_0i_1\cdots i_{m}=n$ is weakly decreasing for $q$ also. Furthermore, this implies that $g \in (t,q) \subset \{{\bf{p}} \in U: p_{i_{j-1}} \geq p_{i_j}\}$, which is a contradiction as $g({i_{j-1}})<g({i_j})$. Hence, $\tilde{B}$ is a bounded cell in $B_{C_n}$. Since $g \in U_{\Pi}$ and $\tilde{B}$ is the inclusion minimal bounded cell of $B_{C_n}$, we deduce that  $\tilde{B}$ lies in $U_{\Pi}$. 
  Now, we are only left with proving that $a_{\tilde{B}} \leq \sigma.$ From \cite[Proposition 4]{ds}, we know that $(a_{\tilde{B}})(i) = \#\{\{i,j\} \in E(C_{n})~|~ g(i) > g(j)\}$. Since $g\in (t,q)$, note that $g(i)>g(j)$ if and only if $ wt({i})+(1-w)q({i})  > wt({j})+(1-w)q({j})$, where $w \in (0,1)$ and $g(k) = wt({k})+(1-w)q({k})$ for all $k$.
 If $g(i) = 0$, then $(a_{\tilde{B}})(i)=0$. Now, assume that $g(i) >0$. Since ${\rm{Supp}}(g)= {\rm{Supp}}(q)$, this implies $i \in W$. If $t(i) = 0$, then $g(i)>g(j)$ implies $q(i) > q(j)$. Hence, $(a_{\tilde{B}})(i) \leq (a_{q})(i)\leq \sigma(i)$. 
  Now, assume that $t({i})>0$. Since ${\rm{Supp}}(t)\subseteq {\rm{Supp}}(q)$, this implies $i\in I$ and $q(i)>0$. Moreover, $q(i)>q(j)$ implies that $t(i)>t(j)$ if $t(i)>0$. 
  In this case, $g(i)>g(j)$ implies that $t({i})>t({j})$. Hence, $(a_{\tilde{B}})(i) \leq (a_{B_{t}})(i)\leq \sigma(i)$. This proves that the complex $\mathcal{H}_{\Pi}$ is exact and 
its minimality follows from the isomorphism (defined in Theorem \ref{complex of G-parking ideal}) between $\mathcal{H}_{\Pi}$ and $\mathcal{F}_{0,\Pi}$.
  \end{proof}

  In the following, we show that  $\mathcal{F}_{0,\Pi}$ is a minimal free resolution of $R_n/M_{\Pi}$. 
First, we define a homomorphism of complexes $\phi$ between the free resolution $\mathcal{F}_{0}(C_n)$ and the cellular resolution $\mathcal{H}_{C_n}$ obtained from $B_{C_n}$. The cellular resolution $$\mathcal{H}_{C_n} :\begin{tikzcd}[cells={nodes={minimum height=2em}}] 0\arrow[r] & H_{n-1} \arrow[r,"\rho"]  &   H_{n-2} \arrow[r,"\rho"]  &  \dots \arrow[r,"\rho"]& H_{1}\arrow[r,"\rho"] & R\arrow[r,"\rho"]& 0
\end{tikzcd}
$$
where $H_{i} = \bigoplus_B {R_{n}}{e_B}$, the direct sum is taken over all the $i-1$ dimensional bounded cells $B$ in $B_{C_n}$ and $R_n{e_B}$ is a free $R_n$-module of rank $1$ for all $B$.

Let $h_{ij}= \{\mathbf{p}=(p_1,\dots,p_n) \in \mathbb{R}^n : p_i = p_j \}$ be the hyperplane corresponding to the edge $\{i,j\} \in E(C_n)$ and let $h^{\geq}_{i,j}= \{{\bf {p}} \in \mathbb{R}^n : p_i \geq p_j \}$. Let ${\Pi}$ 
be a connected $k$-partition of $C_{n}$ and let $\mathcal{C}$ be an $n$-acyclic $k$-partition with the vertex set $\Pi$. Let
$U_{\Pi}=\{ {\bf p} \in \mathbb{R}^n ~|~ p_n=0,p_{1}+\dots+p_{n-1} = 1, p_{k}=p_{l}~ \text{whenever}~  \{k,l\}\subseteq V_i~ \text{for some}~i\}
$
be the affine subspace of $\mathbb{R}^n$ corresponding to $\Pi$.
We define $\phi(e_{\mathcal{C}}) =e_{B(\mathcal{C})} $ where $B(\mathcal{C})$ is the highest dimensional cell in $\tilde{\mathcal{A}}_{C_{n}}$ that lies in $U_{\Pi}$ and is bounded by the half-spaces $\{U_{\Pi} \cap h^{\geq}_{i,j}~|~(V_{j_1},V_{j_2}) ~\text{is a directed edge in} ~\mathcal{C}, i\in V_{j_1}, j\in V_{j_2} ~\text{and}~ \{i,j\}\in E(C_n)\}$.

\begin{proposition}
The map $\phi$ between $\mathcal{F}_0(C_n)$ and $\mathcal{H}_{C_n}$ is well-defined and is an isomorphism.
    \end{proposition}
 \begin{proof}
    Let $\mathcal{C}$ be an $n$-acyclic $k$-partition of $C_n$ with the vertex set $\tilde{\Pi}$. The image of $e_{\mathcal{C}}$ under $\phi$ is $\phi(e_{\mathcal{C}}) =e_{B(\mathcal{C})} $ where $B(\mathcal{C})$ is the highest dimensional cell in $\tilde{\mathcal{A}}_{C_{n}}$ that lies in $U_{\tilde{\Pi}}$ and is bounded by the half-spaces $\{U_{\tilde{\Pi}} \cap h^{\geq}_{i,j }~|~(V_{j_1},V_{j_2}) ~\text{is a directed edge in} ~\mathcal{C}, i\in V_{j_1}, j\in V_{j_2} ~\text{and}~ \{i,j\}\in E(C_n)\}$. 
  For any interior point ${\bf p}$ of $B(\mathcal{C})$, the acyclic graph  $C_n/{\bf p}=\mathcal{C}$. Hence, $B(\mathcal{C}) \in B_{C_n}$ from \cite[Proposition 3]{ds}, and ${\rm{dim}} (B(\mathcal{C})) = {\rm{dim}} (U_{\tilde{\Pi}})=k-2$ implies that $e_{B(\mathcal{C})}\in H_{k-1}$. 
  As two different $n$-acyclic partitions of $C_n$ differ by at least one directed edge, the map $\phi$ takes distinct $n$-acyclic partitions of $C_n$ to distinct cells of $B_{C_n}$.  The set of $n$-acyclic $k$-partitions of $C_n$ and the set of $k-2$ dimensional bounded cells of $B_{C_n}$ have the same number of elements, for all $2 \leq k \leq n$.
Hence, $\phi$ is a well-defined $R_n$-module isomorphism from $\mathcal{F}_{0,k-1}$ onto $H_{k-1}$ for all $k\geq 1$.

Consider the following diagram:
$$
\begin{tikzcd}[cells={nodes={minimum height=2em}}]
0\arrow[r] & F_{0,n-1} \arrow[r,"\delta"] \arrow[d,"\phi"] &   F_{0,n-2} \arrow[r,"\delta"] \arrow[d,"\phi"] &  \dots \arrow[r] & F_{0,1}\arrow[r,"\delta"]\arrow[d,"\phi"] & R\arrow[r,"\delta"]\arrow[d,"Id"]& 0\\
0\arrow[r] & H_{n-1} \arrow[r,"\rho"]  &   H_{n-2} \arrow[r,"\rho"]  &  \dots \arrow[r]&  H_{1}\arrow[r,"\rho"] & R\arrow[r,"\rho"]& 0
\end{tikzcd}.
$$
To prove that $\phi$ is a complex homomorphism, it is enough to show that $\phi \delta(e_{\mathcal{C}}) = \rho \phi(e_{\mathcal{C}})$ where $\mathcal{C} $ is an $n$-acyclic partition of $C_n$. Let $\mathcal{C}$ be an $n$-acyclic $k$-partition of $C_n$ with the associated connected $k$-partition $\tilde{\Pi} =\{V_1,\dots,V_k\}$. The image of $e_{\mathcal{C}}$ under $\phi$ is $e_{B(\mathcal{C})}$. Before proving that $\phi \delta(e_{\mathcal{C}}) = \rho \phi(e_{\mathcal{C}})$, we prove that the $0$-dimensional cells of $B(\mathcal{C})$ are of the form ${\rchi}_{I}/|I|$ where $I \rightarrow \overline{I}$ is an $n$-acyclic $2$-partition of $C_n$ obtained from $\mathcal{C}$. Every $0$-dimensional cell of $B_{C_n}$ that is of the form ${\rchi}_{I}/|I|$ where $I \rightarrow \overline{I}$ is an $n$-acyclic $2$-partition obtained from $\mathcal{C}$ lies in the half-spaces $\{U_{\tilde{\Pi}} \cap h^{\geq}_{i,j}~|(V_i,V_j) ~\text{is a directed edge in} ~\mathcal{C},\newline i\in V_i~, j\in V_j ~\text{and}~ \{i,j\}\in E(C_n)\}$. Hence, such $0$-dimensional cells of $B_{C_n}$ are part of the $0$-dimensional cells of $B(\mathcal{C})$.

Suppose there exists a $0$-dimensional cell ${\rchi}_{J}/|J|$ of $B(\mathcal{C})$ such that $J \rightarrow \overline{J}$ cannot be obtained from the edge contractions of $\mathcal{C}$. As $B(\mathcal{C})$ is a bounded cell in $U_{\tilde{\Pi}}$, the set $J$ is a connected-parset of $\Pi$ (discussed in the proof of Lemma \ref{star convexity theorem}). Suppose $J  =V_{j_1}\cup\dots\cup V_{j_l}$ and $V_{j_1} = \{i_{j_{1}},\dots,i_{j_{2}}-1\}$. Since $J \rightarrow \overline{J}$ cannot be obtained from the edge contractions of $\mathcal{C}$, the source of $\mathcal{C}$ cannot be in $J$. Without loss of generality, assume that $(V_{j_{i}},V_{j_{i+1}})$ are directed edges of $\mathcal{C}$ for all $i$. As $J$ does not contain the source vertex of $\mathcal{C}$, there exists a vertex $V_{j_{0}} = \{i_{j_0},\dots,i_{j_1}-1\} \in \Pi$ of $\mathcal{C}$ such that $(V_{j_{0}},V_{j_{1}})$ is a directed edge of $\mathcal{C}$ and is contractible. This implies that the $0$-dimensional cell ${\rchi}_{J}/|J|$ lies in the half-space $U_{\tilde{\Pi}} \cap h^{\geq}_{i_{j_1}-1,i_{j_1}}$, which gives a contradiction as $({\rchi}_{J}/|J|)(i_{j_1})> ({\rchi}_{J}/|J|)(i_{j_1}-1) = 0$. Hence, every $0$-dimensional cell of $B(\mathcal{C})$ is of the form ${\rchi}_{I}/|I|$ where $I \rightarrow \overline{I}$ is an $n$-acyclic $2$-partition obtained from $\mathcal{C}$. This implies that the $k-3$ dimensional bounded cells in $B(\mathcal{C})$ are of the form $B(\mathcal{C}/f)$ where $f$ is a contractible edge of $\mathcal{C}$.

Let $f =(V_j,V_{j+1})$ be a contractible edge of $\mathcal{C}$. 
The image $$\phi\delta(e_{\mathcal{C}}) =\phi\delta_{0,k-1}(e_{\mathcal{C}})= \sum_{f}{\rm sign}_{\mathcal{C}}(f)m_{\mathcal{C}}(f)\phi(e_{(\mathcal{C}/f)}),$$ where the sum is taken over all the contractible edges $f$ of $\mathcal{C}$. From the above description of the $0$-cells of $B(\mathcal{C})$, the coefficient of $B(\mathcal{C}/f)$ in $\rho(\phi(e_{\mathcal{C}}))$ is $m_{B(\mathcal{C})}/m_{B((\mathcal{C}/f))} = m_{\mathcal{C}}(f) $, where $m_{B(\mathcal{C})}$ and $m_{B((\mathcal{C}/f))}$ are monomials corresponding to the cells $B(\mathcal{C})$ and $ B(\mathcal{C}/f)$, respectively. 
This implies $$\rho(\phi(e_{\mathcal{C}}))=\sum_{f}{\varepsilon(B(\mathcal{C}),B(\mathcal{C}/f))m_{B(\mathcal{C})}/m_{B(\mathcal{C}/f)}\cdot e_{B(\mathcal{C}/f)}},$$ where $\varepsilon(B(\mathcal{C}),B(\mathcal{C}/f))\in \{0,1, -1\}$ is the incidence function described in \cite[Section 1]{bayer1998cellular} and the sum is taken over all the contractible edges $f$ of $\mathcal{C}$. From the construction of ${\rm{sign}}_{\mathcal{C}}(f)$ (\cite[Proposition 2.11]{manjunathschwil}), we can make the ${\rm{sign}}$ functions ${\rm{sign}}_{\mathcal{C}}(f)$ compatible with the incidence functions $\varepsilon(B(\mathcal{C}),B(\mathcal{C}/f))$ for all $\mathcal{C}$ and $f$. As a result the commutativity of the diagram holds and it shows that $\phi$ is a complex homomorphism between $\mathcal{F}_{0}(C_n)$ and $\mathcal{H}_{C_{n}}$. Since $\phi$ is also a bijective homomorphism between the free modules $F_{0,i} ,H_{i}$ for all $i$, the map $\phi$ is a complex isomorphism.
\end{proof}

    \begin{theorem}\label{complex of G-parking ideal}
The complex $\mathcal{F}_{0,\Pi}$ of $R_{n}$-modules minimally resolves $R_n/M_{\Pi}$.
\end{theorem}
\begin{proof}
We view the complex $\mathcal{F}_{0,\Pi}$ as a subcomplex of $\mathcal{F}_0(C_n)$ and $\mathcal{H}_{\Pi}$ as a subcomplex of $\mathcal{H}_{C_n}$, where $\mathcal{H}_{\Pi}$ and $ \mathcal{H}_{C_n}$ are the cellular resolutions obtained from $B_{\Pi}$ and $ B_{C_n}$, respectively. We restrict the complex map $\phi:\mathcal{F}_0(C_n) \longrightarrow \mathcal{H}_{C_n}$  to the subcomplex  $\mathcal{F}_{0,\Pi}$. The restriction of $\phi$ on $\mathcal{F}_{0,\Pi}$ maps distinct basis elements of the free module $F_{0,\Pi,i}$ to distinct basis elements of $H_{\Pi,i}$, where $F_{0,\Pi,i}$ and $H_{\Pi,i}$ are the $i$-th homological degree modules of $\mathcal{F}_{0,\Pi}$ and $\mathcal{H}_{\Pi}$, respectively. Hence, $\phi$ is a homomorphism of complexes from the complex $\mathcal{F}_{0,\Pi}$ to the complex $\mathcal{H}_{\Pi}$. To prove that $\phi$ is an isomorphism of complexes, it is enough to prove that $\phi_{i}:F_{0,\Pi,i} \rightarrow H_{\Pi,i}$ is a surjective map for all $i$.
Suppose that there exists an $(i-1)$ dimensional bounded cell $B$ in $H_{\Pi,i}$ that is not the image of any $V_k$-acyclic $(i+1)$-partition of $(C_{n},\Pi)$. 
From the surjectivity of $\phi_{i}:F_{0,i} \longrightarrow H_i$, there exists an $n$-acyclic $(i+1)$-partition $\mathcal{C}'$ of $C_n$ such that $\phi_{i}(e_{\mathcal{C}'}) =e_{B}$ and $\mathcal{C}'$ cannot be obtained from any $V_k$-acyclic $k$-partition of $(C_n,\Pi)$ by edge contractions. This implies that there exists a vertex $i\in V(C_n)$ such that $\{i,i+1\} \subseteq V_j$ for some $j$, but $\{i,i+1\}$ is not contained in any vertex of $\mathcal{C}^{'}$, i.e. if $\mathcal{C}^{'}=(\tilde{\Pi},\tilde{\mathcal{A}})$, then there exists $V_{j}^{'},V_{j+1}^{'} \in \tilde{\Pi}$ such that $i \in V_{j}^{'},i+1 \in V_{j+1}^{'}$ for some $j$. Without loss of generality, we assume that in the acyclic orientation $\tilde{\mathcal{A}}$ of  $\mathcal{C}^{'}$, there exists an oriented edge  $V_{j+1}^{'} \longrightarrow V_{j}^{'}$. This implies that there exists a $0$-dimensional cell $v={\rchi}_{I}/|I|$ of $B =\phi_{i}(e_{\mathcal{C}^{'}})$ such that
$v(i+1) \neq 0$ but $v(i)=0$. This gives us a contradiction as $B$ lies in $U_{\Pi}$. Hence, $\phi_{i}$ is an isomorphism from the free module $F_{0,\Pi,i}$ onto $H_{\Pi,i}$, and $\phi $ is an isomorphism of complexes from the complex $\mathcal{F}_{0,\Pi}$ to the complex $\mathcal{H}_{\Pi}$. As $\mathcal{H}_{\Pi}$ is an exact complex, we infer that $\mathcal{F}_{0,\Pi}$ is also an exact complex. Since the elements of the differential matrices of $\mathcal{F}_{0,\Pi}$ are contained in $\langle x_1,\dots,x_n\rangle$, the complex $\mathcal{F}_{0,\Pi}$ is a minimal free resolution of $M_{\Pi}$.
\end{proof}

\subsection{Minimal Free Resolution of the Toppling Ideal of a Parcycle}\label{resolution of toppling ideal}
Consider the sequence $\mathcal{F}_{1,\Pi}$ that is defined as follows 
\begin{equation}\label{reseq}
 \mathcal{F}_{1,\Pi} : F_{1,\Pi,k-1} \overset{\delta_{1,k-1}}\longrightarrow \dots \overset{\delta_{1,2}} \longrightarrow F_{1,\Pi,1} \overset{\delta_{1,1}} \longrightarrow F_{1,\Pi,0}, 
 \end{equation}
 where $F_{1,\Pi,j} = \bigoplus_\mathbf{c} R_n(-D(\mathbf{c}))$ and the direct sum is taken over all chip-firing equivalence classes of acyclic $(j+1)$-partitions of $(C_{n},\Pi)$, which are identified with the standard basis elements of $F_{1,\Pi,j}$.
 The free module $R_n(-D(\mathbf{c}))$ denotes the twist of $R_{n}$ by $-D(\mathbf{c})$, where $D(\mathbf{c})\in \mathbb{Z}^n/L_{\Pi}$ denotes the chip-firing equivalence class of divisors corresponding to the elements of $\mathbf{c}$. The differential $\delta_{1,j}  : F_{1,\Pi,j} \longrightarrow F_{1,\Pi,j-1}$ of $\mathcal{F}_{1,\Pi}$ is defined as follows: 
$$\delta_{1,j}(e_\mathbf{c}) = \sum_f{{\rm sign}_{\mathbf{c}}(f) m_{\mathbf{c}}(f)\cdot e_{({\mathbf{c}}/f)}},
 $$
 where the sum is taken over all the contractible edges $f$ of $\mathbf{c}$ and  $e_{\mathbf{c}}$ denotes the basis element of $R_n(-D(\mathbf{c}))$.

Recall that every connected $j$-partition of $(C_{n},\Pi)$ is a connected $j$-partition of $C_n$ that coarsens $\Pi$. Hence, every basis element of $F_{1,\Pi,j}$ corresponds to a basis element of $F_{1,j}$ where $F_{1,j}$ is the $j$-th homological degree free module of $\mathcal{F}_{1}(C_n)$ (defined in \cite[Section 3]{manjunathschwil}). The images of the differentials of $\mathcal{F}_{1}(C_n)$ and $\mathcal{F}_{1,\Pi}$ are the same on the basis elements of $F_{1,\Pi,j}$ for all $j$. This makes the complex $\mathcal{F}_{1,\Pi}$ a subcomplex of $\mathcal{F}_{1}(C_n)$. Hence, $\mathcal{F}_{1,\Pi}$ is a complex of $R_{n}$-modules and is $\mathbb{Z}^n/L_{\Pi}$-graded.

\begin{remark}\label{remsubcom}{\rm For a sub-parcycle $(C_n,\tilde{\Pi})$ of $(C_n,\Pi)$, i.e. $\tilde{\Pi}$ is a connected partition of the vertex set of $C_n$ that coarsens $\Pi$, the complex $\mathcal{F}_{1,\tilde{\Pi}}$ is a subcomplex of $\mathcal{F}_{1,\Pi}$.}\qed\end{remark}

 \begin{example}\label{homorunex}\rm
    For the parcycle $C_{\Pi,\mathcal{M}}$ illustrated in Figure \ref{paroeg}, the complex $\mathcal{F}_{1,\Pi}$ is as follows:

    $$\mathcal{F}_{1,\Pi}:\tilde{R}^3\overset{\delta_{1,3}}\longrightarrow \tilde{R}^8 \overset{\delta_{1,2}}\longrightarrow \tilde{R}^6 \overset{\delta_{1,1}} \longrightarrow \tilde{R}^1
$$
where $\tilde{R}=\mathbb{K}[x_1,x_2,x_3,x_{v_0},x_{v_1},x_{s_1},x_{s_2}]$. The matrices of the differentials are
\[
\delta_{1,1} = 
\resizebox{\textwidth}{!}{$
\begin{bmatrix}
    x_1 x_{v_0} - x_3 x_{s_2} & x_3^{2} - x_2 x_{v_0} & x_2 x_{v_1} - x_3 x_{s_1} & 
    x_1 x_3 - x_2 x_{s_2} & x_3 x_{v_1} - x_{v_0}x_{s_1} & x_1 x_{v_1} - x_{s_1} x_{s_2}
\end{bmatrix}
$}
\]

    $\delta_{1,2}=
\begin{bmatrix}
    0&0&-x_{v_1}&-x_{s_1}&0&0&-x_3&-x_2\\
    0&0&0&0&-x_{v_1}&-x_{s_1}&-x_{s_2}&-x_1\\
    -x_{s_2}&-x_1&0&0&-x_{v_0}&-x_3&0&0\\
    -x_{v_1}&-x_{s_1}&0&0&0&0&x_{v_0}&x_3\\
    0&0&-x_{s_2}&-x_1&x_3&x_2&0&0\\
    x_3&x_2&x_{v_0}&x_3&0&0&0&0\\
\end{bmatrix}$
\\
\\

$
\delta_{1,3}=
\begin{bmatrix}
x_{v_0}&x_3&0\\
0&x_{v_0}&-x_3\\
-x_3&-x_2&0\\
0&-x_3&x_2\\
-x_{s_2}&-x_1&0\\
0&-x_{s_2}&x_1\\
x_{v_1}&x_{s_1}&0\\
0&x_{v_1}&-x_{s_1}\\
\end{bmatrix}
$

The basis elements of the free modules in $\mathcal{F}_{1,\Pi}$ correspond to six chip-firing equivalence classes of acyclic $2$-partitions, eight chip-firing equivalence classes of acyclic $3$-partitions, and three chip-firing equivalence classes of acyclic $4$-partitions, in the order (from left to right) illustrated in Figure \ref{acyclic partitions}.
  
    In the following, we describe the complex $\mathcal{F}_{t,\Pi}$ (a one-parameter family of complexes of free $R_n$-modules).  For the parcycle $C_{\Pi,\mathcal{M}}$ illustrated in Figure \ref{paroeg}, consider the integral weight vector $\lambda_{\mathcal{M}}=(\lambda(x_1), \lambda(x_{v_0}),\lambda(x_3),\lambda(x_2),\lambda(x_{v_1}),\newline \lambda(x_{s_1}),\lambda(x_{s_{2}}))=(1,7,5,2,11,7,2)$ (defined in Subsection \ref{weights}) with the indeterminate $t$ having weight $1$. The complex $\mathcal{F}_{t,\Pi}$ is as follows:

    $$\mathcal{F}_{t,\Pi}:\tilde{R}[t]^3\overset{\delta_{t,3}}\longrightarrow \tilde{R}[t]^8 \overset{\delta_{t,2}}\longrightarrow \tilde{R}[t]^6 \overset{\delta_{t,1}} \longrightarrow \tilde{R}[t]^1
$$
The matrices of the differentials are 
\[
\delta_{t,1} = 
\resizebox{\textwidth}{!}{$
\begin{bmatrix}
    x_1 x_{v_0}-x_3 x_{s_2} t& x_3^{2}-x_2 x_{v_0} t &x_2 x_{v_1}-x_3 x_{s_1} t & x_1 x_3-x_2 x_{s_2} t^2 &x_3 x_{v_1}-x_{v_0}x_{s_1} t^2&x_1 x_{v_1}-x_{s_1} x_{s_2} t^3
\end{bmatrix}
$}
\]\\
 $\delta_{t,2}=
\begin{bmatrix}
    0&0&-x_{v_1}&-x_{s_1} t^2&0&0&-x_3&-x_2 t\\
    0&0&0&0&-x_{v_1}&-x_{s_1} t&-x_{s_2} t&-x_1\\
    -x_{s_2} t^2&-x_1&0&0&-x_{v_0} t&-x_3&0&0\\
    -x_{v_1}&-x_{s_1} t&0&0&0&0&x_{v_0}&x_3\\
    0&0&-x_{s_2} t&-x_1&x_3&x_2&0&0\\
    x_3&x_2&x_{v_0}&x_3&0&0&0&0\\
\end{bmatrix}$
\\
\\

$
\delta_{t,3}=
\begin{bmatrix}
x_{v_0}&x_3&0\\
0&x_{v_0} t&-x_3\\
-x_3&-x_2 t&0\\
0&-x_3&x_2\\
-x_{s_2} t&-x_1&0\\
0&-x_{s_2} t^2&x_1\\
x_{v_1}&x_{s_1} t^2&0\\
0&x_{v_1}&-x_{s_1} t\\
\end{bmatrix}
$

Substituting $t=1$ and $t=0$ in the complex $\mathcal{F}_{t,\Pi}$ gives us complexes $\mathcal{F}_{1,\Pi}$ and $\mathcal{F}_{0,\Pi}$, respectively.\qed
\end{example}

Next, we construct the complex $\mathcal{F}_{t,\Pi}$ of $R_{n}[t]$-modules that is the homogenisation of the complex $\mathcal{F}_{1,\Pi}$ with respect to the weight vector $\lambda$. This process starts with homogenisations of free modules $F_{1,\Pi,j}=\bigoplus_\mathbf{c} R_n\cdot e_{\mathbf{c}}$ for all $j\geq 0$ (Equation (\ref{reseq})). 
  With the help of the weight vector $\lambda$, we assign weights to the variables of $R_n$ and the basis elements $e_{\mathbf{c}}$ of $\mathcal{F}_{1,\Pi,j}$. The weight assigned to the basis element $e_{\mathbf{c}}$ is   $\epsilon_{\mathbf{c}}=\lambda\cdot D(\mathcal{C})$ where $\mathcal{C}\in \mathbf{c}$ is the $V_k$-acyclic $(j+1)$-partition, and the weights assigned to the variables $x_i$ are $\lambda_{i}$ for all $1\leq i \leq n$. With respect to these weights, an element $m={\bf{x}}^{\bf u}e_{\mathbf{c}}\in F_{1,\Pi,j}$ has the weight $w({\bf{x}}^{\bf u}e_{\mathbf{c}})=\lambda\cdot {\bf u}+\epsilon_{\bf{c}}$. The homogenisation of  $F_{1,\Pi,j}$ with respect to $\lambda$ is the free module $F_{t,\Pi,j}=\bigoplus_\mathbf{c} R_n[t]\cdot \tilde{e}_{\mathbf{c}}$ where $R_n[t]\cdot \tilde{e}_{\mathbf{c}}$ is the free $R_{n}[t]$-module of rank one whose basis element $\tilde{e}_{\mathbf{c}}$ has degree
 $(D(\mathbf{c}),\epsilon_{\mathbf{c}})
 $ and ${\rm deg}(t)=({\bf{0}},1)$. Note that $F_{1,\Pi,j}$ is $\mathbb{Z}^n/L_{\Pi}$-graded $R_n$-module. The process of homogenisation converts it into $\mathbb{Z}^n/L_{\Pi}\times \mathbb{Z}$-graded $R_{n}[t]$-module.
 For a detailed discussion on homogenisations of free modules, we refer to  \cite[Chapter 8, Section 3]{miller2005combinatorial}.

 In the process of homogenisation, the complex $\mathcal{F}_{1,\Pi}$ converts into a complex $\mathcal{F}_{t,\Pi}$ of $R_n[t]$-modules where $\mathcal{F}_{t,\Pi}$ is the $\mathbb{Z}^n/L_{\Pi}\times \mathbb{Z}$-graded complex whose $j$-th homological degree module $F_{t,\Pi,j} = \bigoplus_{\mathbf{c}} R_{n}[t]\tilde{e}_{\mathbf{c}}$. The $j$-th differential map $\delta_{t,j}:F_{t,\Pi,j}\rightarrow F_{t,\Pi,j-1}$ of $\mathcal{F}_{t,\Pi}$ is defined as
\begin{equation}\label{homdiftop}
    \delta_{t,j}(\tilde{e}_{\bf{c}})=t^{\epsilon _{\bf{c}}}\sum_f{{\rm sign}_{{\mathbf{c}}}(f)t^{-w(m_{\bf{c}}(f)e_{{\bf{c}}/f})} m_{\mathbf{c}}(f)\cdot\tilde{e}_{{\mathbf{c}}/f}}
\end{equation}
where the sum is taken over all the contractible edges $f$ of $\bf{c}$. The complex $\mathcal{F}_{t,\Pi}$ becomes $\mathcal{F}_{1,\Pi}$ when evaluated at $t=1$. Moreover, for any nonzero $t_0\in \mathbb{K}$, the complex $\mathcal{F}_{t_0,\Pi}$ is isomorphic to $\mathcal{F}_{1,\Pi}$. For more details on homogenisation of $\mathcal{F}_{1,\Pi}$ (in the case of connected $n$-partitions of $C_n)$, we refer to \cite[Section 6]{manjunathschwil}.

\begin{lemma}\label{Lemma 6.1}
The complex $\mathcal{F}_{t,\Pi}/({(t)\mathcal{F}_{t,\Pi}})$  is isomorphic to $\mathcal{F}_{0,\Pi}$ as complexes of $R_{n}$-modules.
 \end{lemma}
\begin{proof}
 Let $\mathbf{c}$ be an equivalence class of acyclic $j$-partitions of $(C_n,\Pi).$ Let $\mathcal{C},  \mathcal{B}\in \mathbf{c}$ and $\mathcal{B}\neq \mathcal{C}$, where $\mathcal{C}$ is the $V_k$-acyclic $j$-partition. From \cite[Lemma 6.2]{manjunathschwil}, to prove that $\mathcal{F}_{t,\Pi}/({(t)\mathcal{F}_{t,\Pi}})$ is isomorphic to $\mathcal{F}_{0,\Pi}$, it is enough to show that $\lambda\cdot D(\mathcal{C})>\lambda \cdot D(\mathcal{B})$. Since $\mathcal{C}, \mathcal{B}$ lie in the same equivalence class, we can obtain $\mathcal{B}$ from $\mathcal{C}$ by a sequence of main vertex firings not including $V_k$. Let $\sigma \in \mathbb{N}^n$ be such that the entries of $\sigma$ represent the number of chip-firings from the main vertices of $(C_n,\Pi)$, then $\sigma(j)$ is $0$ for all $j\in V_k$ and is constant for all basic vertices in a given $V_i$, for all $i$ from $1$ to $k$. This implies that $\lambda\cdot (D(\mathcal{C})- D(\mathcal{B}))=\lambda\cdot (\sigma\Lambda_{C_n})=\lambda \cdot u $, where $u=\sum_{i=1}^{k}{\sigma(i) b_{V_i}}$. Since $\lambda\cdot b_{V_i}>0$ for all $i$ from $1$ to $k-1$ and  $\sigma(i)\geq0$, this implies that $\lambda\cdot (D(\mathcal{C})- D(\mathcal{B}))>0$.
\end{proof}

Now, we prove the exactness of the complex $\mathcal{F}_{t,\Pi}$. The proof of Theorem \ref{minrehomtop} closely follows \cite[Theorem 6.4]{manjunathschwil}.
\begin{theorem}\label{minrehomtop}Given any parcycle $(C_n,\Pi)$ and any $G$-parking function ideal $M_{\Pi}$ of $(C_n,\Pi)$, there exists a Gr\"obner degeneration $R_{n}[t]/I_{\Pi,t}$ of $R_{n}/I_{\Pi}$ to $R_{n}/M_{\Pi}$ such that $\mathcal{F}_{t,\Pi}$ is a minimal free resolution of $R_{n}[t]/I_{\Pi,t}$. 
\end{theorem}
 \begin{proof}
 Since $\mathcal{F}_{1,\Pi}$ is a complex and $\mathcal{F}_{t,\Pi}$ is $\mathbb{Z}$-graded (in the weighted sense), the sequence $\mathcal{F}_{t,\Pi}$ is a complex. In the following, we prove the exactness of the complex $\mathcal{F}_{t,\Pi}$. Consider the short exact sequence 
$$0 \longrightarrow \mathbb{K}[t] \overset{t}\longrightarrow \mathbb{K}[t] \longrightarrow \mathbb{K}[t]/(t)\mathbb{K}[t] \longrightarrow 0
$$ and tensor it with $H_j(\mathcal{F}_{t,\Pi})$ over $R_n[t]$. The module 
$H_j(\mathcal{F}_{t,\Pi}) \otimes R_n[t]/(t)R_n[t]$ is contained in $H_j(\mathcal{F}_{t,\Pi} \otimes R_n[t]/(t)R_n[t]) = H_j(\mathcal{F}_{0,\Pi})$ and the equality holds from Lemma \ref{Lemma 6.1}. From Theorem \ref{complex of G-parking ideal}, the homology module
$H_j(\mathcal{F}_{0,\Pi})$ is $0$ and hence, $H_j(\mathcal{F}_{t,\Pi}) \otimes R_n[t]/(t)R_n[t]=0$  for all $j\geq 1$. The right exactness of the tensor functor and the vanishing of the homology module
$H_j(\mathcal{F}_{0,\Pi})$ implies that the sequence
$$H_j(\mathcal{F}_{t,\Pi}) \otimes \mathbb{K}[t]\overset{t}\longrightarrow  H_j(\mathcal{F}_{t,\Pi})\otimes \mathbb{K}[t] \longrightarrow 0
$$ is exact.
Hence, from Nakayama's lemma  $H_j(\mathcal{F}_{t,\Pi}) = 0$ for all $j\geq 1$. This implies that the complex $\mathcal{F}_{t,\Pi}$ is exact. From \cite[Proposition 3.2.2]{herzog2011monomial}, the homogenisation of a Gr\"obner basis of an ideal is a Gr\"obner basis of the homogenised ideal. From Lemma \ref{Lemma 6.1} and Lemma \ref{construction of parking ideal}, the set 
\begin{equation*}
\{ \delta_{1,1}(e_{\mathbf{c}}) :\mathbf{c} \text{ is a  chip-firing equivalence class of acyclic}~ 2\text{-partitions of}~(C_{n},\Pi)
\} \end{equation*} 
is a Gr\"obner basis of $I_{\Pi}$ under $m_{\rm{rev}}$ and this monomial order is graded with respect to the weight order $\lambda$. Hence, the ideal $I_{\Pi,t}$ is the image of $\delta_{t,1}$ and the complex $\mathcal{F}_{t,\Pi}$ is a free resolution of $R_n[t]/I_{\Pi,t}$. Since the coefficients in the image of the differentials $\delta_{t,i}$ are non-units, it is a minimal free resolution of $R_n[t]/I_{\Pi,t}$.
\end{proof}

As a corollary to Theorem \ref{minrehomtop}, we prove that $\mathcal{F}_{1,\Pi}$ is a minimal free resolution of $R_n/I_{\Pi}$.
\begin{corollary}\label{corollary 6.3}
The complex $\mathcal{F}_{1,\Pi}$ is a minimal free resolution of $R_n/I_{\Pi}$.
\end{corollary}
\begin{proof}
From \cite[Corollary 3.2.5]{herzog2011monomial},  $(t-1)$ is a nonzero divisor of $R_n[t]/I_{\Pi,t}$ and the exactness of $\mathcal{F}_{1,\Pi}$ follows from the exactness of $\mathcal{F}_{t,\Pi}$ \cite[Proposition 8.28]{miller2005combinatorial}. Since the cokernel of $\delta_{1,1}$ is $R_n/I_{\Pi}$ and the image of $\delta_{1,j}$ is contained in $\langle x_1,\dots,x_n \rangle F_{1,\Pi,j-1}$ for all $j$ from $1$ to $k-1$, the complex $\mathcal{F}_{1,\Pi}$ is a minimal free resolution of $R_n/I_{\Pi}$. 
\end{proof}
Now, with the help of Corollary \ref{corollary 6.3} and Theorem \ref{complex of G-parking ideal}, we define a relation between the Betti numbers of the toppling ideal $I_{\Pi}$ and the Betti numbers of the $G$-parking function ideal $M_{\Pi}$.  Let $\alpha(\tilde{\Pi})$ be the number of acyclic orientations on $(C_{n})_{\tilde{\Pi}}$ with a unique sink at the set containing $V_k$ and let $\mathcal{P}_{i+1,j}$ be the set of all connected $(i+1)$-partitions $\tilde{\Pi}$ of $(C_n,\Pi)$ such that $(C_{n})_{\tilde{\Pi}}$ has $j$ edges in total.
\begin{corollary}\label{equality of betti numbers}
Let $(C_n,\Pi)$ be a parcycle where $\Pi$ is a connected $k$-partition of the vertex set of $C_n$. The graded Betti numbers $\beta_{i,j}(R_n/I_{\Pi})=\sum_{\tilde{\Pi} \in \mathcal{P}_{i+1,j}}{\alpha(\tilde{\Pi})}\newline =\beta_{i,j}(R_n/M_{\Pi})$ for all $i,j.$
\end{corollary}
\begin{proof}
From Corollary \ref{corollary 6.3} and Theorem \ref{complex of G-parking ideal}, we know that the complexes $\mathcal{F}_{1,\Pi}$ and $\mathcal{F}_{0,\Pi}$ are minimal free resolutions of $R_n/I_{\Pi}$ and $R_n/M_{\Pi}$, respectively.
Hence, the Betti numbers of $I_{\Pi}$ and $M_{\Pi}$ are $\beta_{i}(R_n/I_{\Pi})=\sum_{\tilde{\Pi} \in \mathcal{P}_{i+1}}{\alpha(\tilde{\Pi})}=\beta_{i}(R_n/M_{\Pi})$ where $\mathcal{P}_{i+1}$ is the set of all connected $(i+1)$-partitions of $(C_n,\Pi)$ and $\alpha(\tilde{\Pi})$ denotes the number of acyclic orientations on $(C_{n})_{\tilde{\Pi}}$ with a unique sink at the set containing $V_k$. Hence, the graded Betti numbers of $I_{\Pi}$ and $M_{\Pi}$ are $\beta_{i,j}(R_n/I_{\Pi})=\sum_{\tilde{\Pi} \in \mathcal{P}_{i+1,j}}{\alpha(\tilde{\Pi})}=\beta_{i,j}(R_n/M_{\Pi})$.
\end{proof}

\section{From Parcycles to Rational Normal Curves}\label{PartoRNT}

In this section, we relate the two main objects of this article: the parcycle and the rational normal curve. 
\subsection{Toppling Ideal of a Cycle Graph and the Rational Normal Curve} \label{obtaining the rational normal curve from toppling ideal}

As we mentioned in the introduction, the starting point of our investigation is 
the observation that the rational normal curve of degree $n$ and the toppling ideal of an $n$-cycle have the same graded Betti table. In the following, we elaborate on this connection.
We start by defining a relation between the defining ideal $P_n$ of the rational normal curve $\Gamma_n$ and the toppling ideal $I_{C_n}$ of $C_n$. Proposition \ref{thm7} is a first step towards uncovering the relation between the rational normal curve and the parcycle.
 \begin{proposition}\label{thm7} Let $\phi : R_{n+1} \longrightarrow R_n$ be the graded ring homomorphism defined by
$\phi(x_{n+1}) = x_1$ and 
$
\phi(x_i)= x_i
$ for all $i \neq n+1$. The image of the ideal $P_{n}+\langle x_{n+1}-x_{1}\rangle$ under $\phi$ is $I_{C_n}$, and $(x_{n+1}-x_1)$ is a nonzero divisor on $R_{n+1}/P_n$. 
\end{proposition}
\begin{proof}
The defining ideal $P_n$ of $\Gamma_{n}$ is generated by the set of all $2\times 2$ minors of the matrix $T_n$ (Equation (\ref{matTn})). In other words, $P_n$ is generated by the set $\{P_{i,j}=x_{i}x_{j}-x_{i-1}x_{j+1} \mid 2\leq i\leq j\leq n\}$. From Lemma \ref{one genereic of A}, $T_n$ is a $1$-generic matrix and hence, $P_n$ is a prime ideal. This implies that $(x_{n+1}-x_1)$ is a nonzero divisor on $\mathbb{K}[x_1,\dots,x_{n+1}]/P_n$. 
 From Lemma \ref{generating lemma}, the toppling ideal $I_{C_n}$ is generated by the set $\{\mathbf{x}^{S \rightarrow \overline{S}}-\mathbf{x}^{\overline{S}\rightarrow S}\mid S\subseteq V(C_n)\setminus\{n\}\,\text{ and the subgraph induced by} ~ S ~\text{on } \,C_n \,\text{is connected}
\}$. 
To prove that $\phi(P_n+\langle x_{n+1}-x_{1}\rangle) =I_{C_n}$, it is enough to show that $\phi(P_n)=I_{C_n}.$
First, we prove that $\phi(P_n) \subseteq I_{C_n}$ by showing that the image of every generator of $P_{n}$ under $\phi$ lies in $I_{C_n}$.
Let $P_{i,j}=x_{i}x_{j}-x_{i-1}x_{j+1}\in P_{n}$ be such that $2 \leq i\leq j \leq n-1.$ The image $\phi(P_{i,j}) =x_{i}x_{j}-x_{i-1}x_{j+1}= \mathbf{x}^{S \rightarrow \overline{S}}-\mathbf{x}^{\overline{S}\rightarrow S}$ where $S= \{i,\dots,j \} \subset V(C_n)\setminus\{n\}$ and the subgraph induced by $S$ on $C_n$ is connected. In the case of $j=n$ and $2\leq i\leq n$, the image $\phi(P_{i,n}) = x_ix_n -x_{i-1} x_{1} =  -(\mathbf{x}^{S \rightarrow \overline{S}}-\mathbf{x}^{\overline{S}\rightarrow S})$, where $S =\{1,\dots,i-1\} \subseteq V(C_n)\setminus\{n\}$ and the subgraph induced by $S$ on $C_n$ is connected. This implies that $\phi(P_n)\subseteq I_{C_n}$.

We now prove $I_{C_n} \subseteq \phi(P_{n})$.
Let $S$ be a subset of $V(C_n)-\{n\}$ such that the induced subgraph by $S$ on $C_n$ is connected. This implies that either $S =\{i,\dots,j\} $ where $2\leq i \leq j \leq n-1$, or $S = \{1,\dots,j \}$ where $1\leq j \leq n-1 $. In the first case, $\mathbf{x}^{S \rightarrow \overline{S}}-\mathbf{x}^{\overline{S}\rightarrow S}  =x_i x_j - x_{i-1}x_{j+1}=\phi(x_i x_j - x_{i-1}x_{j+1}) \in \phi(P_{n}),$ and in the second case $\mathbf{x}^{S \rightarrow \overline{S}}-\mathbf{x}^{\overline{S}\rightarrow S} =x_1 x_j - x_{n}x_{j+1}=\phi(x_{j}x_{n+1}-x_{n}x_{j+1}) \in \phi(P_n)$. This implies that $I_{C_n}\subseteq \phi(P_n)$. Hence, the image of the ideal $P_n+\langle x_{n+1}-x_{1}\rangle$ under $\phi$ is $I_{C_n}$. 
\end{proof}
As a consequence of the above result, we have the following corollary.
\begin{corollary}\label{ratntoppcyc_cor}
The map $\phi$ defined in Proposition \ref{thm7} induces an isomorphism $\tilde{\phi}:R_{n+1}/\langle x_{n+1}-x_1\rangle \rightarrow R_n$ and $R_{n+1}/P_{n}\otimes_{R_{n+1}} R_{n+1}/\langle x_{n+1}-x_1\rangle \simeq R_n/I_{C_n}.$  
\end{corollary}

\subsection{Cohen-Macaulay Initial Monomial Ideals of the Rational Normal Curve and Parcycles}\label{construction of parcycles}

Given a Cohen-Macaulay initial monomial ideal $\mathcal{M}$ of the defining ideal $P_n$ of $\Gamma_n$, we construct a parcycle $C_{\Pi, \mathcal{M}}$ as discussed in the outline of the thesis. Recall from the preliminaries (Theorem \ref{Concanegrirossi main theorem}) that every Cohen-Macaulay initial monomial ideal 
of $P_n$ is generated by certain diagonal and anti-diagonal elements of the matrix $T_n$ (Equation (\ref{matTn})).

The following two extreme cases (with the minimum and maximum values for $k$ in Theorem \ref{Concanegrirossi main theorem}) shed light on the general construction. This is a key technical contribution of this article.

\begin{example}\label{cminitpar_ex} \rm

Consider the case $k=1$, i.e. $i_0=1$ and $i_1=n+1$. The corresponding initial monomial ideal of $P_n$ is  $\langle x_2^2,x_2 x_3,\dots, x_2 x_{n}, x_3^2,x_3x_4,\newline \dots x_3 x_{n}, \dots, x_{n-1}x_n, x_{n}^2\rangle$.
This ideal is the $G$-parking function ideal of $C_{\Pi,\mathcal{M}}=(C_{n+1},\Pi)$ where $\Pi=\{\{1,n+1\},\{2\},\{3\},\dots,\{n\}\}$ and $\{1,n+1\}$ is the sink.

The other extreme is when $k=n$, i.e. $i_0=1,i_1=2,\dots,i_n=n+1$. The corresponding initial ideal is \begin{center} $\langle x_1x_3, x_1x_4,\dots,x_1x_{n+1},x_2x_4,x_2x_5,\dots,x_2x_{n+1},\dots,x_{n-1} x_{n+1} \rangle.$ \end{center}
In this case, the corresponding parcycle consists of $2n$ basic vertices with two adjacent basic vertices labelled $s_1,s_2,$ and the other $2n-2$ basic vertices labelled $1,\dots,n-1$ and $v_1,\dots,v_{n-1}$. The corresponding partition consists of $n$ subsets of size two. The subset $\{s_{1},s_{2}\}$ designated to be the sink, and the other subsets of size two are of the form $\{i,v_i\}$ for  $i$ from $1$ to $n-1$. The basic edges of the parcycle are $\{i,v_i\},\{s_1,s_2\}$ for  $i$ from $1$ to $n-1$ and the relevant edges of the parcycle are $\{v_i,i+1\}$ for $i$ from $1$ to $n-2$, along with $\{1,s_2\}$ and $\{v_{n-1},s_1\}$.
We refer to Figure \ref{extremeneqk} for an illustration.
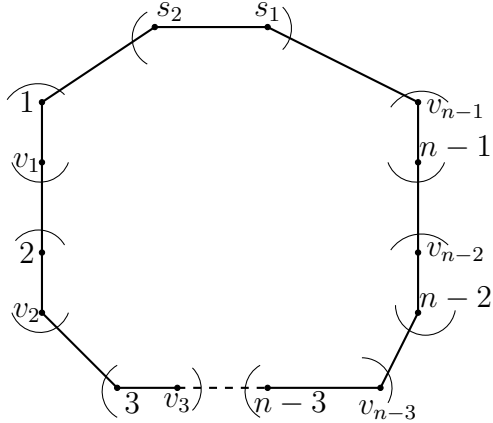
\begin{figure}[h]
\centering
\begin{tikzpicture}
\draw[ -] (0.3,0.1) arc (45:145:0.5);
\draw[fill=black](0,0) circle (1pt);
\node at (-0.2,0){$1$};
\draw[fill=black](0,-0.8) circle (1pt);
\node at (-0.2,-0.8){$v_{1}$};
\draw[thick] (0,0)--(0,-0.8);
\draw[-] (-0.4,-0.8) arc (200:340:0.4);
\draw[fill=black](0,-2) circle (1pt);
\node at (-0.2,-2){$2$};
\draw[thick] (0,-0.8)--(0,-2);
\draw[ -] (0.3,-1.9) arc (30:140:0.4);
\draw[fill=black](0,-2.8) circle (1pt);
\node at (-0.2,-2.8){$v_{2}$};
\draw[thick] (0,-2)--(0,-2.8);
\draw[-] (-0.4,-2.8) arc (200:340:0.4);
\draw[fill=black](1,-3.8) circle (1pt);
\node at (1.2,-4){$3$};
\draw[thick] (0,-2.8)--(1,-3.8);
\draw[-] (1,-3.5) arc (120:240:0.4);
\draw[fill=black](1.8,-3.8) circle (1pt);
\draw[thick] (1,-3.8)--(1.8,-3.8);
\draw[-] (2,-4.1) arc (-45:45:0.4);
\node at (1.8,-4){$v_{3}$};
\draw[fill=black](3,-3.8) circle (1pt);
\node at (3.3,-4){$n-3$};
\draw[thick] (1.8,-3.8)--(3,-3.8)[dashed];
\draw[ -] (2.9,-3.5) arc (120:240:0.4);
\draw[fill=black](4.5,-3.8) circle (1pt);
\draw[thick] (3,-3.8)--(4.5,-3.8);
\node at (4.6,-4.1){$v_{n-3}$};
\draw[ -] (4.6,-4) arc (-30:90:0.4);
\draw[fill=black](5,-2.8) circle (1pt);
\node at (5.5,-2.6){$n-2$};
\draw[thick] (4.5,-3.8)--(5,-2.8);
\draw[-] (4.7,-2.7) arc (180:350:0.4);
\draw[fill=black](5,-2) circle (1pt);
\node at (5.5,-2){$v_{n-2}$};
\draw[thick] (5,-2.8)--(5,-2);
\draw[-] (5.4,-1.9) arc (45:145:0.5);
\draw[fill=black](5,-0.8) circle (1pt);
\node at (5.5,-0.6){$n-1$};
\draw[thick] (5,-2)--(5,-0.8);
\draw[-] (4.6,-0.8) arc (200:340:0.4);
\draw[fill=black](5,0) circle (1pt);
\node at (5.5,-0.1){$v_{n-1}$};
\draw[thick] (5,-0.8)--(5,0);
\draw[ -] (5.4,0) arc (45:145:0.5);
\draw[fill=black](3,1) circle (1pt);
\node at (3,1.2){$s_{1}$};
\draw[fill=black](1.5,1) circle (1pt);
\node at (1.7,1.2){$s_{2}$};
\draw[thick] (1.5,1)--(3,1);
\draw[thick] (3,1)--(5,0);
\draw[thick] (1.5,1)--(0,0);
\draw[-] (1.4,1.2) arc (120:240:0.4);
\draw[ -] (3.2,0.7) arc (-45:45:0.4);
\end{tikzpicture}
\caption{Parcycle in the case $k=n$.}
\label{extremeneqk}
\end{figure}

The $G$-parking function ideal of the above parcycle is the square-free monomial ideal $M_{\Pi,\mathcal{M}}=\langle x_{i}x_{v_{j}}\mid 1\leq i\leq j\leq n-1\rangle$. From Theorem \ref{regularity of A in R/M lemma}, the sequence $\mathcal{A}=(x_{v_j}-x_{i_{j+2}},x_{s_1}-x_{i_{k-1}},x_{s_2}-x_{i_1}\mid 0\leq j\leq k-2)$ is a regular sequence on $\mathcal{S}/M_{\Pi,\mathcal{M}}$ where  $\mathcal{S}=\mathbb{K}[x_1,\dots,x_{n+1},x_{v_0},\dots,x_{v_{k-2}},x_{s_1},x_{s_2}]$, and $\mathcal{S}/M_{\Pi,\mathcal{M}}\otimes_{\mathcal{S}} \mathcal{S}/\langle \mathcal{A} \rangle \simeq \mathbb{K}[x_1,\dots,x_{n+1}]/\mathcal{M}$. 
\qed
\end{example}

These two extreme cases lead to the following ``interpolation'' for the general case.

\textbf{Construction of Parcycles in the General Case:} For a Cohen-Macaulay initial monomial ideal $\mathcal{M}$ of $P_n$ defined by the sequence $1=i_0<i_1<\dots <i_k=n+1$, we associate a parcycle $C_{\Pi,\mathcal{M}}$ with $n+k$ vertices.

\begin{figure}[h]
\centering
\begin{tikzpicture}
\draw[fill=black](0,0) circle (2pt);
\node at (-0.9,0){$i_{1}-1$};
\draw[fill=black](-0.5,-0.8) circle (2pt);
\node at (-1.2,-0.8){${i_{1}-2}$};
\draw[thick] (0,0)--(-0.50,-0.8);

\draw[fill=black](-0.7,-2) circle (2pt);
\node at (-1.2,-2){$2$};
\draw[thick] (-0.5,-0.8)--(-0.7,-2)[dashed];

\draw[fill=black](-0.7,-2.8) circle (2pt);
\node at (-1.2,-2.8){1};
\draw[thick] (-0.7,-2)--(-0.7,-2.8);
\draw[thick, -](-0.4,-2.6) arc (45:145:0.5);

\draw[fill=black](-0.6,-3.5) circle (2pt);
\node at (-1,-3.5){$v_0$};
\draw[thick] (-0.7,-2.8)--(-0.6,-3.5);
\draw[thick, -](-0.9,-3.7) arc (240:350:0.5);

\draw[fill=black](-0.3,-4) circle (2pt);
\node at (-0.5,-4.4){$i_{2}-1$};
\draw[thick] (-0.6,-3.5)--(-0.3,-4);

\draw[fill=black](0.5,-4.8) circle (2pt);
\node at (0,-5.1){$i_{2}-2$};
\draw[thick] (-0.3,-4)--(0.5,-4.8);

\draw[fill=black](1.5,-5.2) circle (2pt);
\node at (1.5,-5.6){$i_{1}+1$};
\draw[thick] (0.5,-4.8)--(1.5,-5.2)[dashed];

\draw[fill=black](2.2,-5.2) circle (2pt);
\node at (2.4,-5.5){$i_{1}$};
\draw[thick] (1.5,-5.2)--(2.2,-5.2);
\draw[thick, -](2,-4.8) arc (145:220:0.5);

\draw[fill=black](2.8,-5.2) circle (2pt);
\node at (2.9,-5.5){$v_{1}$};
\draw[thick] (2.2,-5.2)--(2.8,-5.2);
\draw[thick, -](2.9,-5.4) arc (-30:40:0.5);

\draw[fill=black](3.8,-5.2) circle (2pt);
\node at (4.1,-5.6){$i_{k-4}$};
\draw[thick] (2.8,-5.2)--(3.8,-5.2)[dashed];
\draw[thick, -](3.7,-4.8) arc (145:220:0.5);

\draw[fill=black](4.3,-5) circle (2pt);
\node at (4.8,-5.3){$v_{k-4}$};
\draw[thick] (3.8,-5.2)--(4.3,-5);
\draw[thick, -](4.6,-5.2) arc (-30:90:0.5);

\draw[fill=black](5,-4.5) circle (2pt);
\node at (5.8,-4.5){$i_{k-2}-1$};
\draw[thick] (4.3,-5)--(5,-4.5);

\draw[fill=black](5.5,-3.8) circle (2pt);
\node at (6.4,-3.8){$i_{k-3}+1$};
\draw[thick] (5,-4.5)--(5.5,-3.8)[dashed];

\draw[fill=black](5.7,-3) circle (2pt);
\node at (6.2,-2.8){$i_{k-3}$};
\draw[thick] (5.5,-3.8)--(5.7,-3); 
\draw[thick, -](5.2,-2.9) arc (180:340:0.5);

\draw[fill=black](5.8,-2.2) circle (2pt);
\node at (6.3,-2.2){$v_{k-3}$};
\draw[thick] (5.7,-3)--(5.8,-2.2);
\draw[thick, -](6.2,-2) arc (45:145:0.5);

\draw[fill=black](5.8,-1.6) circle (2pt);
\node at (6.6,-1.6){$i_{k-1}-1$};
\draw[thick] (5.8,-2.2)--(5.8,-1.6); 

\draw[fill=black](5.8,-0.8) circle (2pt);
\node at (6.8,-0.8){$i_{k-2}+1$};
\draw[thick] (5.8,-1.6)--(5.8,-0.8)[dashed];

\draw[fill=black](5.6,-0) circle (2pt);
\node at (6.3,0.2){$i_{k-2}$};
\draw[thick] (5.8,-0.8)--(5.6,0);
\draw[thick, -](5.2,0) arc (180:340:0.5);

\draw[fill=black](5,0.8) circle (2pt);
\node at (5.8,0.8){$v_{k-2}$};
\draw[thick] (5.6,-0)--(5,0.8);
\draw[thick, -](5,1.1) arc (75:180:0.5);

\draw[fill=black](4.2,1.2) circle (2pt);
\node at (4.2,1.5){$n$};
\draw[thick] (5,0.8)--(4.2,1.2);

\draw[fill=black](3,1.2) circle (2pt);
\node at (3,1.5){$i_{k-1}+1$};
\draw[thick] (4,1.2)--(3,1.2)[dashed];

\draw[fill=black](1.8,1.2) circle (2pt);
\node at (1.8,1.5){$s_{1}$};
\draw[thick] (3,1.2)--(1.8,1.2);
\draw[thick, -](2,0.8) arc (-45:45:0.5);

\draw[fill=black](0.8,0.8) circle (2pt);
\node at (0.8,1.2){$s_{2}$};
\draw[thick] (1.8,1.2)--(0.8,0.8);
\draw[thick] (0.8,0.8)--(0,0);
\draw[thick, -](0.4,1) arc (180:270:0.5);
\end{tikzpicture}
\caption{The parcycle $C_{\Pi,\mathcal{M}}$.}
\label{agenparcy}
\end{figure}
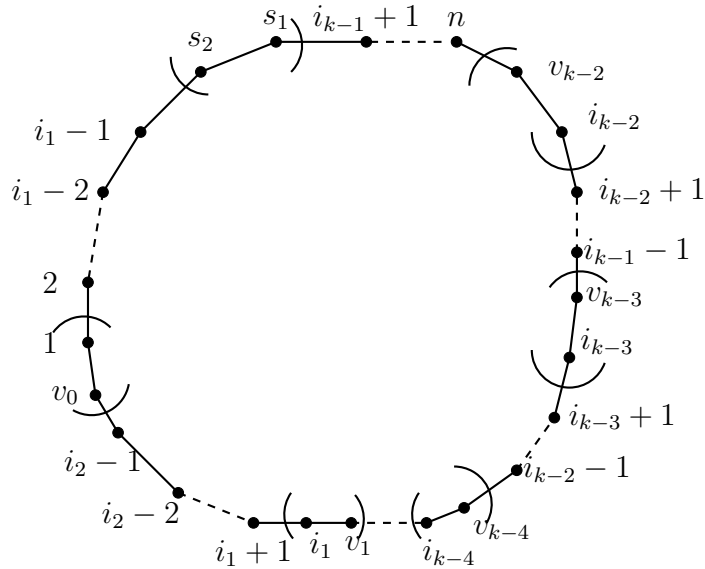

 The vertices of $C_{\Pi,\mathcal{M}}$ are labelled as follows:  two distinguished basic vertices, namely $s_1,s_2$, one collection of vertices called vertices of  {\emph{type I}} labelled $i_r+1,\dots,i_{r+1}-1$ for all $r$ from $0$ to $k-1$ where $i_{r+1}-i_r \geq 2$  and another collection called vertices of \emph{type II} labelled $i_0,v_0,i_1,v_1,\dots,i_{k-2},v_{k-2}$. The partition $\Pi$ consists of the vertices of type I as singletons (hence, there are $n-k$ singletons in total), the vertices of type II  appear in subsets of size two as $\{i_r,v_r\}$ for $r$ from $0$ to $k-2$
and $\{s_1,s_2\}$ is designated to be the sink (hence, there are $k$ subsets of size two in total). Hence, the partition consists of $n$ subsets. We refer to Figure \ref{agenparcy} for an illustration.

We now turn to the relevant edges of $C_{\Pi,\mathcal{M}}$.  Assume that $k \geq 2$ (the case $k=1$ has already been treated in  Example \ref{cminitpar_ex}),  the relevant edges of $C_{\Pi,\mathcal{M}}$ are of the following types: 

\begin{enumerate}
\item  $\{i_{1}-1,s_2\}$. The edge $\{i_{k-1}+1,s_1\}$ if $ i_{k}-i_{k-1}\geq 2 $, and  $\{v_{k-2},s_1\}$ if $ i_{k}-i_{k-1}=1.$ 

\item $\{p,p-1\}$  if both $p$ and $p-1$ are vertices of type I. The edges $\{i_1-1,i_1-2\},\{i_2-1,i_2-2\}$ in Figure \ref{agenparcy} are examples of such edges.

\item  $\{i_{r}+1,i_{r}\}$ if $i_{r+1}-i_{r}\geq 2$, for all $r$ from $0$ to $k-2$. The case $i_{r+1}-i_{r}=1$ is considered in item \ref{fourth}.  The edges $\{2,1\},\{i_{k-3}+1,i_{k-3}\},\{i_{k-2}+1,i_{k-2}\}$ in Figure \ref{agenparcy} are examples of such edges.

\item\label{fourth}  $\{v_r,i_{r+2}-1\}$ where $r$ varies from $0$ to $k-3$, and $\{v_{k-2},i_{k}-1\}$ if $i_{k}-i_{k-1}\geq 2$. The edges $\{v_0,i_{2}-1\}, \{v_{k-4},i_{k-2}-1\}$ in Figure \ref{agenparcy} are examples of such edges.

\end{enumerate}
\begin{example}\label{cminmidp4}\rm
Consider the Cohen-Macaulay initial monomial ideal $\mathcal{M}=\langle  x_1x_3,x_1x_4,x_1x_5, x_2x_5,x_3^{2},x_3x_5 \rangle$ of $P_4$. The ideal $\mathcal{M}$ corresponds to the sequence $i_0=1<i_1=2<i_2=4<i_3=5$. The parcycle $C_{\Pi,\mathcal{M}}$ associated with $\mathcal{M}$ is depicted in Figure \ref{paroeg}.
    \qed
    \end{example}

\subsubsection{Edge Ideals and Graph Quotients}
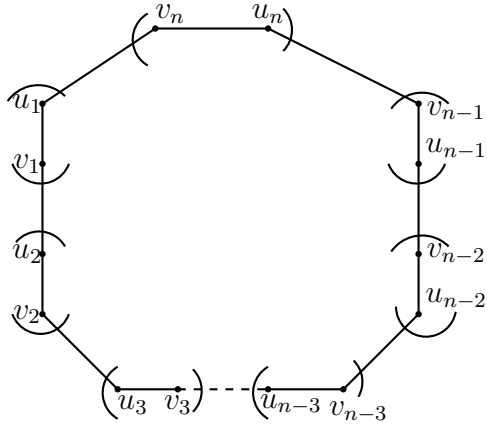
\begin{figure}[h]
    \centering
    \begin{tikzpicture}
\draw[thick, -] (0.3,0.1) arc (45:145:0.5);
\draw[fill=black](0,0) circle (1pt);
\node at (-0.2,0){$u_{1}$};
\draw[fill=black](0,-0.8) circle (1pt);
\node at (-0.2,-0.8){$v_{1}$};
\draw[thick] (0,0)--(0,-0.8);
\draw[thick, -] (-0.4,-0.8) arc (200:340:0.4);
\draw[fill=black](0,-2) circle (1pt);
\node at (-0.2,-2){$u_{2}$};
\draw[thick] (0,-0.8)--(0,-2);
\draw[thick, -] (0.3,-1.9) arc (30:140:0.4);
\draw[fill=black](0,-2.8) circle (1pt);
\node at (-0.2,-2.8){$v_{2}$};
\draw[thick] (0,-2)--(0,-2.8);
\draw[thick, -] (-0.4,-2.8) arc (200:340:0.4);
\draw[fill=black](1,-3.8) circle (1pt);
\node at (1.2,-4){$u_{3}$};
\draw[thick] (0,-2.8)--(1,-3.8);
\draw[thick, -] (1,-3.5) arc (120:240:0.4);
\draw[fill=black](1.8,-3.8) circle (1pt);
\draw[thick] (1,-3.8)--(1.8,-3.8);
\draw[thick, -] (2,-4.1) arc (-45:45:0.4);
\node at (1.8,-4){$v_{3}$};
\draw[fill=black](3,-3.8) circle (1pt);
\node at (3.3,-4){$u_{n-3}$};
\draw[thick] (1.8,-3.8)--(3,-3.8)[dashed];
\draw[thick, -] (3,-3.5) arc (120:240:0.4);
\draw[fill=black](4,-3.8) circle (1pt);
\draw[thick] (3,-3.8)--(4,-3.8);
\node at (4.2,-4.1){$v_{n-3}$};
\draw[thick, -] (4.2,-3.9) arc (-30:45:0.4);
\draw[fill=black](5,-2.8) circle (1pt);
\node at (5.5,-2.6){$u_{n-2}$};
\draw[thick] (4,-3.8)--(5,-2.8);
\draw[thick, -] (4.7,-2.7) arc (180:350:0.4);
\draw[fill=black](5,-2) circle (1pt);
\node at (5.5,-2){$v_{n-2}$};
\draw[thick] (5,-2.8)--(5,-2);
\draw[thick, -] (5.4,-1.9) arc (45:145:0.5);
\draw[fill=black](5,-0.8) circle (1pt);
\node at (5.5,-0.6){$u_{n-1}$};
\draw[thick] (5,-2)--(5,-0.8);
\draw[thick, -] (4.6,-0.8) arc (200:340:0.4);
\draw[fill=black](5,0) circle (1pt);
\node at (5.5,-0.1){$v_{n-1}$};
\draw[thick] (5,-0.8)--(5,0);
\draw[thick, -] (5.4,0) arc (45:145:0.5);
\draw[fill=black](3,1) circle (1pt);
\node at (3,1.2){$u_{n}$};
\draw[fill=black](1.5,1) circle (1pt);
\node at (1.7,1.2){$v_{n}$};
\draw[thick] (1.5,1)--(3,1);
\draw[thick] (3,1)--(5,0);
\draw[thick] (1.5,1)--(0,0);
\draw[thick, -] (1.4,1.2) arc (120:240:0.4);
\draw[thick, -] (3.2,0.7) arc (-45:45:0.4);
\end{tikzpicture}
    \caption{The Uniform Parcycle.}
    \label{graphU}
\end{figure}

   In the following, we shed more light on the construction of the parcycle $C_{\Pi,\mathcal{M}}$.  To this end, we consider the \emph{uniform} parcycle $U_{2n}$ on $2n$ basic vertices, see Figure \ref{graphU}. Note that  $\mathcal{M}, M_{U_{2n}}$ (the $G$-parking function ideal of $U_{2n}$) and $M_{\Pi, \mathcal{M}}$ are monomial ideals that are quadratically generated.
   We associate a graph (possibly with loops) with each of them, denoted as  $\Sigma_{\mathcal{M}}, G_{U_{2n}}$ and $G_{\mathcal{M}}$ respectively, and realise the ideals as edge ideals of the associated graphs.
 The construction of $C_{\Pi,\mathcal{M}}$ is based on an isomorphism between a certain quotient graph of $G_{U_{2n}}$ and $\Sigma_{\mathcal{M}}$. 

\begin{itemize}
    \item {\bf The Graph $G_{U_{2n}}$:}
The $G$-parking function ideal $M_{U_{2n}}$ of $U_{2n}$ is the ideal $\langle x_{u_i}x_{v_j} \mid 1\leq i\leq j\leq n-1\rangle$. The ideal $ M_{U_{2n}}$ is the edge ideal of the bipartite graph whose vertex set is the union of the disjoint sets $\{u_1, \dots, u_{n-1}\}$ and $\{v_1, \dots, v_{n-1}\}$, and whose edge set is $\{\{u_i, v_j\} \mid 1 \leq i \leq j \leq n-1\}$.
Note that the generating set $\{x_{u_i}x_{v_j} \mid 1\leq i\leq j\leq n-1\}$ of $M_{U_{2n}}$ is the set of all anti-diagonals of matrix
\[T_{U}=
\begin{bmatrix}
x_{v_n} & x_{v_1} & x_{v_2} & \dots & x_{v_{n-1}} \\
x_{u_1} & x_{u_2} & x_{u_3} & \dots & x_{u_{n}}
\end{bmatrix}.
\]

\item {\bf The Graph ${\Sigma}_{\mathcal{M}}$:} 
Let $i_0=1<i_1<\dots<i_{k}=n+1$ be the sequence associated with $\mathcal{M}$.
From Theorem \ref{Concanegrirossi main theorem}, the ideal $\mathcal{M}$ is generated by certain diagonals and anti-diagonals of the matrix 
\[
T_n = 
\resizebox{\linewidth}{!}{$
\left[
\begin{array}{ccccccccccccccc}
x_1 & x_2 & \cdots & x_{i_1-1} & \vdots & x_{i_1} & \cdots & x_{i_2-1} & \vdots & \cdots & x_{i_{k-1}-1} & \vdots & x_{i_{k-1}} & \cdots & x_n \\
x_2 & x_3 & \cdots & x_{i_1} & \vdots & x_{i_1+1} & \cdots & x_{i_2} & \vdots & \cdots & x_{i_{k-1}} & \vdots & x_{i_{k-1}+1} & \cdots & x_{n+1} \\
\end{array}
\right]
$}
\]
We rearrange the columns and rows of \(T_n\) so that the ideal \(\mathcal{M}\) is generated by the set of anti-diagonals of the new matrix $\tilde{T}_{n}$ where 
\[
\tilde{T}_{n} = 
\resizebox{\linewidth}{!}{$
\left[
\begin{array}{ccccccccccccccc}
x_{i_1}&x_{i_1-1}&\cdots&x_2&\vdots&x_{i_2}&\cdots&x_{i_{1}+1}&\vdots&\cdots&\vdots&x_{n+1}&\cdots&x_{i_{k-1}+1}\\
x_{i_{1}-1}&x_{i_1-2}&\cdots&x_1&\vdots&x_{i_2-1}&\cdots&x_{i_{1}}&\vdots&\cdots&\vdots&x_{n}&\cdots&x_{i_{k-1}}\\
\end{array}
\right].
$}
\]
For more details on the column and row operations that we perform on $T_n$ to obtain the matrix $\tilde{T}_{n}$ is discussed in Lemma \ref{lemma parkingtensorint}.
In the following, we define the graph ${\Sigma}_{\mathcal{M}}$.
The vertex set of ${\Sigma}_{\mathcal{M}}$ is $\{    
(\tilde{T}_{n})_{(1,i)},(\tilde{T}_{n})_{(2,j)}\mid 2\leq i\leq n, 1\leq j \leq n-1
\}$ where $(\tilde{T}_{n})_{(i,j)}$ is the index of the entry corresponding to the $i$-th row and $j$-th column of $\tilde{T}_{n}$. The edge set of ${\Sigma}_{\mathcal{M}}$ is $\{\{(\tilde{T}_{n})_{(1,i)},(\tilde{T}_{n})_{(2,j)}\}\mid 1\leq j< i\leq n\}.$

\item {\bf The Graph ${G}_{\mathcal{M}}$:} To construct the graph ${G}_{\mathcal{M}}$, we first define a map $\phi_{\mathcal{M}}$ from $\{u_1,\dots,u_n,v_1,\dots,v_n\}$ to $\{1,\dots,n+1\}$. The map $\phi_{\mathcal{M}}$ maps $u_i$ and $v_j$ to elements of $\{1, \dots, n+1\}$ in such a way that
$$\phi_{\mathcal{M}}(T_{U})=\begin{bmatrix}
x_{\phi_{\mathcal{M}}(v_n)} & x_{\phi_{\mathcal{M}}(v_1)} & x_{\phi_{\mathcal{M}}(v_2)} & \dots & x_{\phi_{\mathcal{M}}(v_{n-1})} \\
x_{\phi_{\mathcal{M}}(u_1)} & x_{\phi_{\mathcal{M}}(u_2)} & x_{\phi_{\mathcal{M}}(u_3)} & \dots & x_{\phi_{\mathcal{M}}(u_{n})}
\end{bmatrix}=\tilde{T}_{n}.$$

The vertex set of $G_{\mathcal{M}}$ is defined as $V(G_{\mathcal{M}})=V(G_{U_{2n}})/\sim$, where $u\sim v$ if $\phi_{\mathcal{M}}(u)=\phi_{\mathcal{M}}(v)$. The total number of edges between two vertices $[u]$ and $[v]$ of $G_{\mathcal{M}}$ is the total number of edges in $G_{U_{2n}}$ between the vertices in $[u]$ and the vertices in $[v]$. The total number of loops at the vertex $[u]$ is the total number of edges in $G_{U_{2n}}$  between the vertices in $[u]$. \end{itemize}

Note that the map $\phi_{\mathcal{M}}$ induces an isomorphism of graphs from $G_{\mathcal{M}}=G_{U_{2n}}/\sim$ onto $\Sigma_{\mathcal{M}}$ where $\phi_{\mathcal{M}}([u])=\phi_{\mathcal{M}}(u)$. 

With the help of the map $\phi_{\mathcal{M}}$, the parcycle $C_{\Pi,\mathcal{M}}$ can be realised as a certain quotient graph of the uniform parcycle $U_{2n}$.
We identify a main vertex $\{u_i,v_i\}$ of $U_{2n}$, where $1\leq i \leq n-1$, with a type $I$ vertex, say $l$, if $\phi_{\mathcal{M}}(u_i)=l=\phi_{\mathcal{M}}(v_i)$. The main vertices $\{u_j,v_j\}$ of $U_{2n}$ for which $\phi_{\mathcal{M}}(u_j)\neq \phi_{\mathcal{M}}(v_j)$, we keep them as they are.
While constructing the quotient graph of $U_{2n}$, we assume that the edges of the underlying cycle that are not relevant edges do not produce any edge in the quotient graph. The total number of edges between two vertices (Type $I$, Type $II$) in this quotient graph is the total number of relevant edges between their corresponding main vertices in $U_{2n}$. In this way, we obtain the parcycle $C_{\Pi,\mathcal{M}}$ as a quotient graph of the uniform parcycle $U_{2n}$.

\section{Regular Sequences}\label{CompandConn}We start by extending Proposition \ref{thm7} to relate the rational normal curve and the toppling ideal of a parcycle. From now onwards, we work with the following setting.  Let $\mathcal{M}$ be a Cohen-Macaulay initial monomial ideal of $P_n$ with the associated sequence $i_0=1<i_1<\dots<i_{k}=n+1$ (Theorem \ref{Concanegrirossi main theorem}). Let $C_{\Pi,\mathcal{M}}$ be the parcycle associated with $\mathcal{M}$ as defined in Section \ref{construction of parcycles} and $I_{\Pi,\mathcal{M}}$ be its toppling ideal (Section \ref{Gb and construction of G-parking}). Let $A_{(C,\Pi)}$ be the matrix associated with $C_{\Pi,\mathcal{M}}$ with the fixed main vertex $\{s_1,s_2\}$ (defined in Subsection \ref{relation between one-generic matrix and toppling ideal}). The matrix $A_{(C,\Pi)}$ associated with $C_{\Pi,\mathcal{M}}$ is as follows:
\begin{equation}\label{eq5}
 \resizebox{\linewidth}{!}{%
$
\left[  \begin{array}{cccccccccccccc}
x_{s_2}&x_{i_{1}-1}&\dots&x_2&x_{v_{0}}&x_{i_{2}-1}&\dots&x_{i_{1}+1}&x_{v_1}&\dots&x_{i_{k-2}+1}&x_{v_{k-2}}&\dots&x_{i_{k-1}+1}\\
     x_{i_{1}-1}&x_{i_{1}-2}&\dots&x_{1}&x_{i_{2}-1}&x_{i_{2}-2}&\dots&x_{i_{1}}&x_{i_{3}-1}&\dots&x_{i_{k-2}}&x_{n}&\dots&x_{s_1}
     \end{array}\right].
$}
\end{equation}
      Consider the polynomial ring $\mathcal{S}=\mathbb{K}[x_1,\dots,x_{n+1},x_{v_0},\dots,x_{v_{k-2}},x_{s_1},x_{s_2}]$ in $n+k+2$-variables over $\mathbb{K}$ and the sequence $\mathcal{A}=(x_{v_j}-x_{i_{j+2}},x_{s_1}-x_{i_{k-1}},x_{s_2}-x_{i_1}\mid 0\leq j\leq k-2)$ in $\mathcal{S}$. Note that $\mathcal{S}$ is a $R_{n+1}$-module. 

      \begin{theorem}\label{parcycle to defining ideal}
Let $\mathcal{M}$ be a Cohen-Macaulay initial monomial ideal of $P_{n}$ with the associated sequence $i_0=1<i_1<\dots<i_{k}=n+1$. The $\mathbb{Z}$-graded modules $\mathcal{S}/I_{\Pi,\mathcal{M}} \otimes_{\mathcal{S}} \mathcal{S}/\langle\mathcal{A}\rangle$ and  $R_{n+1}/P_{n}$ are isomorphic as $R_{n+1}$-modules.
\end{theorem}
\begin{proof}
Consider the $\mathbb{K}$-algebra homomorphism $\psi$ that is defined as
$$\psi:\mathcal{S}\rightarrow R_{n+1}
$$
$$\psi(x_i)=x_i,\, \forall \,1\leq i\leq n+1,
$$
$$\psi(x_{v_{j}})=x_{i_{j+2}},\, \forall \,0\leq j\leq k-2,
$$
$$
\psi(x_{s_1})=x_{i_{k-1}}, \psi(x_{s_2})=x_{i_1}.
$$

The kernel of $\psi$ is $\langle\mathcal{A}\rangle$ \cite[1.8.4]{geck2013introduction}. To prove that $\mathcal{S}/I_{\Pi,\mathcal{M}} \otimes_{\mathcal{S}} \mathcal{S}/\langle\mathcal{A}\rangle \simeq R_{n+1}/P_{n}$, it is enough to show that $\psi(I_{\Pi,\mathcal{M}}+\langle\mathcal{A}\rangle)=P_n$.
From Corollary \ref{prime and cohen macaulayness of toppling ideal}, the toppling ideal $I_{\Pi,\mathcal{M}}$ is generated by the set of all $2\times 2$ minors of $A_{(C,\Pi)}$. 
 The image ideal $\psi(I_{\Pi,\mathcal{M}}+\langle \mathcal{A}\rangle)$ is generated by the set of all $2\times 2$ minors of the matrix $\psi(A_{(C,\Pi)})$, where $\psi(A_{(C,\Pi)})$ is the matrix obtained from $A_{(C,\Pi)}$ by replacing the variables $x_{v_{j}}$ with $x_{i_{j+2}}$, for all $0\leq j \leq k-2$ and $x_{s_1},x_{s_2}$ with $x_{i_{k-1}}$ and $x_{i_1}$, respectively. From the construction of $A_{(C,\Pi)}$, the matrix $\psi(A_{(C,\Pi)})$
  is equal to the matrix $\tilde{T}$ that is obtained from  
$T_n\label{Tn} = \begin{bmatrix}
x_1 &x_2 &\dots &x_{n}\\
x_2 & x_3 &\dots &x_{n+1} 
\end{bmatrix}$ by first rearranging the columns of $T_n$ according to the following sequence $(i_{1}-1,i_{1}-2,\dots,1,i_{2}-1,\dots,i_{1},i_3-1,\dots,i_2,\dots,i_{k-2}, i_{k}-1=n,i_{k}-2\dots,i_{k-1})$ (here, $i_j-l$ is an index of a column of $T_n$) and then interchanging the rows.
 Hence, the image of the ideal $I_{\Pi,\mathcal{M}}+\langle \mathcal{A}\rangle $ under $\psi$ is the defining ideal $P_n$ of  $\Gamma_n$.
\end{proof}

\subsection{Gr\"obner Weights}\label{weights}
Let $\mathcal{M}$ be a Cohen-Macaulay initial monomial ideal of $P_n$ and let $C_{\Pi,\mathcal{M}}$ be the parcycle associated with it. Let $\tilde{\mathcal{S}}$ be the polynomial ring over $\mathbb{K}$ whose variables correspond to the basic vertices of $C_{\Pi,\mathcal{M}}$. Note that the polynomial ring $\mathcal{S}=\tilde{\mathcal{S}}[x_{i_{k-1}},x_{i_k}]$.
In this section, we construct an integral weight vector $\lambda\in \mathbb{N}^{n+k}$ such that $\lambda\cdot b_{V_i}>0$ for all $V_i \in \Pi$ that do not contain the sink. We further extend the weight vector $\lambda$ to a weight vector $\tilde{\lambda}\in\mathbb{N}^{n+k+2}$ on $\mathcal{S}$, 
so that the sequence $\mathcal{A}=(x_{v_j}-x_{i_{j+2}},x_{s_1}-x_{i_{k-1}},x_{s_2}-x_{i_1}\mid 0\leq j\leq k-2)$ becomes a homogeneous sequence on the polynomial ring $\mathcal{S}$ with the grading induced by $\tilde{\lambda}$.
We will use these properties of $\lambda$ and $\tilde{\lambda}$ to prove the regularity of $\mathcal{A}$ on $S/I_{\Pi,\mathcal{M}}$ where $I_{\Pi,\mathcal{M}}$ is the toppling ideal of $C_{\Pi,\mathcal{M}}$.

The construction follows a greedy strategy where we traverse the basic vertices of $C_{\Pi,\mathcal{M}}$ in the cyclic order $s_2\dots 1 v_0\dots i_{k-2}v_{k-2}\dots s_1$ starting from the sink vertex $s_2$. We assign ``partial'' values for  $\lambda$ so as to satisfy those constraints that involve the vertices encountered till then. By ``partial'', we mean only for those coordinates of $\lambda$ that correspond to the vertices we have encountered. The ``local'' nature of the constraints ensures the success of such a greedy strategy. The details are as follows. 

We start the computation of $\lambda$ by first fixing an integral value for $\lambda$ at the vertex $s_2$.  Assume that there is at least one type I vertex between $s_2$ and $i_0=1$. Let $(a_0,a_1,\dots,a_m)=(s_2,i_{1}-1,\dots,i_0)$ be the arrangement (according to the cyclic order) of basic vertices of $C_{\Pi,\mathcal{M}}$ from $s_2$ to $i_0$. Note that $a_j$ is a type I vertex between $s_2$ and $i_0$, for all $j$ from $1$ to $m-1$. For every type I vertex $V_i=\{a_i\}$, the inequality $\lambda\cdot b_{V_i}>0$ specialises to $2\lambda_{a_i} > \lambda_{a_{i-1}}+\lambda_{a_{i+1}}$. 
We define 
\begin{center}
$\lambda_{a_i}=\lambda_{a_0}+m(m+1)/2-(m-i)(m-i+1)/2,\, \forall\, 1\leq i\leq m-1,$ \end{center}
 and choose a suitable integer value for $\lambda_{a_m}$ that satisfies the inequality $2\lambda_{a_{m-1}}>\lambda_{a_{m-2}}+\lambda_{a_m}$. Note that
$$
     2\lambda_{a_{i}}-\lambda_{a_{i-1}}-\lambda_{a_{i+1}}=-2 \sum_{j=1}^{m-i} j + \sum_{j=1}^{m-i+1} j + \sum_{j=1}^{m-i-1} j >0,\, \forall\, 1\leq i\leq m-2.$$
 Hence, $\lambda_{a_i}$ satisfies the inequality $\lambda \cdot b_{V_{i}}>0$ for all $1\leq i \leq m-1$. 

Suppose there are no type I vertices between $s_2$ and $i_0$. In this case, the inequality $\lambda\cdot b_{V_{1}}>0$ for $V_1=\{i_0,v_0\}\in \Pi$ becomes either $\lambda_{i_0}+\lambda_{v_0}>\lambda_{s_2}+\lambda_{i_{1}}$, or $\lambda_{i_0}+\lambda_{v_0}>\lambda_{s_2}+\lambda_{i_{2}-1}$. The first case arises if there are no type I vertices between $v_0$ and $i_1$ and the second case arises if there is at least one type I vertex between $v_0$ and $i_1$. In order to homogenise $x_{s_2}-x_{i_1}\in \mathcal{A}$  with respect to $\tilde{\lambda}$, we set  $\lambda_{i_1}:=\lambda_{s_2}$. In the first case, we select suitable values for $\lambda_{i_0}$ and $\lambda_{v_0}$ so as to satisfy the inequality $\lambda_{i_0}+\lambda_{v_0}>\lambda_{s_2}+\lambda_{i_{1}}$, and in the second case, we choose a value for $\lambda_{i_0}$ and obtain the inequality $\lambda_{v_0}>\lambda_{i_{2}-1}+l$
 where $l=\lambda_{s_2}-\lambda_{i_0}\in \mathbb{Z}$.
 
 We now know the values of $\lambda_{s_2},\dots,\lambda_{i_0}$. 
 Let $p-1$ be the number of type I vertices between $v_0$ and $i_1$.
 Let $(u_0,u_1,\dots,u_p)=(v_0,i_{2}-1,\dots,i_1)$ be the arrangement of the vertices of $C_{\Pi,\mathcal{M}}$ from $v_0$ to $i_1$. Note that $u_j$ is a type I vertex between $v_0$ and $i_1$, for all $j$ from $1$ to $p-1$. 
 In the case when there is at least one type I vertex between $s_2$ and $i_0$,
 we know the values of $\lambda_{i_0}=\lambda_{a_m}, \lambda_{i_0+1}=\lambda_{a_{m-1}}$ and recall that $\lambda_{i_1}=\lambda_{s_2}$.
  The inequalities $\lambda_{i_0}+\lambda_{v_0}>\lambda_{i_{0}+1}+\lambda_{i_{2}-1}$ and $\lambda_{i_0}+\lambda_{v_0}>\lambda_{s_2}+\lambda_{i_{2}-1}$ transform into an inequality of the form $\lambda_{u_0}>\lambda_{u_1}+w$ where $w\in \mathbb{Z}$. We define
 \begin{equation}\label{solution of ui}
\lambda_{u_i}=(p-i)\cdot w+p(p+1)/2-i(i+1)/2+\lambda_{i_1}, \, 0\leq i\leq p. 
 \end{equation}

A simple calculation shows that $\lambda_{u_0} > \lambda_{u_1}+w$ and 
$2\lambda_{u_{i}}-\lambda_{u_{i-1}}-\lambda_{u_{i+1}}=-2 \sum_{j=1}^{i} j + \sum_{j=1}^{i-1} j + \sum_{j=1}^{i+1} j >0,$
for all $1\leq i \leq p-1$. Hence,  $\lambda_{u_i}$ satisfies the inequality $\lambda\cdot b_{V_i}>0$ for all $1\leq i \leq p-1$ where $V_{i}=\{u_i\}$.

For any $j\geq 1$, the process of computing values of $\lambda_{s_2},\dots,\lambda_{i_0},\lambda_{v_0},\dots,\lambda_{i_{j-1}},\newline \lambda_{v_{j-1}},\dots,\lambda_{i_j}$ is the same as above.  We compute the values of $\lambda_{v_{j-1}},\dots,\lambda_{i_j}$ with the help of the known values of $\lambda_{s_2},\dots,\lambda_{i_0},\lambda_{v_0},\dots,\lambda_{i_{j-1}}$ and $\lambda_{i_j}:=\lambda_{v_{j-2}}, \lambda_{s_1}:=\lambda_{v_{k-3}}$ as in Equation \ref{solution of ui}. 
We choose $\kappa \in \mathbb{N}$ sufficiently large such that $\lambda':= \lambda+ \kappa \cdot (1,\dots,1)\in \mathbb{N}^{n+k}$. Since $(1,1,\dots,1)\cdot b_{V_i}=0$ for all $i$, $\lambda' \cdot b_{V_i}>0$ for all $V_i$ not containing $\{s_1,s_2\}$. This  yields a Gr\"obner weight as required. For the description of the Gr\"obner cone of $\mathcal{M}$, the reader is urged to refer to  \cite[Definition 4.12]{conca2007contracted}.

\subsection{Properties of $\mathcal{A}$}
Consider the polynomial ring $R_p=\mathbb{K}[x_1,\dots,x_p]$ where $p\geq n+1$. For a monomial $u=x_1^{d_1}\cdots x_{p}^{d_p}\in R_p$, we define ${\rm{Supp}}(u):=\{i\mid d_{i}\neq 0\}$. 
Let $A$ be a matrix of order $2\times n$ with entries of the form $kx_i$ where $k\in \mathbb{K}\setminus\{0\}$ and $x_i\in R_p$ is an indeterminate over $\mathbb{K}$.  In every row of $A$, let all the entries be distinct. 
We assume that there are two entries $k_1x_a$ and $k_2x_b$ in $A$ that are the only entries in $A$ with support $a$ and $b$, respectively and they lie in different rows of $A$. Let $\overline{{\rm{ind}}}(kx_i)$ and ${\rm{ind}}(k'x_j)$ be the indices of the columns in $A$ that contain $kx_i$ and $k'x_j$ in their first and second rows, respectively. Without loss of generality, we assume that $k_1x_a$ lies in the first row of $A$ and $k_2x_b$ lies in the second row of $A$. Let $\overline{{\rm{ind}}}(k_1x_a)<{\rm{ind}}(k_2x_b)$ and ${\rm{Supp}}(kx_{c})\neq {\rm{Supp}}(k'x_{t})$ for all entries of $A$ that satisfy ${\rm{ind}}(kx_c)<\overline{{\rm{ind}}}(k_1x_a)<{\rm{ind}}(k_2x_b)<\overline{{\rm{ind}}}(k'x_t)$.
Let $M$ be the monomial ideal generated by the set of all anti-diagonals $\tilde{\mathcal{D}}$ of the $2\times 2$ minors of $A$. We prove the regularity of the element $k_{1}x_{a}-k_{2}x_{b}$ on $R_p/M$.

\begin{lemma}\label{reglemma}
    The polynomial $k_1x_a-k_2x_b$ is a regular element on $R_p/M$.
\end{lemma}
\begin{proof}
   Let $m\in R_p$ be a polynomial such that $(k_1x_a-k_2x_b)m\in M$. Since $M$ is a monomial ideal, we can take $m$ to be a monomial (Proposition \ref{moidpro}). Furthermore, $(k_1x_a-k_2x_b)m\in M$ implies that $x_a m$ and $x_b m\in M$. From \cite[Proposition 1.1.5]{herzog2011monomial}, there exist monomials, say $u$ and $v$, in the generating set $\tilde{\mathcal{D}}$ of $M$ such that $u$ and $ v$ divide $x_{a}m$ and $x_{b}m$, respectively. If $a \notin {\rm{Supp}}(u)$, then $u$ divides $m$ and this implies $m \in M$. Suppose that $a \in {\rm{Supp}}(u)$, this implies that $m$ is divisible by some $x_c$ where $x_c$ (with a unit multiple, say $k$) is a second row entry of $A$ that satisfies ${\rm{ind}}(kx_c)<\overline{{\rm{ind}}}(k_1x_a)$. This is because the only elements in $\tilde{\mathcal{D}}$ that contain $x_a$ in their support are of the form $x_a x_c$ such that ${\rm{ind}}(kx_c)<\overline{{\rm{ind}}}(k_1x_a)$.
   
Since $v\in \tilde{\mathcal{D}}$ divides $x_{b}m$, we have the following cases:
 $(1)~ b \notin {\rm{Supp}}(v)$ and in this case $v$ divides $m$. $(2)~ b \in {\rm{Supp}}(v)$ and in this case $m$ is divisible by some $x_t$ where $x_t$ (with a unit multiple, say $k'$) is a first row entry of $A$ that satisfies ${\rm{ind}}(k_2x_b)<\overline{{\rm{ind}}}(k'x_t)$. 
 
 In the case when both $a$ and $b$ do not lie in the support of $u$ and $v$, respectively, the variables $x_c$ and $x_t$ divide $m$. Since ${\rm{ind}}(kx_c)<\overline{{\rm{ind}}}(k_1x_a)<{\rm{ind}}(k_2x_b)<\overline{{\rm{ind}}}(k'x_t)$, the monomial $x_c x_t \in \tilde{\mathcal{D}}$ and divides $m$. Hence, $k_1x_a-k_2x_b$ is a regular element on $R_p/M$.
\end{proof}

 In the following, we prove the regularity of the sequence $\mathcal{A}$ on the quotient rings of $M_{\Pi,\mathcal{M}}$ and $I_{\Pi,\mathcal{M}}$. Let $\mathcal{D}$ be the set of all anti-diagonals of the $2\times 2$ minors of the matrix $A_{(C,\Pi)}$. From Corollary \ref{GparantiD}, the $G$-parking function ideal $M_{\Pi,\mathcal{M}}$ of $C_{\Pi,\mathcal{M}}$ is generated by the set  $\mathcal{D}$.

\begin{theorem}\label{regularity of A in R/M lemma}
The sequence $\mathcal{A}=(x_{v_j}-x_{i_{j+2}},x_{s_1}-x_{i_{k-1}},x_{s_2}-x_{i_1}\mid 0\leq j\leq k-2)$ is a regular sequence on $\mathcal{S}/M_{\Pi,\mathcal{M}}$ and $\mathcal{S}/I_{\Pi,\mathcal{M}}$. \end{theorem}
\begin{proof}We begin by proving the regularity of the sequence $(x_{v_j}-x_{i_{j+2}}\mid 0\leq j\leq k-2)$ on $\mathcal{S}/M_{\Pi,\mathcal{M}}$.
First, we consider the case where $k\geq 4$ and the result for the case $k\leq 3$ will be proved later. We prove the regularity of the sequence $(x_{v_j}-x_{i_{j+2}}\mid 0\leq j\leq k-4)$ on $\mathcal{S}/M_{\Pi,\mathcal{M}}$ with the help of Lemma \ref{reglemma} and
by applying induction on $j$. The base case $j=0$ corresponds to showing that $x_{v_0}-x_{i_2}$ is a regular element on $\mathcal{S}/M_{\Pi,\mathcal{M}}$.
The variables $x_{v_0}$ and $x_{i_2}$ appear only once in $A_{(C,\Pi)}$ (as defined in (\ref{eq5})) and lie in the first and second rows, respectively. Note that in every row of $A_{(C,\Pi)}$, all the entries are distinct. Since
$\overline{\rm{ind}}(x_{v_0})<{\rm{ind}}(x_{i_2})$ in $A_{(C,\Pi)}$ and $x_{c}\neq x_{t}$ whenever ${\rm{ind}}(x_c)<\overline{{\rm{ind}}}(x_{v_0})<{\rm{ind}}(x_{i_2})<\overline{{\rm{ind}}}(x_t)$, the polynomial $x_{v_0}-x_{i_2}$ is a regular element on $\mathcal{S}/M_{\Pi,\mathcal{M}}$ (Lemma \ref{reglemma}).

The induction hypothesis is that $\mathcal{A}'=(x_{v_j}-x_{i_{j+2}} \mid 0\leq j \leq l-1)$ is a regular sequence on $\mathcal{S}/M_{\Pi,\mathcal{M}}$.
We now prove that $x_{v_l}-x_{i_{l+2}}$ is a regular element on $\mathcal{S}/(M_{\Pi,\mathcal{M}}+\langle \mathcal{A}'\rangle)$.
Consider the $\mathbb{K}$-algebra homomorphism 
$$\tilde{\psi}:\mathcal{S}\rightarrow \mathbb{K}[x_1,\dots,x_{n+1},x_{v_{l}},\dots,x_{v_{k-2}},x_{s_1},x_{s_2}]
$$
$$\tilde{\psi}(x_{v_j})=x_{i_{j+2}}, \, 0\leq j \leq l-1,
$$
$$
\tilde{\psi}(x_i)=x_i,\, \text{otherwise}.
$$
The kernel of $\tilde{\psi}$ is $\langle \mathcal{A}' \rangle$. Note that $\mathcal{S}/M_{\Pi,\mathcal{M}}\otimes_{\mathcal{S}}\mathcal{S}/\langle \mathcal{A}'\rangle \simeq \mathcal{S}'/\tilde{\psi}(M_{\Pi,\mathcal{M}})$ where $\mathcal{S}'=\mathbb{K}[x_1,\dots,x_{n+1},x_{v_{l}},\dots,x_{v_{k-2}},x_{s_1},x_{s_2}]$ and $\tilde{\psi}(M_{\Pi,\mathcal{M}})$ is the image of $M_{\Pi,\mathcal{M}}+\langle \mathcal{A}'\rangle$ under $\tilde{\psi}$. The ideal $\tilde{\psi}(M_{\Pi,\mathcal{M}})$ is generated by the set of all anti-diagonals of the $2\times 2$ minors of the matrix $\tilde{\psi}(A_{(C,\Pi)})$, where $\tilde{\psi}(A_{(C,\Pi)})$ is obtained from $A_{(C,\Pi)}$ by replacing the variables $x_{v_j}$ with $x_{i_{j+2}}$ for all $j$ from $0$ to $l-1$. 
By the construction of $A_{(C,\Pi)}$, the variables $x_{v_l}$ and $x_{i_{l+2}}$ appear only once in $\tilde{\psi}(A_{(C,\Pi)})$, in the first and second rows, respectively. From Lemma \ref{reglemma}, $(x_{v_{l}}-x_{i_{l+2}})$ is a regular element on $\mathcal{S}'/\tilde{\psi}(M_{\Pi,\mathcal{M}})$. Hence, the sequence $(x_{v_j}-x_{i_{j+2}}\mid 0\leq j\leq k-4)$ is regular on $\mathcal{S}/M_{\Pi,\mathcal{M}}$.

We now prove that the sequence $(x_{s_1}-x_{i_{k-1}}, x_{s_2}-x_{i_1})$ is regular on $\mathcal{S}/(M_{\Pi,\mathcal{M}}+\langle x_{v_j}-x_{i_{j+2}}\mid 0\leq j\leq k-4 \rangle)$.
Let $\tilde{A}_{(C,\Pi)}$ be the matrix obtained from $A_{(C,\Pi)}$ by replacing the variables $x_{v_j}$ with $x_{i_{j+2}}$ for all $j$ from $0$ to $k-4$. Let $ \tilde{\mathcal{D}}$ be 
the set of all anti-diagonals of the $2\times 2$ minors of $\tilde{A}_{(C,\Pi)}$.
Since the variables $x_{i_{k-1}}$ and $x_{i_{k}}$ are not entries of $A_{(C,\Pi)}$, they are also not entries of $\tilde{A}_{(C,\Pi)}$. This implies that $i_{k-1}$ and $i_{k}$ do not lie in the support of any element of $\tilde{\mathcal{D}}$.
 Since $\langle \tilde{\mathcal{D}} \rangle$ is a monomial ideal in $\mathbb{K}[x_1,\dots,x_{n+1},x_{v_{k-3}},x_{v_{k-2}},x_{s_1},x_{s_2}]$, the sequence $(x_{v_{k-3}}-x_{i_{k-1}},x_{v_{k-2}}-x_{i_{k}})$ is a regular sequence on $\mathbb{K}[x_1,\dots,x_{n+1},x_{v_{k-3}},x_{v_{k-2}},x_{s_1},x_{s_2}]/\langle \tilde{\mathcal{D}} \rangle \simeq \mathcal{S}/(M_{\Pi,\mathcal{M}}+\langle x_{v_j}-x_{i_{j+2}}\mid 0\leq j\leq k-4 \rangle)$. 
 
 In the case where $k\leq 3$, the matrix $\tilde{A}_{(C,\Pi)}=A_{(C,\Pi)}$ and $M_{\Pi,\mathcal{M}}=\langle \tilde{\mathcal{D}} \rangle$.
This proves the regularity of the sequence $(x_{v_j}-x_{i_{j+2}}\mid 0\leq j\leq k-2)$ on $\mathcal{S}/M_{\Pi,\mathcal{M}}$.  
 Since $\{s_1,s_2\}$ is the sink of $C_{\Pi,\mathcal{M}}$, the basic vertices $s_1,s_2$ do not lie in the support of any element of 
$\mathcal{D}$. Following the same steps as in the proof of the regularity of the sequence $(x_{v_{k-3}}-x_{i_{k-1}},x_{v_{k-2}}-x_{i_{k}})$ on $\mathcal{S}/(M_{\Pi,\mathcal{M}}+\langle x_{v_j}-x_{i_{j+2}}\mid 0\leq j\leq k-4 \rangle)$, the sequence $(x_{s_1}-x_{i_{k-1}}, x_{s_2}-x_{i_1})$ is a regular sequence on $\mathcal{S}/(M_{\Pi,\mathcal{M}}+\langle x_{v_j}-x_{i_{j+2}}\mid 0\leq j\leq k-2 \rangle)$. This proves that $\mathcal{A}$ is a regular sequence on $\mathcal{S}/M_{\Pi,\mathcal{M}}$.


The sequence $\mathcal{A}$ is a homogeneous sequence with respect to the weight vector $\lambda$ (described in Section \ref{weights}) associated with $C_{\Pi,\mathcal{M}}$. From Theorem \ref{construction of parking ideal}, the initial ideal of $I_{\Pi,\mathcal{M}}$ with respect to $\lambda$ is $M_{\Pi,\mathcal{M}}.$ As $\mathcal{A}$ is a regular sequence on $\mathcal{S}/M_{\Pi,\mathcal{M}}$, the regularity of $\mathcal{A}$ on $\mathcal{S}/I_{\Pi,\mathcal{M}}$ holds from \cite[Proposition 15.15]{eisenbud2013commutative}.
\end{proof}

As a consequence of the above theorem, we prove that the sequence $\mathcal{A}$ is a regular sequence on $\mathcal{S}[t]/I_{\Pi,\mathcal{M},t}$.
\begin{corollary}\label{Corollary 5.3}
The sequence $\mathcal{A}$ is a regular sequence on $\mathcal{S}[t]/I_{\Pi,\mathcal{M},t}$.
    \end{corollary}
\begin{proof}
By Theorem \ref{construction of parking ideal}, the set $\mathcal{G}=\{\mathbf{x}^{S\rightarrow \overline{S}}-\mathbf{x}^{\overline{S}\rightarrow S}\mid~S~ \text{is a connected-parset of } \newline C_{\Pi,\mathcal{M}} \text{ not containing} \,\{s_1,s_2\} \}$ is a Gr\"obner basis of $I_{\Pi,\mathcal{M}}$ with respect to $m_{\rm{rev}}$. The monomial order $m_{\rm{rev}}$ naturally extends to $\mathcal{S}[t]$ as follows: ${\bf{x}}^{\bf a}t^{c}<_{m_{\rm{rev}}}'{\bf{x}}^{\bf b}t^{d}$ if and only if (i) ${\bf{x^{a}}}<_{m_{\rm{rev}}}{\bf{x^{b}}} $, or (ii) ${\bf{x}^{a}}={\bf{x}^{b}}$ and $c<d.$ From \cite[Proposition 3.1.2]{herzog2011monomial}, the set $\mathcal{G}^{h}=\{g^{h}\mid g\in \mathcal{G}\}$ is a Gr\"obner basis of $I_{\Pi,\mathcal{M},t}$ with respect to $<_{m_{\rm{rev}}}'$.
Consider the weight vector $\lambda^{'}$ in $\mathcal{S}[t]$ such that $\lambda^{'}(y)=\lambda({y})$ for all variables $y\in \mathcal{S}$, and $\lambda^{'}(t)=0.$ The monomial order $<_{m_{\rm{rev}}}'$ is graded with respect to $\lambda^{'}$ and ${\rm{in}}_{<_{m_{\rm{rev}}}'}(I_{\Pi,\mathcal{M},t})\subseteq {\rm{in}}_{\lambda^{'}}(I_{\Pi,\mathcal{M},t}).$ From \cite[Theorem 3.1.2]{herzog2011monomial}, ${\rm{in}}_{{\lambda}'}(I_{\Pi,\mathcal{M},t})={\rm{in}}_{<_{m_{\rm{rev}}}'}(I_{\Pi,\mathcal{M},t})=M_{\Pi,\mathcal{M}}
.$ As $\mathcal{A}$ is a regular sequence on $\mathcal{S}/M_{\Pi,\mathcal{M}}$, the sequence $\mathcal{A}$ is regular on $\mathcal{S}[t]/M_{\Pi,\mathcal{M}}$. Moreover, $\mathcal{A}$ is a homogeneous sequence with respect to $\lambda^{'}$. The regularity of the sequence $\mathcal{A}$ in $\mathcal{S}[t]/I_{\Pi,\mathcal{M},t}$ follows from \cite[Proposition 15.15]{eisenbud2013commutative}. \end{proof}

Next, we show that the quotient ring of $\mathcal{M}$ is the quotient ring of the $G$-parking function ideal of the parcycle $C_{\Pi,\mathcal{M}}$ modulo $\mathcal{A}$.

\begin{lemma} \label{lemma parkingtensorint}
Let $\mathcal{M}$ be a Cohen-Macaulay initial monomial ideal of $P_n$ with the associated sequence $i_0=1<i_1<\dots<i_{k}=n+1$. The $\mathbb{Z}$-graded modules $\mathcal{S}/M_{\Pi,\mathcal{M}} \otimes_{\mathcal{S}} \mathcal{S}/\langle\mathcal{A}\rangle$ and  $R_{n+1}/\mathcal{M}$ are isomorphic as $R_{n+1}$-modules.
\end{lemma}
\begin{proof}
Consider the $\mathbb{K}$-algebra homomorphism $\psi$ from $\mathcal{S}$ to $R_{n+1}$ that is defined in Theorem \ref{parcycle to defining ideal}. The kernel of the $\mathbb{K}$-algebra homomorphism $\psi$ is $\langle \mathcal{A}\rangle$.
To prove that $\mathcal{S}/M_{\Pi,\mathcal{M}}\otimes_{\mathcal{S}}\mathcal{S}/\langle \mathcal{A}\rangle \simeq R_{n+1}/\mathcal{M}$, it is enough to show that $\psi(M_{\Pi,\mathcal{M}})=\mathcal{M}$. From Corollary \ref{GparantiD}, the $G$-parking function ideal $M_{\Pi,\mathcal{M}}$ of $C_{\Pi,\mathcal{M}}$ (with the sink fixed at $\{s_1,s_2\}$) is generated by the set of all anti-diagonals $\mathcal{D}$ of the $2\times 2$ minors of the matrix $A_{(C,\Pi)}$. 
 This implies $\psi(M_{\Pi,\mathcal{M}})$ is generated by the set of all anti-diagonals of the $2\times 2$ minors of the matrix $\psi(A_{(C,\Pi)})$, where $\psi(A_{(C,\Pi)})$ is the matrix obtained from $A_{(C,\Pi)}$ by replacing the variables $x_{v_j}$ with $x_{i_{j+2}}$ for all $j$ from $0$ to $k-2$ and $x_{s_1}, x_{s_2}$ with $x_{i_{k-1}}$ and $x_{i_1}$, respectively. 
Now, we use the following observations to prove that the ideal $\psi(M_{\Pi,\mathcal{M}})=\mathcal{M}$:
\begin{enumerate}
   
    \item{\label{1}} The set of all anti-diagonals $z_{s}y_{t}$ of the $2\times 2$ minors of the matrix $$Y=\begin{bmatrix}
    y_1&y_2&\dots&y_n\\
    z_1&z_2&\dots&z_n
    \end{bmatrix}$$ is equal to the set of all diagonals $z_{u}y_{v}$  of the $2\times 2$ minors of the matrix $$Z=\begin{bmatrix}
    z_1&z_2&\dots&z_n\\
    y_1&y_2&\dots&y_n
    \end{bmatrix}$$ where $s,t$ are column indices of $Y$ and $u,v$ are column indices of $Z$ that satisfy $s<t$ and $u<v$. 
    \item{\label{2}} The set of all diagonals $y_{s}z_{t}$ of the $2\times 2$ minors of the matrix $$Y=\begin{bmatrix}
    y_1&y_2&\dots&y_n\\
    z_1&z_2&\dots&z_n
    \end{bmatrix}$$ is equal to the set of anti-diagonals $z_{u}y_{v}$  of the $2\times 2$ minors of the matrix $$Z=\begin{bmatrix}
    y_n&y_{n-1}&\dots&y_1\\
    z_n&z_{n-1}&\dots&z_1
    \end{bmatrix}$$ where $s,t$ are column indices of $Y$, and $u,v$ are column indices of $Z$ that satisfy $s<t$ and $u<v$.
\end{enumerate}
Let $T_{(C,\Pi)}$ be the matrix obtained by interchanging the rows of $\psi(A_{(C,\Pi)})$.
From the first observation, the ideal $\psi(M_{\Pi,\mathcal{M}})$ is equal to the ideal generated by the set of diagonals of the $2\times 2$ minors of the matrix $T_{(C,\Pi)}$.
    

    From the second observation, the set of diagonals of the $2\times 2$ minors of the matrix 
    \begin{equation}\label{eq7}
        \begin{bmatrix}
    x_{i_{j+1}-1}&x_{i_{j+1}-2}&\dots&x_{i_j}\\
    x_{i_{j+1}}&x_{i_{j+1}-1}&\dots&x_{i_{j}+1}
    \end{bmatrix}
    \end{equation}
    is equal to the set of anti-diagonals of the $2\times 2$ minors of the matrix  $$\begin{bmatrix}
    x_{i_{j}}&x_{i_{j}+1}&\dots&x_{i_{j+1}-1}\\
    x_{i_{j}+1}&x_{i_{j}+2}&\dots&x_{i_{j+1}}
    \end{bmatrix}.$$

    Now, applying Observation \ref{2} to the submatrices of the matrix $T_{(C,\Pi)}$
    that are of the form described in Equation (\ref{eq7}), the ideal $\psi(M_{\Pi,\mathcal{M}})$ is generated by the union of the following sets:

(i) the set of main diagonals $x_{v}x_{r+1}$ of the $2\times 2$ minors of the matrix $T_n$ (defined in Equation (\ref{matTn})) with column indices $v,r$ such that $v+1\leq i_j <r+1$ for some $j$,

(ii) the set of anti-diagonals $x_{v+1}x_{r}$ of the $2\times 2$ minors of the matrix $T_n$ with column indices $v, r$ such that $i_j<v+1<r+1\leq i_{j+1}$ for some $j$.

This proves that $\psi(M_{\Pi,\mathcal{M}})=\mathcal{M}.$
\end{proof}

\section{Theorem \ref{Grodegminres_theo} and Its Applications}\label{mainthmpf}
In this section, we study some special Gr\"obner degenerations of rational normal curves. The special fiber of these degenerations is Cohen-Macaulay initial monomial ideals. We construct explicit minimal free resolutions of these degenerations.


Let $\mathcal{M}$ be a Cohen-Macaulay initial monomial ideal of $P_n$ and $i_0=1<i_1<\dots<i_{k}=n+1$ be the sequence associated with it. Let $C_{\Pi,\mathcal{M}}$ be the parcycle (defined in Section \ref{construction of parcycles}) and $\lambda_{\mathcal{M}} \in \mathbb{N}^{n+k}$ be a weight vector associated with the monomial ideal $\mathcal{M}$ as defined in Section \ref{weights}. Let $I_{\Pi,\mathcal{M}}$ be the toppling ideal of $C_{\Pi,\mathcal{M}}$ and $\mathcal{F}_{1}(C_{\Pi,\mathcal{M}})$ be the complex associated with it (as defined in Section \ref{resolution of toppling ideal}). We construct the complex $\mathcal{F}_{t}(C_{\Pi,\mathcal{M}})$ from $\mathcal{F}_{1}(C_{\Pi,\mathcal{M}})$, as discussed in Section \ref{resolution of toppling ideal}, with the help of the weight vector $\lambda_{\mathcal{M}}$ associated with $\mathcal{M}$. Let $\delta_{t,j}$ be the $j$-th differential map of $\mathcal{F}_{t}(C_{\Pi,\mathcal{M}})$. Let $\tilde{e}_{\bf{c}}$ be the basis element of the free module at the homological degree $j$ of $\mathcal{F}_{t}(C_{\Pi,\mathcal{M}})$ that corresponds to the chip-firing equivalence class ${\mathbf{c}}$ of acyclic $(j+1)$-partitions of $C_{\Pi,\mathcal{M}}$. The image of $\tilde{e}_{\bf{c}}$ under $\delta_{t,j}$ (as in Equation (\ref{homdiftop})) is 
 \begin{equation}\label{hommapratnom}
   \delta_{t,j}(\tilde{e}_{\bf{c}})=t^{\epsilon _{\bf{c}}}\sum_f{{\rm sign}_{{\mathbf{c}}}(f)t^{-w(m_{\bf{c}}(f)e_{{\bf{c}}/f})} m_{\mathbf{c}}(f)\cdot\tilde{e}_{{\mathbf{c}}/f}}.  
 \end{equation}

 The complex $\mathcal{F}_{t}(C_{\Pi,\mathcal{M}})$ is a complex of free modules over the polynomial ring $\tilde{\mathcal{S}}$
whose variables correspond to the basic vertices of $C_{\Pi,\mathcal{M}}$. We define a new complex $\mathcal{F}_{t} :=\mathcal{F}_{t}(C_{\Pi,\mathcal{M}})\bigotimes\mathcal{S}[t]$ where $\mathcal{S}=\mathbb{K}[x_1,\dots,x_{n+1},x_{v_0},\dots,x_{v_{k-2}},\newline x_{s_1},x_{s_2}]$. The complex $\tilde{\mathcal{F}}_t(C_{\Pi,\mathcal{M}}):=\mathcal{F}_{t}\bigotimes \mathcal{S}[t]/\langle \mathcal{A}\rangle$ where $\mathcal{A}=(x_{v_j}-x_{i_{j+2}},x_{s_1}-x_{i_{k-1}},x_{s_2}-x_{i_1}\mid 0\leq j\leq k-2)$ is a homogeneous sequence on $\mathcal{S}$ with respect to the grading induced by $\tilde{\lambda}_{\mathcal{M}}\in \mathbb{N}^{n+k+2}$. The weight vector $\tilde{\lambda}_{\mathcal{M}}$ is an extension of $\lambda_{\mathcal{M}}\in \mathbb{N}^{n+k}$ with some extra properties (defined in Section \ref{weights}).

The complex $\tilde{\mathcal{F}}_t(C_{\Pi,\mathcal{M}})$ is a complex of $R_{n+1}[t]$-modules since $\mathcal{S}[t]/\langle \mathcal{A}\rangle \simeq R_{n+1}[t]$. In particular, 
\begin{equation}\label{consFtil}
  \tilde{\mathcal{F}}_{t}(C_{\Pi,\mathcal{M}}) : \tilde{F}_{t,\Pi,n-1} \overset{\tilde{\delta}_{t,n-1}}\longrightarrow \dots \overset{\tilde{\delta}_{t,2}} \longrightarrow \tilde{F}_{t,\Pi,1} \overset{\tilde{\delta}_{t,1}} \longrightarrow \tilde{F}_{t,\Pi,0},   
\end{equation} 
 
  where $\tilde{F}_{t,\Pi,j} = \bigoplus_{\mathbf{c}} R_{n+1}[t]\bar{e}_{\mathbf{c}}$ and the direct sum is taken over all
 chip-firing equivalence classes $\mathbf{c}$ of acyclic $(j+1)$-partitions of $C_{\Pi,\mathcal{M}}$. The basis element $\bar{e}_{\mathbf{c}}$ has weight $\epsilon _{\bf{c}}$ where $\epsilon _{\bf{c}}=\lambda_{\mathcal{M}}\cdot D(\mathcal{C})$ and $\mathcal{C}\in \mathbf{c}$ is the $V_k$-acyclic $(j+1)$-partition of $C_{\Pi,\mathcal{M}}$.

  The $j$-th differential map $\tilde{\delta}_{t,j}:\tilde{F}_{t,\Pi,j}\rightarrow \tilde{F}_{t,\Pi,j-1}$ of $\tilde{\mathcal{F}}_{t}(C_{\Pi,\mathcal{M}})$ is defined as
\begin{equation}\label{difmaprnc}
    \tilde{\delta}_{t,j}(\bar{e}_{\bf{c}})=t^{\epsilon _{\bf{c}}}\sum_f{{\rm sign}_{{\mathbf{c}}}(f)t^{-w(m_{\bf{c}}(f)e_{{\bf{c}}/f})} \psi(m_{\mathbf{c}}(f))\cdot\bar{e}_{{\mathbf{c}}/f}},
\end{equation}
where the sum is taken over all the contractible edges $f$ of $\bf{c}$. The map $\tilde{\delta}_{t,j}$ is obtained from Equation (\ref{hommapratnom}) with the help of the map $\psi:\mathcal{S}\rightarrow R_{n+1}$ defined in Theorem \ref{parcycle to defining ideal}. 
Note that the complex $\tilde{\mathcal{F}}_t(C_{\Pi,\mathcal{M}})$ is 
$(\lambda,1)$-graded where $\lambda\in \mathbb{N}^{n+1}$ is the restriction of $\tilde{\lambda}_{\mathcal{M}}$ on $R_{n+1}$, i.e. $\lambda(x_j)=\tilde{\lambda}_{\mathcal{M}}(x_j)$ for all $j$ from $1$ to $n+1$. 

\begin{example}\rm
Consider the Cohen-Macaulay initial monomial ideal $\mathcal{M}=\langle  x_1x_3,x_1x_4,x_1x_5, x_2x_5,x_3^{2},x_3x_5 \rangle$ of $P_4$. The ideal $\mathcal{M}$ corresponds to the sequence $i_0=1<i_1=2<i_2=4<i_3=5$. The complex $\tilde{\mathcal{F}}_{t}(C_{\Pi,\mathcal{M}})$ of $R_{5}[t]=\mathbb{K}[x_1,\dots,x_{5}][t]$-modules is as follows:
$$\tilde{\mathcal{F}}_{t}(C_{\Pi,\mathcal{M}}):R_5[t]^3\overset{\tilde{\delta}_{t,3}}\longrightarrow R_5[t]^8 \overset{\tilde{\delta}_{t,2}}\longrightarrow R_5[t]^6 \overset{\tilde{\delta}_{t,1}} \longrightarrow R_5[t]^1
$$
The matrices of the differentials are
\[
\tilde{\delta}_{t,1} = 
\resizebox{\linewidth}{!}{$
\begin{bmatrix}
    x_1 x_{4}-x_3 x_{2} t& x_3^{2}-x_2 x_{4} t &x_2 x_{5}-x_3 x_{4} t & x_1 x_3-x_2^{2} t^2&x_3 x_{5}-x_{4}^{2} t^2&x_1 x_{5}-x_{4} x_{2} t^3
\end{bmatrix}
$}
\]
\\
 $\tilde{\delta}_{t,2}=
\begin{bmatrix}
    0&0&-x_{5}&-x_{4} t^2&0&0&-x_3&-x_2 t\\
    0&0&0&0&-x_{5}&-x_{4} t&-x_{2} t&-x_1\\
    -x_{2} t^2&-x_1&0&0&-x_{4} t&-x_3&0&0\\
    -x_{5}&-x_{4} t&0&0&0&0&x_{4}&x_3\\
    0&0&-x_{2} t&-x_1&x_3&x_2&0&0\\
    x_3&x_2&x_{4}&x_3&0&0&0&0\\
\end{bmatrix}$
\\
\\

$
\tilde{\delta}_{t,3}=
\begin{bmatrix}
x_{4}&x_3&0\\
0&x_{4} t&-x_3\\
-x_3&-x_2 t&0\\
0&-x_3&x_2\\
-x_{2} t&-x_1&0\\
0&-x_{2} t^2&x_1\\
x_{5}&x_{4} t^2&0\\
0&x_{5}&-x_{4} t\\
\end{bmatrix}
$

Substituting $t=1$ and $t=0$ in the complex $\tilde{\mathcal{F}}_{t}(C_{\Pi,\mathcal{M}})$ gives us minimal free resolutions of the quotient rings of $P_4$ and $\mathcal{M}$, respectively.\qed
\end{example}
\subsection{Proof of Theorem \ref{Grodegminres_theo}}
\begin{proof}

Let $\mathcal{M}$ be a Cohen-Macaulay initial monomial ideal of $P_n$. Let $i_0=1<i_1<\dots<i_{k}=n+1$ be the sequence associated with $\mathcal{M}$.
In the case $k=1$, i.e. $i_0=1<i_1=n+1$, the parcycle associated with $\mathcal{M}$ is $C_{\Pi,\mathcal{M}}=(C_{n+1},\Pi)$ where $\Pi=\{\{1,n+1\},\{2\},\{3\},\dots,\{n\}\}$ and $\{1,n+1\}$ is the sink. In this case, the proof of the theorem follows from Lemma \ref{Lemma 6.1} and Theorem \ref{minrehomtop} as $M_{\Pi,\mathcal{M}}=\mathcal{M}$.

Assume that $k\geq 2$. Let $C_{\Pi,\mathcal{M}}$ be the parcycle associated with $\mathcal{M}$ (defined in Subsection \ref{construction of parcycles}).
Let $I_{\Pi,\mathcal{M}}$  and $M_{\Pi,\mathcal{M}}$ be the toppling ideal and the $G$-parking function ideal (with sink $\{s_1,s_2\}$) of the parcycle $C_{\Pi,\mathcal{M}}$, respectively. In Section \ref{weights}, we construct a weight vector $\lambda_{\mathcal{M}}\in \mathbb{N}^{n+k}$ that satisfies $\lambda_{\mathcal{M}}\cdot b_{V_i}>0$ for all main vertices $V_i$ of $C_{\Pi,\mathcal{M}}$ that are different from the sink vertex $\{s_1,s_2\}$ (recall $b_{V_i}$ from introduction). Furthermore, we construct an extension $\tilde{\lambda}_{\mathcal{M}}\in \mathbb{N}^{n+k+2}$ of $\lambda_{\mathcal{M}}$ for which the sequence $\mathcal{A}=(x_{v_j}-x_{i_{j+2}},x_{s_1}-x_{i_{k-1}},x_{s_2}-x_{i_1}\mid 0\leq j\leq k-2)$ becomes a homogeneous sequence on the polynomial ring $\mathcal{S}=\mathbb{K}[x_1,\dots,x_{n+1},x_{v_0},\dots,x_{v_{k-2}},x_{s_1},x_{s_2}]$ with the grading induced by $\tilde{\lambda}_{\mathcal{M}}$. Let $\lambda$ be the weight vector obtained by restricting $\tilde{\lambda}_{\mathcal{M}}$ to $R_{n+1}$, i.e. $\lambda(x_j)=\tilde{\lambda}_{\mathcal{M}}(x_j)$ for all $j$ from $1$ to $n+1$. From the above mentioned properties of $\tilde{\lambda}_{\mathcal{M}}$ and Theorem \ref{parcycle to defining ideal}, the ideal $\mathcal{M}\subseteq {\rm{in}}_{\lambda}(P_n)$.
Let $\prec$ be a monomial order such that ${\rm{in}}_{\prec}(P_n)=\mathcal{M}$ (\cite[Theorem 4.9]{conca2007contracted}). Since ${\rm{in}}_{\prec}(P_n)=\mathcal{M}\subseteq {\rm{in}_{\lambda}}(P_n)$, the ideal $\mathcal{M}={\rm{in}_{{\prec}_{\lambda}}}(P_n)$ where ${{\prec}_{\lambda}}$ is the monomial order obtained from the weight vector $\lambda$ followed by the monomial order $\prec$ on $R_{n+1}$ (\cite[Proposition 2.2.6]{herzog2011monomial}). The set of all $2\times 2$ minors of the matrix $T_n$ (Equation (\ref{matTn})) forms a Gr\"obner basis of $P_n$ since the set of the initial terms of these minors with respect to ${\prec}_{\lambda}$ generates $\mathcal{M}$ (Definition \ref{gbdef}). Let $P_{n,t}$ be the homogenisation of $P_n$ with respect to $\lambda$. From Corollary \ref{specialfibofgb} and the proof of \cite[Theorem 3.1.2]{herzog2011monomial}, the special fibre (fibre at $t=0$) of the Gr\"obner degeneration $R_{n+1}[t]/P_{n,t}$ of $R_{n+1}/P_n$ with respect to $\lambda$ is the quotient ring $R_{n+1}/\mathcal{M}$.

 Consider the complex $\tilde{\mathcal{F}}_{t}(C_{\Pi,\mathcal{M}})=\mathcal{F}_{t}\bigotimes \mathcal{S}[t]/\langle \mathcal{A} \rangle$ of $R_{n+1}[t]$-modules. Since $\mathcal{F}_{t}$ is exact and minimal and 
$\mathcal{A}$ is a regular sequence on $\mathcal{S}[t]$ and $\mathcal{S}[t]/I_{\Pi,\mathcal{M},t}$, the complex $\tilde{\mathcal{F}}_{t}(C_{\Pi,\mathcal{M}})$ is exact (Corollary \ref{Corollary 5.3}). The minimality of $\tilde{\mathcal{F}}_{t}(C_{\Pi,\mathcal{M}})$ follows from the minimality of the complex $\mathcal{F}_{t}$ and the homogeneity of the sequence $\mathcal{A}$. Since $\lambda$ is the restriction of $\tilde{\lambda}_{\mathcal{M}}$ to $R_{n+1}$ and $\mathcal{A}$ is a homogeneous sequence on $\mathcal{S}$ with the grading induced by $\tilde{\lambda}_{\mathcal{M}}$, the weights $\lambda(\psi(x_{l}))=\tilde{\lambda}_{\mathcal{M}}(x_{l})$ for all the variables $x_{l}$ of $\mathcal{S}$ (recall $\psi$ from Theorem \ref{parcycle to defining ideal}). 
From Theorem \ref{construction of parking ideal}, the set $\{f_{S}=\mathbf{x}^{S\rightarrow \overline{S}}-\mathbf{x}^{\overline{S}\rightarrow S}\mid S\text{ is a connected-parset of}~C_{\Pi,\mathcal{M}} ~ \text{not}\newline \text{containing}~ \{s_1,s_2\} \}$ generates $I_{\Pi,\mathcal{M}}$ and ${\rm{in}}_{\tilde{\lambda}_{\mathcal{M}}}(f_{S})=\mathbf{x}^{S\rightarrow \overline{S}}$. Let $f_{S}^{h}$ be the homogenisation of $f_{S}$ with respect to $\tilde{\lambda}_{\mathcal{M}}$. Consider the map $\bar{\psi}:\mathcal{S}[t]\rightarrow R_{n+1}[t]$ defined as $\bar{\psi}(t)=t$ and $\bar{\psi}(x_{l})=\psi(x_{l})$ for all the variables $x_{l}$ of $\mathcal{S}$. Since $\tilde{\lambda}_{\mathcal{M}}(x_{l})=\lambda(\psi(x_{l}))$ for all the variables $x_{l}$ of $\mathcal{S}$, $\bar{\psi}(f_{S}^{h})=\psi(f_{S})^{h}$ where $\psi(f_{S})^{h}$
is the homogenisation of $\psi(f_{S})$ with respect to $\lambda$. Hence, from Theorem \ref{parcycle to defining ideal}, $ \mathcal{S}[t]/I_{\Pi,\mathcal{M},t} \bigotimes \mathcal{S}[t]/\langle \mathcal{A}\rangle \simeq R_{n+1}[t]/P_{n,t}$. This shows that the complex $\tilde{\mathcal{F}}_{t}(C_{\Pi,\mathcal{M}})$ is a minimal free resolution of $R_{n+1}[t]/P_{n,t}$.
\end{proof}
We compute the (weighted) Hilbert series of $R_{n+1}[t]/P_{n,t}$ where $P_{n,t}$ is the Gr\"obner degeneration of $P_n$ to $\mathcal{M}$ and $\mathcal{M}$ is a Cohen-Macaulay initial monomial ideal of $P_n$. Let $\lambda_{\mathcal{M}}$ be the weight vector associated with $\mathcal{M}$ as defined in Section \ref{weights} and $\tilde{\lambda}_{\mathcal{M}}$ be the extension of $\lambda_{\mathcal{M}}$ as discussed in the proof of Theorem \ref{Grodegminres_theo}.
Let
${\rm{Acy}}_{\mathcal{M}}(i)$ be the set of all $\{s_1,s_2\}$-acyclic $(i+1)$-partitions of $C_{\Pi,\mathcal{M}}$ where $C_{\Pi,\mathcal{M}}$ is the parcycle associated with $\mathcal{M}$. Recall that for any $\mathcal{C}\in {\rm{Acy}}_{\mathcal{M}}(i)$, the vector
$D(\mathcal{C})=\sum_{l\in V(C_{\Pi,\mathcal{M}})}{\rm{out}}_{\mathcal{C}}(l)l$ where $V(C_{\Pi,\mathcal{M}})$ is the set of the basic vertices of $C_{\Pi,\mathcal{M}}$.

\begin{corollary}\label{weiHilGrob}
    Let $P_{n,t}$ be the Gr\"obner degeneration of $P_n$ to $\mathcal{M}$ and $\lambda$ be the weight vector obtained from restricting $\tilde{\lambda}_{\mathcal{M}}$ on $R_{n+1}$. The weighted Hilbert series of $R_{n+1}[t]/P_{n,t}$ corresponding to the weight vector $(\lambda,1)$ is 
$$
        H(R_{n+1}[t]/P_{n,t};z)=\frac{1+\sum_{i=1}^{n-1}{\sum_{\mathcal{C}\in{\rm{Acy}}_{\mathcal{M}}(i)} {(-1)^{i}       z^{\lambda_{\mathcal{M}} \cdot D(\mathcal{C})}}}}{(1-z)\prod_{i=1}^{n+1}(1-z^{\lambda(x_i)})}.
   $$
\end{corollary}
\begin{proof}
    From Theorem \ref{Grodegminres_theo}, the complex $\tilde{\mathcal{F}}_{t}(C_{\Pi,\mathcal{M}})$ is a minimal free resolution of $R_{n+1}[t]/P_{n,t}$. Hence, $H(R_{n+1}[t]/P_{n,t};z)=H(R_{n+1}[t];z)+\sum_{i=1}^{n-1}(-1)^{i}H(\tilde{F}_{t,\Pi,i};z)$ where $\tilde{F}_{t,\Pi,i}$ (Equation (\ref{consFtil})) is the free module of $\tilde{\mathcal{F}}_{t}(C_{\Pi,\mathcal{M}})$ at the homological degree $i$.
\end{proof}

\subsection{Multigraded Betti Numbers}\label{mulgrabnum}
Recall that in Section \ref{construction of parcycles}, we associated a parcycle with each Cohen-Macaulay initial monomial ideal $\mathcal{M}$ of $P_n$. In this section, we derive formulas for the multigraded Betti numbers of $\mathcal{M}$ and their support in terms of certain connected partitions of the parcycle associated with $\mathcal{M}$. We conclude with a discussion on using these formulas to compute the (local) Hilbert series of certain lex-segment ideals.

Let $\mathcal{M}$ be a Cohen-Macaulay initial monomial ideal of $P_n$. Let $1=i_0<i_1<\dots<i_{k}=n+1$ be the sequence associated with $\mathcal{M}$ and let $C_{\Pi,\mathcal{M}}$ be the corresponding parcycle of $\mathcal{M}$. 
From Theorem \ref{complex of G-parking ideal}, the complex $\mathcal{F}_{0,\Pi}$ is a $\mathbb{Z}^{n+k}$-graded minimal free resolution of the quotient ring of the $G$-parking function ideal $M_{\Pi,\mathcal{M}}$ of $C_{\Pi,\mathcal{M}}$. Recall that the complex $\mathcal{F}_{0,\Pi}$ is a complex of free modules over the polynomial ring $\tilde{\mathcal{S}}$ whose variables correspond to the basic vertices of $C_{\Pi,\mathcal{M}}$. Let $\mathcal{F}_0:=\mathcal{F}_{0,\Pi} \otimes_{\tilde{\mathcal{S}}} \mathcal{S}$ where $\mathcal{S}=\mathbb{K}[x_1,\dots,x_{n+1},x_{v_0},\dots,x_{v_{k-2}}, x_{s_1}, x_{s_2}]$. From Theorem \ref{regularity of A in R/M lemma} and Lemma \ref{lemma parkingtensorint}, the complex $\tilde{\mathcal{F}_0}(C_{\Pi,\mathcal{M}})=\mathcal{F}_0 \otimes_{\mathcal{S}} \mathcal{S}/\langle \mathcal{A}\rangle$ is a minimal free resolution of $R_{n+1}/\mathcal{M}$ where $R_{n+1}=\mathbb{K}[x_1,\dots,x_{n+1}]$ and
$\mathcal{A}=(x_{v_j}-x_{i_{j+2}},x_{s_1}-x_{i_{k-1}},x_{s_2}-x_{i_1}\mid 0\leq j\leq k-2)$. The resolution $\tilde{\mathcal{F}_0}(C_{\Pi,\mathcal{M}})$ is $\mathbb{Z}^{n+1}$-graded since it is obtained from the resolution $\mathcal{F}_{0}$ by replacing the variables $x_{v_j}$ with $x_{i_{j+2}}$ for all $j$ from $0$ to $k-2$, and $x_{s_1}, x_{s_2}$ with $x_{i_{k-1}}$ and ${x_{i_1}}$, respectively.
\begin{remark}\label{celres}\rm{The cellular resolution property of the complex $\tilde{\mathcal{F}_0}(C_{\Pi,\mathcal{M}})$ is inherited from  $\mathcal{F}_{0,\Pi}$.
    }\qed
\end{remark}
Let ${\rm{Acy}}_{\mathcal{M}}(i)$ be the set of all $\{s_1,s_2\}$-acyclic $(i+1)$-partitions of $C_{\Pi,\mathcal{M}}$. Note that the complex $\mathcal{F}_{0}$ is a minimal free resolution of $\mathcal{S}/M_{\Pi,\mathcal{M}}$.
The free module of $\mathcal{F}_{0}$ at the homological degree $i$ is given by $F_{0,\Pi,i}\otimes \mathcal{S}=\bigoplus_{\mathcal{C}}\mathcal{S}(-\tilde{D}(\mathcal{C}))$ where $\mathcal{C}\in {\rm{Acy}}_{\mathcal{M}}(i)$ and $\tilde{D}(\mathcal{C})\in \mathbb{Z}^{n+k+2}$. The vector $\tilde{D}(\mathcal{C})$ can be thought of as a $\mathbb{Z}$-valued function on the set of variables of $\mathcal{S}$ as follows:
\begin{equation*}
\tilde{D}(\mathcal{C})(x_l)=
\begin{cases}
    {\rm{out}}_{\mathcal{C}}(l),~\text{if}~ l~\text{is a basic vertex of}~C_{\Pi,\mathcal{M}},\\
        
        0,~\text{otherwise}.
    \end{cases}
\end{equation*}
\begin{proposition}
    Let $\mathcal{M}$ be a Cohen-Macaulay initial monomial ideal of $P_n$ and $C_{\Pi,\mathcal{M}}$ be the parcycle associated with it. The Betti number $\beta_{i,{\bf a}}(\mathcal{S}/M_{\Pi,\mathcal{M}})$ at the homological degree $i$ with support ${\bf a}\in \mathbb{Z}^{n+k+2}$ is the total number of $\{s_1,s_2\}$-acyclic $(i+1)$-partitions $\mathcal{C}$ of $C_{\Pi,\mathcal{M}}$ for which $\tilde{D}(\mathcal{C})={\bf a}$.
\end{proposition}
\begin{proof}
    The proof of the proposition follows from the above discussion on the degrees of the basis elements of $\mathcal{F}_0$.
\end{proof}
By keeping track of the supports of the multigraded Betti numbers of $\mathcal{S}/M_{\Pi,\mathcal{M}}$ as we go modulo the regular sequence $\mathcal{A}$, we obtain corresponding formulas for $R_{n+1}/\mathcal{M}$. By definition, the complex $\tilde{\mathcal{F}_0}(C_{\Pi,\mathcal{M}})=\mathcal{F}_0 \otimes_{\mathcal{S}} \mathcal{S}/\langle \mathcal{A}\rangle$. By examining the form of $\mathcal{A}$ and $\tilde{D}(\mathcal{C})$,
the free module of $\tilde{\mathcal{F}_0}(C_{\Pi,\mathcal{M}})$ at the homological degree $i$ is given by $\bigoplus_{\mathcal{C}}R_{n+1}(-\bar{D}(\mathcal{C}))$ where $\mathcal{C}\in {\rm{Acy}}_{\mathcal{M}}(i)$ 
and $\bar{D}(\mathcal{C})\in \mathbb{Z}^{n+1}$ is as follows:
\begin{equation}\label{mgradbetm}
\bar{D}(\mathcal{C})(x_{l})=
    \begin{cases}
    \tilde{D}(\mathcal{C})(x_l),~\text{if}~l\in \{1,\dots,n+1\}\setminus\{i_2,\dots,i_k\},\\

    \tilde{D}(\mathcal{C})(x_l)+{\rm{out}}_{\mathcal{C}}(v_{j-2}),~\text{if}~l=i_j~\text{and}~ j \in \{2,\dots,k\}.
    \end{cases}
\end{equation}
In the following, the sets $\{e_1,\dots,e_{n+1}\}$ and $\{e_1,\dots,e_{n+1},e_{v_0},\dots,e_{v_{k-2}},e_{s_1},e_{s_2}\}$ are the standard bases of the lattices $\mathbb{Z}^{n+1}$ and $\mathbb{Z}^{n+k+2}$, respectively. Taking cue from the regular sequence $\mathcal{A}$ and $\tilde{D}(\mathcal{C})$, we define the $\mathbb{Z}$-linear map $$\eta:\mathbb{Z}^{n+k+2}\rightarrow \mathbb{Z}^{n+1}$$
$$\eta(e_i)= e_i, 1\leq i\leq n+1,$$
$$\eta(e_{v_{j}})= e_{i_{j+2}}, 0\leq j\leq k-2,$$
$$\eta(e_{s_1})=0=\eta(e_{s_2}).$$
The set $\{e_{i_{j+2}}-e_{v_j},e_{s_1},e_{s_2}\mid 0\leq j\leq k-2\}$ is a basis of ${\rm{ker}}(\eta)$.

\begin{theorem}\label{Betdeg}
Let $\mathcal{M}$ be a Cohen-Macaulay initial monomial ideal of $P_n$ and $
C_{\Pi,\mathcal{M}}$ be the parcycle associated with it. 
A vector ${\bf b}\in \mathbb{Z}^{n+1}$ is the support of a multigraded Betti number $\beta_{i,{\bf b}}$ of $R_{n+1}/\mathcal{M}$ if and only if ${\bf b}=\bar{D}(\mathcal{C})$ where $\mathcal{C}\in {\rm{Acy}}_{\mathcal{M}}(i)$. Moreover, the Betti number $\beta_{i,{\bf b}}$ of $R_{n+1}/\mathcal{M}$ corresponding to the support ${\bf b} \in \mathbb{Z}^{n+1}$ is
    
    \[ \beta_{i,{\bf b}}(R_{n+1}/\mathcal{M})=\sum_{{\bf a}\in\eta^{-1}({\bf b})}\beta_{i,{\bf a}}(\mathcal{S}/M_{\Pi,\mathcal{M}}).
    \]
\end{theorem}
\begin{proof}
    The first part of the theorem follows from the above discussion on the degrees $\bar{D}(\mathcal{C})$ of the basis elements of the complex $\tilde{\mathcal{F}}_0(C_{\Pi,\mathcal{M}})$. The second part follows from the observation that for any two $\mathcal{C}_1,\mathcal{C}_2\in {\rm{Acy}}_{\mathcal{M}}(i)$, the vectors
    $\bar{D}(\mathcal{C}_1)=\bar{D}(\mathcal{C}_2)$ if and only if 
    $\tilde{D}(\mathcal{C}_1)- \tilde{D}(\mathcal{C}_2)\in {\rm{ker}}(\eta)$. 
\end{proof}

\begin{example}\rm{Consider the Cohen-Macaulay initial monomial ideal $\mathcal{M}$ of $P_4$ defined in Example \ref{cminmidp4}. The parcycle $C_{\Pi,\mathcal{M}}$ associated with $\mathcal{M}$ is depicted in Figure \ref{paroeg}. Let $\mathcal{C}_1, \mathcal{C}_2$ and $\mathcal{C}_3$ be the $\{s_1,s_2\}$-acyclic $4$-partitions of $C_{\Pi,\mathcal{M}}$ (illustrated in Figure \ref{acyclic partitions}) with sources at $\{1,v_0\},\{3\}$ and $\{2,v_1\}$, respectively. The supports ${\bf b}\in \mathbb{Z}^{5}$ of the nonzero multigraded Betti numbers $\beta_{3,{\bf b}}$ of $\mathbb{K}[x_1,\dots,x_5]/\mathcal{M}$ are given by $\bar{D}(\mathcal{C}_1)=(1,0,1,1,1), \bar{D}(\mathcal{C}_2)=(1,0,2,0,1)$ and $\bar{D}(\mathcal{C}_3)=(1,1,1,0,1)$. Moreover, $\beta_{3,\bar{D}(\mathcal{C}_i)}=1$ for all $i$.
    }\qed
    \end{example}

\subsubsection{Hilbert Series of Lex-Segment Ideals}
Consider the lexicographic order $>_{{\rm lex}}$ on $\mathbb{K}[x,y]$ induced by the ordering $x>y$. 
Let ${\bf m}=\langle x,y\rangle$ be the homogeneous maximal ideal of $\mathbb{K}[x,y]$.
\begin{definition}
A monomial ideal $L$ of $\mathbb{K}[x,y]$ is called a lex-segment ideal if $q(x,y)\in L$ and $p(x,y)\in \mathbb{K}[x,y]$ are monomials of degree $d$ such that $p(x,y)>_{{\rm lex}}q(x,y)$, then $p(x,y)\in L$.
    \end{definition}

    Lex-segment ideals naturally appear in Zariski's theorem on the factorisation of contracted ideals $I$ of $\mathbb{K}[x,y]$ where $\mathbb{K}$ is an algebraically closed field (\cite[Theorem 3.3]{conca2007contracted}). Recall that an ${\bf m}$-primary ideal $I$ of $\mathbb{K}[x,y]$ is called a \emph{contracted ideal} if there exists a linear form $z\in {\bf m}$ such that $I=J \cap \mathbb{K}[x,y]$ where $J=I (\mathbb{K}[x,y][{\bf m}/z])$.
    The lex-segment ideals $L_i$ that show up in the factorisation of $I$ have several characteristics in common with $I$.
For instance, the Cohen-Macaulayness of the associated graded ring of a contracted ideal $I=L_1 \cdots L_s$ is equivalent to the Cohen-Macaulayness of the associated graded ring of $L_i$ for all $i$ (\cite[Corollary 3.14]{conca2005graded}).

Let $L$ be an ${\bf{m}}$-primary lex-segment ideal of $\mathbb{K}[x,y]$. We say the lex-segment ideal $L$ is of \emph{ initial degree} $n$ if $n={\rm min}\{m\mid x^m\in L\}$. For each lex-segment ideal of initial degree $n$, we associate a vector $w=(w_0(L),\dots,w_{n}(L))\in \mathbb{Z}^{n+1}$ to $L$ where $w_i(L)={\rm{min}}\{j\mid x^{n-i}y^{j}\in I\}$. The entries of $w$ form an increasing sequence, i.e. $0=w_0(L)\leq w_1(L)\leq \cdots\leq w_n(L)$.  Let $\mathcal{M}$ be a Cohen-Macaulay initial monomial ideal of $P_n$ and $C({\mathcal{M}})=\{{\bf b}\in \mathbb{Q}^{n+1}_{\geq 0} \mid {\rm{in}_{\tau}}({\rm{in}}_{{\bf b}}(P_n))=\mathcal{M}
\}$ where ${\tau}$ is a given term order such that ${\rm{in}_{\tau}}(P_n)=\mathcal{M}$. From \cite[Corollary 6.4]{conca2007contracted}, if the vector $w$ associated with $L$ lies in $C({\mathcal{M}})$, then 

\begin{align*}
    H^{1}_{L}(z) &= \sum_{j \geq 0} \operatorname{dim}_{\mathbb{K}}(\mathbb{K}[x,y]/L^{j+1})z^{j} \\
    &= |w|\frac{1+(n-1)z}{(1-z)^3} + w \cdot \sum_{i \geq 1}(-1)^{i}\sum_{\mathcal{C} \in \operatorname{Acy}_{\mathcal{M}}(i)}\beta_{i,\bar{D}(\mathcal{C})}(R_{n+1}/\mathcal{M})\frac{z^{|\bar{D}(\mathcal{C})|-1}}{(1-z)^{n+1}}\bar{D}(\mathcal{C})
\end{align*}
where $\bar{D}(\mathcal{C})\in \mathbb{Z}^{n+1}$ (defined in Equation (\ref{mgradbetm})). In the same vein, the $h$-polynomial and the Hilbert polynomial of $L$ can be computed via the $\mathbb{Z}^{n+1}$-graded Betti numbers of $\mathcal{M}$ if the vector $w$ associated with $L$ lies in $C({\mathcal{M}})$ (\cite[Proposition 6.4]{conca2007contracted}).

\section{Appendix}\label{Appendix}
\appendix{}

 In Appendix \ref{Eagon northcott complex and toppling complex}, we relate the minimal free resolution of the toppling ideal and the Eagon-Northcott complex, both associated with the rational normal curve.

\section{Eagon-Northcott Complex and the Minimal Free Resolution of the Toppling Ideal}\label{Eagon northcott complex and toppling complex}

Recall that in Section \ref{mainthmpf}, we constructed a minimal free resolution of the rational normal curve (viewed as the toppling ideal of a parcycle). In the following, we explicitly construct an isomorphism between the minimal free resolution $\mathcal{F}_{1,\Pi}$ where $\Pi=\{\{1,n+1\}, \{2\}, \dots,\{n\}\}$ and the well-known Eagon-Northcott complex that also minimally resolves the rational normal curve. Note that such an isomorphism exists by general considerations.

We know that the defining ideal $P_{n}$ of $\Gamma_n$ is a determinantal ideal generated by the set of all $2\times 2$ minors of the matrix $T_n$ (Equation (\ref{matTn})).
From \cite[Theorem 6.4]{eisenbud2005geometry}, the Eagon-Northcott complex minimally resolves $R_{n+1}/P_{n}$ where $R_{n+1}=\mathbb{K}[x_1,\dots,x_{n+1}]$. First, we will specialise the Eagon-Northcott complex \cite{EAGONNORTHCOTT}, which resolves more general determinantal varieties to the case of the rational normal curve. 
Let $y_1,\dots,y_n,z_1,z_2$ be indeterminates over $R_{n+1}$. Let $G = \bigoplus_{i=1}^n R_{n+1} y_i$ be a free module with $\{y_i\}^{n}_{i=1}$ as its basis and $\bigwedge G$ be the exterior algebra over $G$.

Each row of the matrix $T_n$ defines a homomorphism on $\bigwedge G$. Let $\Delta_1$ and $\Delta_2$ be the homomorphisms defined by the first and second rows of $T_n$, respectively. The map $\Delta_k:\bigwedge G \longrightarrow \bigwedge G$ is defined as follows: $$\Delta_{k}(y_{i_{1}}\wedge y_{i_{2}}\wedge \cdots \wedge y_{i_{n}}) = \sum_{p=1}^{n} (-1)^{p+1}(T_n)_{k, i_p} (y_{i_1}\wedge \cdots \wedge \hat{y}_{i_p}\wedge \cdots \wedge y_{i_n}).
$$
Let $\phi_q$ be the free $R_{n+1}$-module with $\{z^{v_1}_1z^{v_2}_2\mid v_1,v_2\in \mathbb{N}, v_1+v_2=q\}$ as its basis.
The Eagon-Northcott resolution $E^{A}$ of $P_{n}$ is defined as
$$
E^A:0\longrightarrow E_{n-1}^{A} \overset{d}\longrightarrow E_{n-2}^{A}\overset{d}\longrightarrow \dots \overset{d} \longrightarrow E_{1}^{A} \overset{d} \longrightarrow E_{0}^{A} \longrightarrow 0,
$$ where $E_{q+1}^A = \bigwedge^{2+q} G \otimes \phi_q$ is a free $R_{n+1}$-module, for all $0 \leq q \leq n-2$. The differentials are defined as follows: $d(y_{i_1}\wedge\cdots \wedge  y_{i_{q+2}} \otimes z_1^{v_1}z_2^{v_2} )=$
\begin{equation*}
   \begin{cases}
              
                  \Delta_1(y_{i_1}\wedge\cdots \wedge y_{i_{q+2}}) \otimes z_1^{v_1-1}z_2^{v_2}+  \Delta_2(y_{i_1}\wedge\cdots \wedge y_{i_{q+2}}) \otimes z_1^{v_1}z_2^{v_2-1}            \\
               \qquad \qquad \qquad \qquad \qquad \qquad \qquad\quad, \text{if}~ v_1>0,v_2>0,q>0,
               \\
                \Delta_1(y_{i_1}\wedge\cdots \wedge y_{i_{q+2}}) \otimes z_1^{v_1-1}z_2^{v_2}              \qquad,~ \text{if}~ v_1>0,v_2 =0,q>0,\\
               \Delta_2(y_{i_1}\wedge\cdots\wedge y_{i_{q+2}}) \otimes z_1^{v_1}z_2^{v_2-1} \qquad,~ \text{if}~v_1=0,v_2>0,q>0,
           \end{cases}
           \end{equation*}

and $d(y_{i_1}y_{i_2}\otimes1) = x_{i_1}x_{i_2+1}-x_{i_2}x_{i_1+1}$.

Now, we define a homomorphism of complexes $\Psi$ between the Eagon-Northcott resolution $E^A$ and the minimal free resolution $\mathcal{F}_{1,\Pi}$ of $P_{n}$. The homomorphism $\Psi$ is defined as follows:
\begin{equation}\label{eq9}
    \begin{tikzcd}[cells={nodes={minimum height=2em}}]
0\arrow[r] & E_{n-1}^A \arrow[r,"d"] \arrow[d,"\Psi"] &   E_{n-2}^A \arrow[r,"d"] \arrow[d,"\Psi"] &  \dots \arrow[r] & E_{1}^A\arrow[r,"d"]\arrow[d,"\Psi"] & R_{n+1}\arrow[r,"d"]\arrow[d,"Id"]& 0\\
0\arrow[r]&F_{1,\Pi,n-1}  \arrow[r,"\delta"]&F_{1,\Pi,n-2}  \arrow[r,"\delta"]&\dots \arrow[r]&F_{1,\Pi,1}  \arrow[r,"\delta"] &R_{n+1} \arrow[r,"\delta"] & 0
\end{tikzcd},
\end{equation}
where the map $\Psi:E_{r}^A \longrightarrow F_{1,\Pi,r} $ is defined as follows:
$$\Psi(y_{i_1}\wedge \cdots \wedge y_{i_{r+1}} \otimes z_1^{v_1}z_2^{v_2}) = \mathcal{C}_{r+1}.$$ The element $y_{i_1}\wedge\cdots \wedge y_{i_{r+1}} \otimes z_1^{v_1}z_2^{v_2}$ is a basis element of $E_{r}^A$ and $\mathcal{C}_{r+1}=(\tilde{\Pi},\mathcal{O})$ is a $\{1,n+1\}$-acyclic $(r+1)$-partition of $(C_{n+1},\Pi)$ where $\tilde{\Pi}= \{\{i_1+1,\dots,i_2\}, \{i_2+1,\dots,i_3\},\dots,\{i_{r+1}+1,\dots,i_1\}\}$ and the source is at $V_i\in \tilde{\Pi}$ that contains the vertex $i_{v_2+1}+1$.  
\begin{theorem}
    The morphism of complexes $\Psi$ from the Eagon-Northcott resolution $E^A$ of $P_n$ to the minimal free resolution $\mathcal{F}_{1,\Pi}$  of $P_{n}$ is an isomorphism of complexes.
\end{theorem}
\begin{proof}

Since $\Psi$ is defined on the basis elements $y_{i_1}\wedge\cdots\wedge y_{i_{r+1}} \otimes z_1^{v_1}z_2^{v_2}$ of $E_{r}^A$ where $1 \leq i_1<i_2< \dots <i_{r+1} \leq n$ and $ v_1+v_2 =r-1$, the map
$\Psi:E_{r}^A \longrightarrow F_{1,\Pi,r}$ is a well-defined $R$-module homomorphism. As the modules $E_{r}^A$ and $F_{1,\Pi,r} $ have the same rank $r \cdot \binom{n}{r+1}$ and $\Psi$ maps each basis element of $E_{r}^A$ to a basis element of $F_{1,\Pi,r} $, we conclude that the map $\Psi:E_{r}^A \longrightarrow F_{1,\Pi,r} $ is an isomorphism for all $1 \leq r \leq n-1$. 

We now prove that the diagram (in Equation (\ref{eq9isocomp})) is commutative, i.e. $\Psi_{r} d_{r+1} = \delta_{r} \Psi_{r+1}$ for all $r\geq0$. To prove the commutativity of the diagram, it is enough to prove that the maps $\Psi_{r} d_{r+1}:E_{r+1}^{A}\rightarrow F_{1,\Pi,r}$ and $\delta_{r} \Psi_{r+1}: E_{r+1}^{A}\rightarrow F_{1,\Pi,r}$ agree on the basis elements of $E_{r+1}^{A}$ for all $r\geq 0$. We know that the basis elements of $E_{r+1}^{A}$ are of the form $y_{i_1}\wedge \cdots \wedge y_{i_{r+2}} \otimes z_1^{v_1}z_2^{v_2}$, where $1\leq i_1<i_2<\dots<i_{r+2}\leq n$ and $v_1,v_2\in \mathbb{N}$ satisfy $v_1+v_2=r$.
\begin{itemize}
    \item {\bf Case 1:}\label{ca1} We assume that $v_1=0$. The basis element $y_{i_1}\wedge \cdots \wedge y_{i_{r+2}} \otimes z_1^{v_1}z_2^{v_2}$ can be written as $y_{i_1}\wedge \cdots \wedge y_{i_{r+2}} \otimes z_2^{v_2}$, where $1\leq i_1<i_2<\dots<i_{r+2}\leq n$ and $v_2=r$. The image of the basis element under the map $\Psi_{r} d_{r+1}$ is 

    \begin{align*}
\Psi_{r} d_{r+1}(y_{i_1} \wedge \cdots \wedge y_{i_{r+2}} \otimes z_2^{r})
&= \Psi_{r}(\Delta_2(y_{i_1} \wedge \cdots \wedge y_{i_{r+2}}) \otimes z_2^{r-1}) \\
&= \sum_{p=1}^{r+2} (-1)^{p+1} x_{i_p+1} \Psi_{r}(y_{i_1} \wedge \cdots \wedge \hat{y}_{i_p} \wedge \cdots \wedge y_{i_{r+2}}\otimes z_2^{r-1})
\end{align*}
Before defining the image of $y_{i_1} \wedge \cdots \wedge \hat{y}_{i_p} \wedge \cdots \wedge y_{i_{r+2}}\otimes z_2^{r-1}$ under $\Psi_{r}$, we reindex the subscripts of $y$. The basis element $y_{i_1} \wedge \cdots \wedge \hat{y}_{i_p} \wedge \cdots \wedge y_{i_{r+2}}\otimes z_2^{r-1}=y_{j_1} \wedge \cdots \wedge y_{j_{r+1}}\otimes z_2^{r-1}$ where 
\begin{equation}\label{eq85}
y_{j_l} =
\begin{cases} 
y_{i_l}, & \text{if } 1 \leq l < p; \\
y_{i_{l+1}}, & \text{if } p \leq l \leq r+1.
\end{cases}
\end{equation}
The map $\Psi_r :E_{r}^{A}\rightarrow F_{1,\Pi,r}$ maps $y_{i_1} \wedge \cdots \wedge \hat{y}_{i_p} \wedge \cdots \wedge y_{i_{r+2}}\otimes z_2^{r-1}$ to $[{\mathcal{C}}_{r+1,p}]$ where ${\mathcal{C}}_{r+1,p}$ is the $\{1,n+1\}$-acyclic $(r+1)$-partition of $(C_{n+1},\Pi)$ whose main vertices are $V_1=\{j_1+1,\dots,j_2\},V_2= \{j_2+1,\dots,j_3\},\dots,V_{r+1}=\{j_{r+1}+1,\dots,j_1\}$ and the source is at the main vertex $V_{r}$.

Let $[\mathcal{C}_{r+2}]$ be the image of $y_{i_1} \wedge \cdots \wedge y_{i_{r+2}} \otimes z_2^{r}$ under $\Psi_{r+1}$. From the construction of the map $\Psi_{r+1}:E_{r+1}^{A}\rightarrow F_{1,\Pi,r+1}$, $\mathcal{C}_{r+2}$ is the $\{1,n+1\}$-acyclic $(r+2)$-partition of $(C_{n+1},\Pi)$ whose main vertices are $V'_1=\{i_1+1,\dots,i_2\},V'_2= \{i_2+1,\dots,i_3\},\dots,V'_{r+2}=\{i_{r+2}+1,\dots,i_1\}$ and the source is at the main vertex $V'_{r+1}$. Note that the set of contractible edges in $[\mathcal{C}_{r+2}]$ is  $\{f_1=(V^{'}_1,V^{'}_{r+2}),f_{l+1}=(V^{'}_{l+1},V^{'}_{l})\mid 1\leq l \leq r+1\}$. The image of $y_{i_1} \wedge \cdots \wedge y_{i_{r+2}} \otimes z_2^{r}$ under $\delta_{r}\Psi_{r+1}$ is
\begin{align*}
\delta_{r}\Psi_{r+1}(y_{i_1} \wedge \cdots \wedge y_{i_{r+2}} \otimes z_2^{r}) &= \delta_{r}([\mathcal{C}_{r+2}]) \\
&= \sum_{p=1}^{r+2} {\rm sign}_{[\mathcal{C}_{r+2}]}(f_p) m_{[\mathcal{C}_{r+2}]}(f_p) \cdot ([\mathcal{C}_{r+2}]/f_p) \\
&= \sum_{p=1}^{r+2} {\rm sign}_{[\mathcal{C}_{r+2}]}(f_p) x_{i_p+1} \cdot [{\mathcal{C}}_{r+1,p}].
\end{align*}

\item {\bf Case 2:} Let $v_2=0$. In this case, the basis element $y_{i_1}\wedge \cdots \wedge y_{i_{r+2}} \otimes z_1^{v_1}z_2^{v_2}=y_{i_1}\wedge \cdots \wedge y_{i_{r+2}} \otimes z_1^{v_1}$, where $1\leq i_1<i_2<\dots<i_{r+2}\leq n$ and $v_1=r$. The image of the basis element under the map $\Psi_{r} d_{r+1}$ is \begin{align*}
&\Psi_{r} d_{r+1}(y_{i_1} \wedge \cdots \wedge y_{i_{r+2}} \otimes z_1^{r})\\&
= \Psi_{r}(\Delta_1(y_{i_1} \wedge \cdots \wedge y_{i_{r+2}}) \otimes z_1^{r-1}) \\
&= \sum_{p=1}^{r+2} (-1)^{p+1} x_{i_p} \Psi_{r}(y_{i_1} \wedge \cdots \wedge \hat{y}_{i_p} \wedge \cdots \wedge y_{i_{r+2}}\otimes z_1^{r-1})
\end{align*}
Similar to Case 1,
we reindex the subscripts of $y$. The basis element $y_{i_1} \wedge \cdots \wedge \hat{y}_{i_p} \wedge \cdots \wedge y_{i_{r+2}}\otimes z_1^{r-1}=y_{j_1} \wedge \cdots \wedge y_{j_{r+1}}\otimes z_1^{r-1}$ where $y_{j_l}$ are as in Equation \ref{eq85}. 
The image of $y_{i_1} \wedge \cdots \wedge \hat{y}_{i_p} \wedge \cdots \wedge y_{i_{r+2}}\otimes z_1^{r-1}$ under $\Psi_r$ is $[\tilde{\mathcal{C}}_{r+1,p}]$ where $\tilde{\mathcal{C}}_{r+1,p}$ is the $\{1,n+1\}$-acyclic $(r+1)$-partition of $(C_{n+1},\Pi)$ whose main vertices are $V_1=\{j_1+1,\dots,j_2\},V_2= \{j_2+1,\dots,j_3\},\dots,V_{r+1}=\{j_{r+1}+1,\dots,j_1\}$ and the source is at the main vertex $V_{1}$.

Let $[\tilde{\mathcal{C}}_{r+2}]$ be the image of $y_{i_1} \wedge \cdots \wedge y_{i_{r+2}} \otimes z_1^{r}$ under $\Psi_{r+1}$. From the construction of the map $\Psi_{r+1}$, the acyclic partition $\tilde{\mathcal{C}}_{r+2}$ is the $\{1,n+1\}$-acyclic $(r+2)$-partition of $(C_{n+1},\Pi)$ whose main vertices are $V'_1=\{i_1+1,\dots,i_2\},V'_2= \{i_2+1,\dots,i_3\},\dots,V'_{r+2}=\{i_{r+2}+1,\dots,i_1\}$ and the source is at $V'_{1}$. The set of contractible edges in $[\tilde{\mathcal{C}}_{r+2}]$ is  $\{\tilde{f}_{1}=(V^{'}_{r+2},V^{'}_1)),\tilde{f}_{l+1}=(V^{'}_{l},V^{'}_{l+1})\mid 1\leq l \leq r+1\}$. 
The image of $y_{i_1} \wedge \cdots \wedge y_{i_{r+2}} \otimes z_1^{r}$ under $\delta_{r}\Psi_{r+1}$ is
\begin{align*}
\delta_{r}\Psi_{r+1}(y_{i_1} \wedge \cdots \wedge y_{i_{r+2}} \otimes z_1^{r}) &= \delta_{r}([\tilde{\mathcal{C}}_{r+2}]) \\
&= \sum_{p=1}^{r+2} {\rm sign}_{[\tilde{\mathcal{C}}_{r+2}]}(\tilde{f}_p) m_{\tilde{\mathcal{C}}_{r+2}}(\tilde{f}_p) \cdot ([\tilde{\mathcal{C}}_{r+2}]/\tilde{f}_p) \\
&= \sum_{p=1}^{r+2} {\rm sign}_{[\tilde{\mathcal{C}}_{r+2}]}(\tilde{f}_p) x_{i_p} \cdot [\tilde{\mathcal{C}}_{r+1,p}].
\end{align*}

\item {\bf Case 3:} Let $v_1$ and $v_2$ be both nonzero. The basis elements of $E_{r+1}^{A}$ for which $v_1$ and $v_2$ both are nonzero  are of the form $y_{i_1}\wedge \cdots \wedge y_{i_{r+2}} \otimes z_1^{v_1}z_2^{v_2}$, where $1\leq i_1<i_2<\dots<i_{r+2}\leq n$ and $v_1,v_2\in \mathbb{N}$ satisfy $v_1+v_2=r$. Let $v_2=t$ where $1\leq t\leq r-1$. The image
\begin{equation}\label{eq86psid}
\begin{split}
&\Psi_{r} d_{r+1}\left(y_{i_1} \wedge \cdots \wedge y_{i_{r+2}} \otimes z_1^{v_1} z_2^{v_2}\right)\\& 
= \Psi_{r} \bigg( 
    \Delta_1 \left(y_{i_1} \wedge \cdots \wedge y_{i_{r+2}}\right) \otimes z_1^{r-t-1} z_2^{t} + \Delta_2 \left(y_{i_1} \wedge \cdots \wedge y_{i_{r+2}}\right) \otimes z_1^{r-t} z_2^{t-1} 
\bigg) \\
&= \sum_{p=1}^{r+2} (-1)^{p+1} x_{i_p} \Psi_{r} \left( y_{i_1} \wedge \cdots \wedge \hat{y}_{i_p} \wedge \cdots \wedge y_{i_{r+2}} \otimes z_1^{r-t-1} z_2^{t} \right) \\
& \quad + \sum_{p=1}^{r+2} (-1)^{p+1} x_{i_p+1} \Psi_{r} \left( y_{i_1} \wedge \cdots \wedge \hat{y}_{i_p} \wedge \cdots \wedge y_{i_{r+2}} \otimes z_1^{r-t} z_2^{t-1} \right) \\
&= \sum_{p=1}^{r+2} (-1)^{p+1} x_{i_p} \left[ \mathcal{C}^{(1)}_{r+1,p} \right] + \sum_{p=1}^{r+2} (-1)^{p+1} x_{i_p+1} \left[ \mathcal{C}^{(2)}_{r+1,p} \right]
\end{split}
\end{equation}
As in Case 1 and 2, we reindex the subscripts of $y$ and rewrite $y_{i_1} \wedge \cdots \wedge \hat{y}_{i_p} \wedge \cdots \wedge y_{i_{r+2}}$ as $y_{j_1} \wedge \cdots \wedge y_{j_{r+1}}$ where $y_{j_l}$ are as in Equation \ref{eq85}. The acyclic partitions $\mathcal{C}^{(1)}_{r+1,p}$ and $\mathcal{C}^{(2)}_{r+1,p}$ (in Equation \ref{eq86psid}) are $\{1,n+1\}$-acyclic $(r+1)$-partitions of $(C_{n+1},\Pi)$ with main vertices $V_1=\{j_1+1,\dots,j_2\},V_2= \{j_2+1,\dots,j_3\},\dots,V_{r+1}=\{j_{r+1}+1,\dots,j_1\}$ and the sources are at $V_{t+1}$ and $V_{t}$, respectively.

Let $[\mathcal{C}'_{r+2}]$ be the image of $y_{i_1} \wedge \cdots \wedge y_{i_{r+2}} \otimes z_1^{r-t} z_2^{t}$ under $\Psi_{r+1}$. The acyclic partition $\mathcal{C}'_{r+2}$ is a $\{1,n+1\}$-acyclic $(r+2)$-partition of $(C_{n+1},\Pi)$ with the main vertices $V'_1,\dots,V'_{r+2}$ (as defined in Case 1) and the source is at $V'_{t+1}$. Since $1\leq t\leq r-1$ and the source of $\mathcal{C}'_{r+2}$ is at $V'_{t+1}$, the set of contractible edges in $[\mathcal{C}'_{r+2}]$ is $\{f_1=(V^{'}_1,V^{'}_{r+2}),\tilde{f}_{1}=(V^{'}_{r+2},V^{'}_1)),f_{l+1}=(V^{'}_{l+1},V^{'}_{l}),\tilde{f}_{l+1}=(V^{'}_{l},V^{'}_{l+1})\mid 1\leq l \leq r+1\}$.
The image
\begin{equation}
\begin{split}
&\delta_{r}\Psi_{r+1}(y_{i_1} \wedge \cdots \wedge y_{i_{r+2}} \otimes z_1^{r-t} z_2^{t}) 
= \delta_{r}([\mathcal{C}'_{r+2}]) \\
&= \sum_{p=1}^{r+2} {\rm sign}_{[\mathcal{C}'_{r+2}]}(\tilde{f}_p) m_{[\mathcal{C}'_{r+2}]}(\tilde{f}_p) \cdot ([\mathcal{C}'_{r+2}]/\tilde{f}_p) \\
&\quad+ \sum_{p=1}^{r+2} {\rm sign}_{[\mathcal{C}'_{r+2}]}(f_p) m_{[\mathcal{C}'_{r+2}]}(f_p) \cdot ([\mathcal{C}'_{r+2}]/f_p) \\
&= \sum_{p=1}^{r+2} {\rm sign}_{[\mathcal{C}'_{r+2}]}(\tilde{f}_p) x_{i_p} \left[ \mathcal{C}^{(1)}_{r+1,p} \right] + \sum_{p=1}^{r+2} {\rm sign}_{[\mathcal{C}'_{r+2}]}(f_p) x_{i_p+1} \left[ \mathcal{C}^{(2)}_{r+1,p} \right].
\end{split}
\end{equation}
From the construction of the sign functions (discussed in Section \ref{preparcyc_subsect}), we can ensure that the sign functions of the Eagon-Northcott complex and the toppling resolution are compatible. From Cases 1, 2, and 3, the map $\psi$ is an isomorphism of complexes.\end{itemize}
\end{proof}

\footnotesize
\noindent {\bf Author's address:}

\smallskip

\noindent Department of Mathematics,\\
Indian Institute of Technology Bombay,\\
Powai, Mumbai,
India 400076.\\

\noindent {\bf Email id:}

Rahul Karki: rahulkarki877@gmail.com.

Madhusudan Manjunath: madhu@math.iitb.ac.in,  madhusudan73@gmail.com.\\

\begin{thebibliography}{9}


\bibitem{baker2013chip}
Matthew Baker and Farbod Shokrieh,\textsl{Chip-Firing Games, Potential Theory on Graphs, and Spanning Trees}, Journal of Combinatorial Theory, Series A {\bf 120.1} (2013), pp. 164--182.

\bibitem{bayer1998cellular}
Dave Bayer and Bernd Sturmfels, \textsl{Cellular Resolutions of Monomial Modules}, Journal f\"ur die Reine und Angewandte Mathematik {\bf 502}
(1998), pp. 123--140.

\bibitem{MR1896344}
Robert Cori, Dominique Rossin and Bruno Salvy, \textsl{Polynomial Ideals for Sandpiles and Their {G}r\"{o}bner Bases}, Theoretical Computer Science {\bf 276.1-2} (2002), pp. 1--15.

\bibitem{ds}
Anton Dochtermann and Raman Sanyal, \textsl{Laplacian Ideals, Arrangements, and Resolutions}, Journal of Algebraic Combinatorics {\bf 40.3}
(2014), pp. 805--822.


\bibitem{eisenbud2005geometry}
David Eisenbud, \textsl{The Geometry of Syzygies: A Second Course in Commutative Algebra and Algebraic Geometry}, Vol. 229. Springer, 2005.

\bibitem{greene1983interpretation}
Curtis Greene and Thomas Zaslavsky, \textsl{On the Interpretation of Whit-
ney Numbers Through Arrangements of Hyperplanes, Zonotopes, Non-
Radon Partitions, and Orientations of Graphs}, Transactions of the American Mathematical Society {\bf 280.1} (1983), pp. 97--126.

\bibitem{herzog2011monomial}
J{\"u}rgen Herzog and Takayuki Hibi, \textsl{Monomial Ideals}, Springer, 2011, pp. 3--22.

\bibitem{manjunathschwil}
Madhusudan Manjunath, Frank-Olaf Schreyer and John Wilmes, \textsl{Minimal Free Resolutions of the {$G$}-Parking Function Ideal and the Toppling Ideal}, Transactions of the American Mathematical Society {\bf 367.4}
(2015), pp. 2853–2874.

\bibitem{manjunath2013monomials}
Madhusudan Manjunath and Bernd Sturmfels, \textsl{Monomials, Binomials and Riemann--Roch}, Journal of Algebraic Combinatorics {\bf 37.4}
(2013), pp. 737–756.

\bibitem{miller2005combinatorial} Ezra Miller and Bernd Sturmfels, \textsl{Combinatorial Commutative Algebra}, Vol. 227. Springer Science \& Business Media, 2005.

\bibitem{mohammadi2016divisors} Fatemeh Mohammadi and Farbod Shokrieh, \textsl{Divisors on Graphs, Binomial and Monomial Ideals, and Cellular Resolutions}, Mathematische Zeitschrift {\bf 283.1} (2016), pp. 59--102.

\bibitem{stanley2004introduction} Richard P. Stanley, \textsl{An Introduction to Hyperplane Arrangements}, Geometric Combinatorics, IAS/Park
City Mathematics Series, {\bf 13.389-
496} (2007), p. 24.

\bibitem{conca2007contracted} Aldo Conca, Emanuela De Negri and Maria Evelina Rossi, \textsl{Contracted Ideals and the Gr{\"o}bner fan of the Rational Normal Curve}, Algebra \& Number Theory {\bf{1.3}} (2007), pp. 239--268.

\bibitem{conca2005graded} Aldo Conca, Emanuela De Negri, AV Jayanthan, and Maria Evelina Rossi, \textsl{Graded Rings Associated With Contracted Ideals}, Journal of Algebra {\bf 2.284} (2005), pp. 593--626.

\bibitem{de2023semigroup}
Hern{\'a}n de Alba Casillas, Daniel Duarte, and Ra{\'u}l Vargas Antuna, \textsl{A Semigroup Defining the Gr{\"o}bner Degeneration of a Toric Ideal}, Semigroup Forum, Springer (2023), pp. 1--20.

\bibitem{herzog2010binomial} J{\"u}rgen Herzog, Takayuki Hibi, Freyja Hreinsd{\'o}ttir, Thomas Kahle,
and Johannes Rauh, \textsl{Binomial Edge Ideals and Conditional Independence Statements}, Advances in Applied Mathematics {\bf 45.3} (2010), pp. 317--333. 

\bibitem{farkas2022minimal} Gavril Farkas and Eric Larson, \textsl{The Minimal Resolution Conjecture for Points on General Curves}, arXiv preprint arXiv:2209.11308 (2022).


\bibitem{li2022minimal}  Yupeng Li, Ezra Miller, and Erika Ordog, \textsl{Minimal Resolutions of Lattice Ideals}, arXiv preprint arXiv:2208.09557 (2022).


\bibitem{eagon2019minimal} John Eagon, Ezra Miller and Erika Ordog, \textsl{Minimal Resolutions of Monomial Ideals}, arXiv preprint arXiv:1906.08837 (2019).

\bibitem{o2018minimal} Liam O’Carroll and Francesc Planas-Vilanova, \textsl{Minimal Free Resolutions of Lattice Ideals of Digraphs}, Algebraic Combinatorics {\bf 1.2} (2018), pp. 283--326.

\bibitem{ene2011cohen} Viviana Ene, J{\"u}rgen Herzog, and Takayuki Hibi, \textsl{Cohen-Macaulay Binomial Edge Ideals}, Nagoya Mathematical Journal {\bf 204} (2011), pp. 57--68.

\bibitem{corry2018divisors} Scott Corry and David Perkinson, \textsl{Divisors and Sandpiles}, American Mathematical Soc. {\bf Vol. 114} (2018).

\bibitem{harris2013algebraic} Joe Harris, \textsl{Algebraic Geometry: A First Course}, Springer Science \& Business Media {\bf Vol. 133} (2013).

\bibitem{eisenbud2013commutative} David Eisenbud, \textsl{Commutative Algebra: With a View Toward Algebraic Geometry}, Springer Science \& Business Media {\bf Vol. 150} (2013).

\bibitem{EAGONNORTHCOTT} J. A. Eagon and D. G. Northcott, \textsl{Ideals Defined by Matrices and a Certain Complex Associated With Them}, Proceedings of the Royal Society of London. Series A, Mathematical and Physical Sciences {\bf 269.1337} (1962), pp. 188--204.

\bibitem{banks2021subcomplexes} Maya Banks and Aleksandra Sobieska, \textsl{Subcomplexes of Certain Free Resolutions}, arXiv preprint arXiv:2112.15137 (2021).

\bibitem{BAKER2007766} Matthew Baker and Serguei Norine, \textsl{Riemann–Roch and Abel–Jacobi Theory on a Finite Graph}, Advances in Mathematics {\bf 215.2} (2007), pp. 766-788.


\bibitem{amini2010riemann} Omid Amini and Madhusudan Manjunath, \textsl{Riemann-Roch for Sub-Lattices of the Root Lattice $A_n$}, The Electronic Journal of Combinatorics {\bf 17.R124} (2010).

\bibitem{klivans2018mathematics} Caroline J. Klivans, \textsl{The Mathematics of Chip-Firing}, Discrete Mathematics and its Applications, CRC Press (2018), ISBN: 9781351801003.


\bibitem{bruns1998cohen} Winfried Bruns and J{\"u}rgen Herzog, \textsl{Cohen-Macaulay Rings}, Cambridge University Press {\bf 39} (1998).

\bibitem{geck2013introduction} Meinolf Geck, \textsl{An Introduction to Algebraic Geometry and Algebraic Groups}, Oxford Graduate Texts in Mathematics, OUP Oxford (2013), ISBN: 9780199676163.

\bibitem{deopurkar2014grobner} Anand Deopurkar, Maksym Fedorchuk and David Swinarski, \textsl{Gr{\"o}bner Techniques and Ribbons}, Albanian Journal of Mathematics {\bf 8.1} (2014).

\bibitem{perkinson2013primer}
David Perkinson, Jacob Perlman and John Wilmes, \textsl{Primer for the Algebraic Geometry of Sandpiles}, Tropical and Non-Archimedean Geometry {\bf 605} (2013), pp. 211--256.
\end{thebibliography}
\end{document}